\documentclass[11pt]{amsart}

\usepackage{amsmath}
\usepackage{amsxtra}
\usepackage{amscd}
\usepackage{amsthm}
\usepackage{amsfonts}
\usepackage{amssymb}
\usepackage{mathtools}
\usepackage{eucal}
\usepackage[all]{xy}
\usepackage{graphicx}
\usepackage{mathrsfs}
\usepackage{euler}
\usepackage{hyperref}
\usepackage{color}
\usepackage{longtable}
\usepackage{float}
\usepackage{caption}
\usepackage{enumerate}
\usepackage{faktor}
\usepackage{bbm}
\usepackage{lineno}
\usepackage{stmaryrd}
\usepackage{rotating}
\usepackage[T1]{fontenc} % improved font encoding
\usepackage[utf8]{inputenc} % for better handling of non-ASCII characters

\usepackage{tikz}
\usepackage{tikz-cd}
\usetikzlibrary{matrix,arrows,decorations.pathmorphing}

\usepackage[bottom=1.25in, top=1.5in, right=1.5in, left=1.5in]{geometry}

\theoremstyle{plain}
  \newtheorem{theorem}{Theorem}
	\newtheorem{thm}[theorem]{Theorem}
  \newtheorem{proposition}[theorem]{Proposition}
	\newtheorem{prop}[theorem]{Proposition}
  \newtheorem{lemma}[theorem]{Lemma}
  \newtheorem{lem}[theorem]{Lemma}
  \newtheorem{cor}[theorem]{Corollary}

\theoremstyle{definition}
	\newtheorem{definition}[theorem]{Definition}
	
  \newtheorem{example}[theorem]{Example}
	
	\newtheorem{conjecture}[theorem]{Conjecture}
	
\theoremstyle{remark}
	
	\newtheorem{remark}[theorem]{Remark}

    {\cleardoublepage\thispagestyle{empty}\null\vfill\begin{center}%
    \bfseries Acknowledgements\end{center}}%

\DeclareMathAlphabet{\mathcal}{OMS}{cmsy}{m}{n}

\newcommand{\one}{\mathbbm{1}}
\newcommand{\A}{\mathbb{A}}

\newcommand{\C}{\mathbb{C}}

\newcommand{\F}{\mathbb{F}}
\newcommand{\G}{\mathbb{G}}

\newcommand{\N}{\mathbb{N}}

\renewcommand{\P}{\mathbb{P}}
\newcommand{\Q}{\mathbb{Q}}
\newcommand{\R}{\mathbb{R}}
\newcommand{\T}{\mathbb{T}}
\newcommand{\Z}{\mathbb{Z}}

\newcommand{\sca}{\mathscr{A}}

\newcommand{\cala}{\mathcal{A}}

\newcommand{\cald}{\mathcal{D}}
\newcommand{\cale}{\mathcal{E}}
\newcommand{\calh}{\mathcal{H}}
\newcommand{\calk}{\mathcal{K}}
\newcommand{\calm}{\mathcal{M}}
\newcommand{\calo}{\mathcal{O}}
\newcommand{\calp}{\mathcal{P}}
\newcommand{\calq}{\mathcal{Q}}
\newcommand{\cals}{\mathcal{S}}
\newcommand{\calw}{\mathcal{W}}
\newcommand{\calu}{\mathcal{U}}
\newcommand{\calx}{\mathcal{X}}

\newcommand{\mfa}{\mathfrak{a}}
\newcommand{\mfb}{\mathfrak{b}}

\newcommand{\mfd}{\mathfrak{d}}
\newcommand{\mff}{\mathfrak{f}}

\newcommand{\mfl}{\mathfrak{l}}

\newcommand{\mfp}{\mathfrak{p}}
\newcommand{\mfq}{\mathfrak{q}}

\newcommand{\mfs}{\mathfrak{S}}
\newcommand{\rmn}{\mathrm{N}}
\newcommand{\rmo}{\mathrm{O}}
\newcommand{\rmt}{\mathrm{T}}

\newcommand{\ab}{{\mathrm{ab}}}
\newcommand{\Ad}{{\mathrm{Ad}}}

\newcommand{\Aut}{{\mathrm{Aut}}}

\renewcommand{\Im}{{\mathrm{Im}}}

\newcommand{\Id}{\mathrm{I}}
\newcommand{\gal}{\mathrm{Gal}}
\newcommand{\Gal}{{\mathrm{Gal}}}

\newcommand{\Hom}{{\mathrm{Hom}}}

\newcommand{\Ind}{{\mathrm{Ind}}}

\newcommand{\Pic}{\mathrm{Pic}}

\newcommand{\std}{\mathrm{std}}

\newcommand{\Reg}{{\mathrm{Reg}}}

\newcommand{\gl}{\mathrm{GL}}
\newcommand{\GL}{{\mathrm{GL}}}
\newcommand{\go}{\mathrm{GO}}
\newcommand{\GO}{{\mathrm{GO}}}
\newcommand{\GSO}{{\mathrm{GSO}}}

\newcommand{\PGL}{{\mathrm{PGL}}}

\newcommand{\SO}{{\mathrm{SO}}}
\newcommand{\SL}{{\mathrm{SL}}}

\newcommand{\gso}{\mathrm{GSO}}

\newcommand{\disc}{{\mathrm{disc}}}
\newcommand{\Div}{{\mathrm{Div}}}

\newcommand{\et}{{\mathrm{\acute{e}t}}}
\newcommand{\End}{{\mathrm{End}}}
\newcommand{\Frob}{{\mathrm{Frob}}}

\newcommand{\sgn}{{\mathrm{sgn}}}
\newcommand{\ord}{{\mathrm{ord}}}
\newcommand{\Tr}{{\mathrm{Tr}}}
\newcommand{\vol}{{\mathrm{vol}}}

\newcommand{\cusp}{\mathrm{Cusp}}

\newcommand{\ks}{\mathrm{KS}}

\newcommand{\hv}{\mathrm{HV}}

\newcommand{\rs}{\mathrm{RS}}

\newcommand{\new}{{\mathrm{new}}}
\newcommand{\opt}{{\mathrm{opt}}}
\newcommand{\stark}{{\mathrm{Stark}}}

\newcommand{\bs}{\backslash}
\newcommand{\wt}{\widetilde}

\newcommand{\longhookrightarrow}{\, \xhookrightarrow{\quad} \,}
\newcommand{\lra}{{\, \longrightarrow \,}}
\newcommand{\iso}{\, \xrightarrow{\widesim{}} \,}

\newcommand{\paren}[1]{\mathopen{}\left(#1\right)\mathclose{}}
\newcommand{\set}[1]{\mathopen{}\left\{#1\right\}\mathclose{}}
\newcommand{\sbrac}[1]{\mathopen{}\left[#1\right]\mathclose{}}
\newcommand{\abrac}[1]{\mathopen{}\left\langle#1\right\rangle\mathclose{}}
\newcommand{\verts}[1]{\mathopen{}\left\lvert#1\right\rvert\mathclose{}}

\newcommand\restr[2]{{% we make the whole thing an ordinary symbol
  \left.\kern-\nulldelimiterspace % automatically resize the bar with \right
  #1 % the function
  \right|_{#2} % this is the delimiter
  }}

\newcommand{\Mid}{\,\middle|\,}
\newcommand{\dotspace}{{}\cdot{}}

\newcommand{\pair}[1]{\abrac{#1}}
\newcommand{\wh}{\widehat}

\newcommand{\widesim}[2][2]{
  \mathrel{\overset{#2}{\scalebox{#1}[1]{$\sim$}}}
}

\renewcommand{\setminus}{-}

\newcommand{\nocontentsline}[3]{}
\newcommand\stoptoc{%
   \let\origcontentsline\addcontentsline
   \let\addcontentsline\nocontentsline
}
\newcommand\resumetoc{%
   \let\addcontentsline\origcontentsline
}

\makeatletter

\renewcommand{\thepart}{\Roman{part}} 

\numberwithin{subsubsection}{subsection}

\renewcommand\part{%
   \if@noskipsec \leavevmode \fi
   \par
   \addvspace{4ex}%
   \@afterindentfalse
   \secdef\@part\@spart}

\def\@part[#1]#2{%
    \ifnum \c@secnumdepth >\m@ne
      \refstepcounter{part}%
      \addcontentsline{toc}{part}{\thepart\hspace{1em}#1}%
    \else
      \addcontentsline{toc}{part}{#1}%
    \fi
    {\parindent \z@ \centering
     \interlinepenalty \@M
     \normalfont
     \ifnum \c@secnumdepth >\m@ne
       \Large\bfseries \partname\nobreakspace\thepart
       \par\nobreak
     \fi
     \huge \bfseries #2%
     \par}%
    \nobreak
    \vskip 3ex
    \@afterheading}
\def\@spart#1{%
    {\parindent \z@ \centering
     \interlinepenalty \@M
     \normalfont
     \huge \bfseries #1\par}%
     \nobreak
     \vskip 3ex
     \@afterheading}

\AtBeginDocument
 {
    \def\@thm#1#2#3{%
      \ifhmode
        \unskip\unskip\par
      \fi
      \normalfont
      \trivlist
      \let\thmheadnl\relax
      \let\thm@swap\@gobble
      \let\thm@indent\indent % indent
      \thm@headfont{\scshape}% heading font small caps
      \thm@notefont{\fontseries\mddefault\upshape}%
      \thm@headpunct{.}% add period after heading
      \thm@headsep 5\p@ plus\p@ minus\p@\relax
      \thm@space@setup
      #1% style overrides
      \@topsep \thm@preskip               % used by thm head
      \@topsepadd \thm@postskip           % used by \@endparenv
      \def\dth@counter{#2}%
      \ifx\@empty\dth@counter
        \def\@tempa{%
          \@oparg{\@begintheorem{#3}{}}[]%
        }%
      \else
        \H@refstepcounter{#2}%
        \hyper@makecurrent{#2}%
        \let\Hy@dth@currentHref\@currentHref
        \AddToHookNext{para/begin}{\MakeLinkTarget*{\Hy@dth@currentHref}}%
        \def\@tempa{%
          \@oparg{\@begintheorem{#3}{\csname the#2\endcsname}}[]%
        }%
      \fi
      \@tempa
    }%
  \dth@everypar={%
    \@minipagefalse \global\@newlistfalse
    \@noparitemfalse
    \if@inlabel
      \global\@inlabelfalse
      \begingroup \setbox\z@\lastbox
       \ifvoid\z@ \kern-\itemindent \fi
      \endgroup
      \unhbox\@labels
    \fi
    \if@nobreak \@nobreakfalse \clubpenalty\@M
    \else \clubpenalty\@clubpenalty \everypar{}%
    \fi
  }%
}

\makeatother

\title[The Harris--Venkatesh conjecture for derived Hecke operators I]
{The Harris--Venkatesh conjecture for derived Hecke operators I:
imaginary dihedral forms}

\author{Robin Zhang}
\address{Department of Mathematics, Columbia University, USA \newline
	\indent Department of Mathematics, Massachusetts Institute of Technology, USA \newline
  \indent Present address: Institut des Hautes \'{E}tudes Scientifiques, France}
\email{rzhang@ihes.fr}

\date{August 11, 2026}

\begin{document}

% \linenumbers

\begin{abstract}
	\normalsize
	The Harris--Venkatesh conjecture relates the action of
  derived Hecke operators on weight-one modular forms to
  algebraic units. We introduce the Harris--Venkatesh
  period, defined on a product of modular curves, and show
  that it specializes to the derived Hecke action on
  weight-one forms. We also define and construct optimal
  forms, two-variable modular forms that explicitly realize
  a theta lifting from the Jacquet--Langlands
  correspondence.
	Using a positive-characteristic multiplicity-one argument
	comparing the Harris--Venkatesh period on newforms and
	optimal forms, we prove the full Harris--Venkatesh
	conjecture for imaginary dihedral weight-one modular forms.
\end{abstract}

\maketitle

\dedicatory{To the memory of Michael Zhao (1995--2018).}

\setcounter{tocdepth}{2}
\tableofcontents

%%%%%%%%%%%%%%%%%%%%%%%%%%%% Introduction %%%%%%%%%%%%%%%%%%%%%%%%%%%%

\part*{Introduction}
\phantomsection
\addcontentsline{toc}{part}{Introduction}
We study the conjecture of
Harris--Venkatesh~\cite{hv},
which frames the general conjectures of
Galatius, Prasanna, and Venkatesh~\cite{prasanna-venkatesh, venkatesh-2,
galatius-venkatesh} on derived Hecke algebras and
motivic cohomology groups in the
coherent cohomology of the Hodge bundle on the modular curve
(cf. \cite{atanasov, horawa, oh} for higher-dimensional
coherent contexts).
In this setting, the Harris--Venkatesh conjecture \cite[Conjecture~3.1]{hv}
gives a modular analogue of the Stark conjecture
by relating units to the predicted action of
derived Hecke operators on weight-$1$ modular forms.

Let $f$ be a modular form of weight $1$ and
level $\Gamma_1(N)$.
Let $\rho_f: \Gal(\overline{\Q}/\Q) \lra \GL(M)$
be its associated Artin representation
realized on a free module $M$ of rank $2$ over $\Z[\chi_{\rho_f}]$
by Deligne--Serre, where $\chi_{\rho_f}$ is the character of $\rho_f$.
Then $\rho_f$ is realized on the Galois group $\Gal(E/\Q)$ of
a finite Galois extension $E$ of $\Q$.
Fix an embedding $E \hookrightarrow \C$.
There is an associated $3$-dimensional complex representation
$\Ad(\rho_f)$, which can be viewed as the action of $\Gal(E/\Q)$
by conjugation through $\rho_f$ on $2 \times 2$ trace-free matrices.
The adjoint representation factors through the projective image of
$\rho_f$, namely the image of $\Gal(E/\Q)$ in
$\GL_2(\C)/\C^\times = \PGL_2(\C) = \SO_3(\C)$.
This projective image is finite and hence is a finite subgroup of
$\SO_3(\C)$, which must be cyclic, $D_{2n}$, $A_4$, $S_4$, or $A_5$.
One can then classify weight $1$ modular forms
by the projective image of $\rho_f$:
Eisenstein series are the modular forms with cyclic image,
dihedral forms are those with image isomorphic to $D_{2n}$,
and the remaining forms are called ``exotic''.
The Stark conjecture for $\Ad(\rho_f)$ is known in the
Eisenstein and dihedral cases. It is also known in the
tetrahedral and octahedral cases, because the corresponding
adjoint characters are rational-valued
\cite[Chapitre~II, Corollaire~7.4]{tate-stark-1984}.
Only the icosahedral case remains open in general.

The Stark conjecture for $\Ad(\rho_f)$ relates the leading term of the
Artin $L$-function $L(\Ad(\rho_f), s)$ to the regulator
of an element of a dual unit group.
Let $R_f := \Z[\chi_{\Ad(\rho_f)}]$ and
$G := \Gal(E/\Q)$.
Let $M_{\Ad} \subset \Ad(\rho_f)$ be the $G$-stable
$R_f$-lattice generated by
\[
  x_\sigma :=
  2\rho_f(\sigma)
  - \Tr\paren{\rho_f(\sigma)}\Id_{2\times2},
  \qquad \sigma \in G.
\]
The \emph{dual unit group} attached to $\Ad(\rho_f)$ is
\[
  \calu\paren{\Ad(\rho_f)}
    := \Hom_{R_f[G]}
      \paren{M_{\Ad}, \calo_E^\times \otimes_\Z R_f};
\]
this is called the ``Stark unit group'' in \cite{hv}.

For an archimedean place $w$ of $E$, let
$c_w \in G$ be the associated complex conjugation and
set $x_w := x_{c_w}$. Evaluation at $x_w$ gives an
injection
\[
  \begin{tikzcd}[row sep = tiny]
    \calu\big(\Ad(\rho_f)\big) \otimes_\Z \Q
    \arrow[r, hookrightarrow, "\operatorname{ev}_w"] &
    \paren{\calo_E^\times \otimes_\Z R_f} \otimes_\Z \Q,
    \\
    u \arrow[r,mapsto] & u(x_w).
  \end{tikzcd}
\]
Its image is contained in the image of the idempotent
\[
  e_{\Ad(\rho_f)}
  :=
  \frac{1}{\verts{G}}
  \sum_{\sigma \in G}
  \chi_{\Ad(\rho_f)}(\sigma)\sigma^{-1}
\]
on $\paren{\calo_E^\times\otimes_\Z R_f} \otimes_\Z \Q$.
Composing evaluation with the logarithm of the absolute
value at $w$ defines the Stark regulator
$\Reg_\R:\calu(\Ad(\rho_f))\otimes_\Z\Q
\longrightarrow \R\otimes_\Z R_f$.

The Stark conjecture for $\Ad(\rho_f)$ predicts that there
is a unique element
$u_\stark \in \calu(\Ad(\rho_f)) \otimes_\Z \Q$
such that
\[
  L'\big(\Ad(\rho_f),0\big)
  =
  \Reg_\R\paren{u_\stark},
\]
compatibly with Galois conjugation of $\Ad(\rho_f)$
under $\Aut(\C)$.
Equivalently, one may choose a unit
$\epsilon \in \calo_E^\times$ and a positive integer $W$
such that
\begin{align*}
  L'\paren{\Ad(\rho_f),0}
  &=
  \frac{1}{W}
  \sum_{\sigma\in G}
  \chi_{\Ad(\rho_f)}(\sigma)
  \log\verts{\epsilon^{\sigma^{-1}}}, \\
  u_\stark(x_w)
  &=
  \frac{1}{W}
  \sum_{\sigma\in G}
  \chi_{\Ad(\rho_f)}(\sigma)
  \epsilon^{\sigma^{-1}}.
\end{align*}

For the mod $p$ setting,
one can define an $\F_p^\times$ regulator map
by taking $x_w$ for non-archimedean places $w$ of $E$ over $p$.
Evaluation at $x_w$ defines an embedding
$\calu\big(\Ad(\rho_f)\big) \lra \paren{\calo_E^\times}^{\Frob_w}
\otimes \Z\sbrac{\chi_{\Ad(\rho_f)}}$,
whose image in $\calo_{E_w}^\times \otimes R_f$ is in
$\Z_p^\times \otimes R_f$.
Thus, the composition of reduction modulo the ideal corresponding to $w$
and evaluation at $x_w$ defines an $\F_p^\times$ regulator map
(called ``reduction of a Stark unit'' in \cite{hv, dhrv}):
\[
	\Reg_{\F_p^\times}: \calu\big(\Ad(\rho_f)\big) \lra
		\F_p^\times \otimes \Z\sbrac{\chi_{\Ad(\rho_f)}}.
\]
The derivative of the Artin $L$-function is then
replaced by the action of the Shimura class.
To construct the Shimura class on coarse modular curves,
we invert $6$ in all coefficient rings below, thereby
removing the possible $2$- and $3$-primary ramification
of the Shimura covering. Set
$\Delta_p:=\F_p^\times\otimes_\Z\Z\sbrac{1/6}$; thus
$\Delta_p$ is the maximal quotient of $\F_p^\times$
whose order is coprime to $6$.

Let $p \geq 5$ be a prime coprime to $N$. For every
prime $\ell \geq 5$ dividing $p - 1$, the
$\ell$-primary quotient of the Shimura covering is
finite \'{e}tale by
Harris--Venkatesh~\cite[Section~3.1]{hv}.
Taking their product gives a finite \'{e}tale covering
over $\Z\sbrac{1/(6p)}$ with Galois group $\Delta_p$
and class
\[
  \mfs_p
  \in
  H^1_{\et}
  \paren{
    X_0(p)_{\Z[1/(6p)]},
    \Delta_p}.
\]
Set $R_p := \Z\sbrac{1/6}/(p-1)$.
The natural map from $\Delta_p$ to
$\calo_{X_0(p)_{R_p}}\otimes_{\Z[1/6]}\Delta_p$
sends this class to
\[
  \mfs_p
  \in
  H^1\paren{X_0(p)_{R_p},\calo}
  \otimes_{\Z\sbrac{1/6}}\Delta_p.
\]
By Serre duality, it defines a functional on weight-$2$
forms, again denoted by
\[
  \mfs_p:
  S_2\paren{\Gamma_0(p),\Z\sbrac{1/6}}
  \longrightarrow
  \Delta_p.
\]

Harris--Venkatesh~\cite[Section~3.2]{hv} then define
a derived Hecke operator $T_{p,N}$ on weight-$1$ cusp
forms of level $N$.
Set $R_{p,N} := \Z\sbrac{1/(6N)}/(p-1)$.
The operator is the composite
\[
\begin{tikzcd}
  H^0\paren{
    X_{\Gamma_1(N),R_{p,N}},
    \omega}
    \arrow[d, "\pi_1^*"]
    \arrow[r, "T_{p,N}"]
  &
  H^1\paren{
    X_{\Gamma_1(N),R_{p,N}},
    \omega}
    \otimes_{\Z\sbrac{1/6}}\Delta_p
  \\
  H^0\paren{
    X_{\Gamma_1(N)\cap\Gamma_0(p),R_{p,N}},
    \omega}
    \arrow[r, "\cup \mfs_p"]
  &
  H^1\paren{
    X_{\Gamma_1(N)\cap\Gamma_0(p),R_{p,N}},
    \omega}
    \otimes_{\Z\sbrac{1/6}}\Delta_p
    \arrow[u, "\pi_{2*}"]
\end{tikzcd}
\]
We apply this operator to cusp forms through the natural
inclusions $\omega(-\cusp)\hookrightarrow\omega$ on the
two modular curves.

For every prime $\ell \geq 5$ dividing $p-1$, let
$\ell^t\Vert p-1$ and fix a discrete logarithm
$\log_\ell:\F_p^\times\twoheadrightarrow
\Z/\ell^t\Z$. It factors uniquely through $\Delta_p$;
we use the same notation for the resulting map
$\log_\ell:\Delta_p\twoheadrightarrow\Z/\ell^t\Z$.
Let $f^* = \sum_n \overline{a_n} q^n$ denote
the dual newform of $f = \sum_n a_n q^n$,
where $\overline{a_n}$ is the complex conjugate
of the $n$-th Fourier coefficient of $f$.

\begin{conjecture}[Harris--Venkatesh conjecture]
	\label{conj:hv}
	Let $f$ be a Hecke new cusp form of weight $1$ and level $N$.
	There is an element $u \in \calu(\Ad(\rho_f))$ and
	a positive integer $m$
	such that for every prime $p \geq 5$ coprime to $N$ and every
  prime $\ell \geq 5$ dividing $p - 1$ and coprime to
  $N$,
	\[
		m \cdot \log_\ell \mfs_p \paren{\Tr_p^{Np}\paren{f(z)f^*(pz)}}
			= \log_\ell \Reg_{\F_p^\times}(u).
	\]
\end{conjecture}

The evaluation of $\mfs_p$ on Eisenstein series
was studied by Mazur~\cite[p. 103]{mazur} and computed by
Merel~\cite{merel2}; see the discussion in \cite[Section 5.2]{hv}.
The first theoretical steps toward the Harris--Venkatesh conjecture
in the dihedral setting
were done by Darmon--Harris--Rotger--Venkatesh \cite{dhrv},
with a verification of the Harris--Venkatesh conjecture
for special cases of dihedral modular forms.
Dihedral forms are
classified by finite characters of $G_K := \gal(\overline{K} / K)$
with $K/\Q$ quadratic:
given a dihedral form $f$,
there is a quadratic number field $K$ and finite character
$\chi$ of $G_K$ such that
$\rho_f = \Ind_{G_K}^{G_\Q}(\chi)$;
conversely, given a character $\chi$ of $G_K$,
there is a newform $f_\chi \in S_1(\Gamma_1(N), \Z[\chi])$
with $q$-expansion $f_\chi(z) = \sum_{n = 1}^\infty a_n q^n$
such that $\sum_n a_n n^{-s} = L(\chi, s)$.
There are two dihedral cases: the imaginary dihedral case
(also called CM dihedral) when $K$ is imaginary quadratic,
and the real dihedral case (also called RM dihedral)
when $K$ is real quadratic.

\begin{theorem}[{Darmon--Harris--Rotger--Venkatesh~\cite[Theorem~1.2]{dhrv}}]
\label{thm:dhrv}
	Let $f$ be a dihedral modular form of weight $1$ and level $N$.
	Let $\chi$ be its associated character of $G_K$ for
	a quadratic number field $K$ of discriminant $D_K$
	and different $\cald_K$.
	\begin{itemize}
		\item If $K$ is imaginary, assume that
			$D_K$ is an odd prime and that $\chi$ is unramified;
		\item if $K$ is real, assume that $D_K$ is odd
			and that $\chi$ has conductor dividing
			$\cald_K$.
	\end{itemize}
	Then the Harris--Venkatesh conjecture for $f$ is true.
\end{theorem}
\begin{remark}
	Darmon--Harris--Rotger--Venkatesh \cite[Section 1.3]{dhrv}
	also checks that both sides of Conjecture~\ref{conj:hv}
	vanish when $K$ is imaginary quadratic (resp. real quadratic)
	and $p$ splits in $K$ (resp. $p$ is inert in $K$),
	so one can assume that $p$ is inert (resp. splits)
	for the purpose of proving the Harris--Venkatesh conjecture.
\end{remark}

Let $c_K$ denote the nontrivial element of
$\Gal(K/\Q)$. For a character $\chi$ of $G_K$, write
$\chi^{c_K}$ for its conjugate by any lift of $c_K$ to
$G_\Q$; this does not depend on the lift. We call
$\chi^{1-c_K}:=\chi/\chi^{c_K}$ the \emph{antinorm}
of $\chi$.
The methods of Darmon--Harris--Rotger--Venkatesh \cite{dhrv}
cannot be directly generalized to ramified characters.
A recent approach by
Lecouturier~\cite{lecouturier-hv}
replaces some of the theta lifting arguments of \cite{dhrv}
that rely on $\chi$ being unramified by
the central $L$-value formulas of
Waldspurger \cite{waldspurger} and Ichino \cite{ichino}
to prove some ramified cases of
a weak form of the Harris--Venkatesh conjecture.
In particular, \cite[Theorem~1.3]{lecouturier-hv} shows that
if $D_K$ is odd and the antinorm
$\xi := \chi^{1 - c_K}$ is unramified, then there are a
positive integer $m$, an element
$u \in \calu\paren{\Ad(\rho_f)}$, and a sign
$s_{p,N} \in \{\pm 1\}$ depending on $p$ and $N$ such
that
\[
  s_{p,N}m\log_\ell\mfs_p
  \paren{\Tr_p^{Np}\paren{f_\chi(z)f_\chi^*(pz)}}
  =
  \log_\ell\Reg_{\F_p^\times}\paren{u}
\]
for every pair of primes $(\ell,p)$ such that
$\ell^2\nmid p-1$. The sign occurs because the
special-value formulas involve a square, while the
square-freeness of $p-1$ is used to apply Mazur's
congruence for Eisenstein series. Lecouturier
\cite[Section~1]{lecouturier-hv} remarks that the
method might extend to ramified $\xi$ using a more
complicated special-value formula, but that it is unclear
how to remove $s_{p, N}$ or treat non-squarefree $p - 1$ in
general.

The main result of this article is
the resolution of Conjecture~\ref{conj:hv}
for all weight-one forms in the imaginary dihedral case,
including those that arise from ramified $\chi$,
ramified $\xi$, and number fields of
even discriminant $D_K$.
\begin{thm}
	\label{thm:hv}
	\label{THM:HV}
	Let $f$ be an imaginary dihedral modular form
	of weight $1$.
	Then the Harris--Venkatesh conjecture for $f$ is true.
\end{thm}

\begin{remark}
	\label{rem:hv-stark}
	The second paper in this series
	\cite{zhang-hvs} explores the
	connection between the Harris--Venkatesh conjecture for $f$ and the
	Stark conjecture for $\Ad(\rho_f)$ further; in particular, it
	explicates how the unit $u$ relates to the Stark unit $u_\stark$.
  It also provides a simpler formulation of Conjecture \ref{conj:hv}
  without the auxiliary prime $\ell \geq 5$
  and the discrete logarithm $\log_\ell$;
  the discrete logarithm is defined up to multiplication by a scalar
  and both sides of Conjecture \ref{conj:hv} depend linearly
  on the choice of the discrete logarithm,
  so it is only useful for removing considerations
  modulo $2$ and $3$.
	The third paper in this series \cite{zhang-rs} uses the explicit
	constructions of this article to refine the Harris--Venkatesh conjecture
	and explicitly describe the integer $m$.
	We expect the methods developed here to also be applicable to the
	real dihedral case of Conjecture~\ref{conj:hv},
	which is the subject of the forthcoming work \cite{zhang-real}.
\end{remark}

%%%%%%%%%%%%%%%%%%%%%%%%%%%%%%%

\section*{Ingredients of the main theorem}

For a weight-$1$ newform $f$ with dual newform $f^*$,
let $\pi_f$ (resp. $\pi_{f^*}$)
denote the automorphic representation of
$\GL_2(\A)$ generated by $f$ (resp. $f^*$)
under the action of Hecke operators.
To prove Theorem \ref{thm:hv}, we introduce two new objects
and establish two general comparison theorems:
\begin{enumerate}[\bfseries \quad\enspace(A)]
  \item The Harris--Venkatesh period $\calp_\hv$ on
    $\pi_f \otimes \pi_{f^*}$ with values in
    \[
      \Z\sbrac{
        \chi_{\Ad(\rho_f)},
        \zeta_N,
        1/(6N)}
      \otimes_{\Z\sbrac{1/6}}\Delta_p;
    \]
  \item An explicit uniqueness, up to scaling, theorem
    for pairings from $\pi_f \otimes \pi_{f^*}$ to
    \[
      \Z\sbrac{
        \chi_{\Ad(\rho_f)},
        1/(6N)}
    \]
    modulo $(p - 1)$;
  \item The optimal forms $f^\opt$ and $\varphi^\opt$,
    a two-variable modular form of finite full level
    and its adelization on
    $\pi_f\otimes\pi_{f^*}$;
  \item A comparison between $f^\opt(z,pz)$ and the
    $p$-th explicit Jacquet--Langlands transfer of
    $[\one]\otimes[\xi]$.
\end{enumerate}
Strictly speaking, in \textbf{(A)} and \textbf{(B)}
the period and comparison pairing are defined on
suitable integral models of the $\Sigma$-components
$(\bigotimes_{q\in\Sigma}\pi_{f,q})
\otimes
(\bigotimes_{q\in\Sigma}\pi_{f^*,q})$, where
$\Sigma := \set{q \mid N}$. For simplicity, we suppress
these integral models and the $\Sigma$-subscripts in
the rest of the introduction; see
Sections~\ref{sec:comparison}
and~\ref{sec:proof-hv-reform}.

Points \textbf{(A)} and \textbf{(B)} transfer the
Harris--Venkatesh conjecture from the newform tensor
$f\otimes f^*$ to any vector in
$\pi_f\otimes\pi_{f^*}$ with nonzero canonical period.
Points \textbf{(C)} and \textbf{(D)} construct a
distinguished vector
$\varphi^\opt\in\pi_f\otimes\pi_{f^*}$ whose theta
realization is directly related to the elliptic units of
Darmon--Harris--Rotger--Venkatesh~\cite{dhrv}. Unlike
the theta-lifting argument of \cite{dhrv} and the
special-value argument of \cite{lecouturier-hv}, this
construction allows $\chi$ to be arbitrarily ramified.

The latter two ingredients are the most technical. Their
construction uses Atkin--Lehner theory, $q$-expansion
principles, and global theta lifting for
$G=\GL_2\times_{\G_m}\GO(V)$. The resulting newvector
criteria and Whittaker and Kirillov formulas may also be
of independent interest.

\textbf{(A)}
The first main ingredient is the \emph{Harris--Venkatesh period}.
Let $\Sigma$ be the set of primes dividing $N$
and consider the projective system $X_\Sigma = \varprojlim_m X(N^m)$
of modular curves unramified outside of $\Sigma$ with
cuspidal divisor $C_\Sigma$.
The Harris--Venkatesh period is given on two copies
of the space of weight-$1$ cusp forms unramified outside of
$\Sigma$:
\[
  \calp_\hv:
  H^0\paren{X_\Sigma,\omega(-C_\Sigma)}
    \otimes
  H^0\paren{X_\Sigma,\omega(-C_\Sigma)}
  \longrightarrow
  \Z\sbrac{1/(6N)}
    \otimes_{\Z\sbrac{1/6}}\Delta_p.
\]
Moreover, this pairing is invariant under the action of
$\prod_{q \mid N}\GL_2(\Q_q)$ and satisfies
\begin{align*}
  \mfs_p\paren{\Tr_p^{Np}\paren{f(z)f^*(pz)}}
    &= \sbrac{\Gamma(1) : \Gamma_0(N)}
      \calp_\hv\paren{f \otimes f^*}.
\end{align*}

\textbf{(B)}
Using a modulo-$(p - 1)$ multiplicity-one argument
to compare the Harris--Venkatesh period $\calp_\hv$
against the canonical pairing
$\calp: \pi_f \otimes \pi_{f^*} \rightarrow \Q(\chi_{\Ad(\rho_f)})$
on the dual pair $(\pi_f, \pi_{f^*})$
of cuspidal automorphic representations of $\gl_2(\A)$,
we show that the Harris--Venkatesh conjecture for the
form $f$ is equivalent to the Harris--Venkatesh conjecture for
any element of $\pi_f \otimes \pi_{f^*}$.
\begin{theorem}
\label{thm:hv-reform}
  Let $f$ be a newform of level $\Gamma_1(N)$ and
  weight $1$.
  Set $\calo_{f, N} := \calo(f)\sbrac{1/N}$,
  let $\Pi_{f, \Sigma, \calo_{f,N}}$ be the integral
  model of the $\Sigma$-component of
  $\pi_f \otimes \pi_{f^*}$ defined in
  Section~\ref{sec:proof-hv-reform}, and set
  \[
    \calu_f\paren{\Ad(\rho_f)}
    :=
    \calu\paren{\Ad(\rho_f)}
      \otimes_{R_f}\calo(f).
  \]
  The Harris--Venkatesh conjecture for $f$ is
  equivalent to the existence of
  $(\varphi,u,m) \in
      \Pi_{f, \Sigma, \calo_{f,N}}
        \times
      \calu_f\paren{\Ad(\rho_f)}
        \times \N$
  such that $\calp(\varphi) \neq 0$ and
  \[
    m \cdot \log_\ell\calp_\hv\paren{f_\varphi}
    =
    \log_\ell\Reg_{\F_p^\times}\paren{u}
  \]
  for every prime $p \geq 5$ coprime to $N$ and every
  prime $\ell \geq 5$ dividing $p - 1$ and coprime to $N$.

  Moreover, if the Harris--Venkatesh conjecture
  for $f$ holds, then the above property holds for every
  $\varphi \in \Pi_{f, \Sigma, \calo_{f,N}}$ such that
  $\calp(\varphi) \neq 0$.
\end{theorem}
 
\textbf{(C)}
For imaginary dihedral $f$,
the adjoint representation $\Ad(\rho_f)$
decomposes as $\eta \oplus \Ind_{G_K}^{G_\Q}(\xi)$,
where $\eta$ is the quadratic character of
$G_\Q$ for the imaginary quadratic extension $K/\Q$
and the ``antinorm'' $\xi := \chi^{1 - c_K}$
is a ring class character of conductor $c$.
In this setting, we introduce the second main ingredient:
the two-variable \emph{optimal form} $f^\opt(z_1, z_2)$
of finite full level, generated by $f(z_1) f^*(z_2)$
under the action of Hecke operators.
The optimal form is defined explicitly in terms of
theta series using characteristic functions of optimal orders.

\textbf{(D)}
We define a map $[\dotspace]$ from Hecke characters of
$K$ to automorphic forms on a definite quaternion
algebra, extending the construction of
Darmon--Harris--Rotger--Venkatesh
\cite[Section~2.2]{dhrv}. Enlarge the splitting field
$E$ so that it contains the ring class field $H_c$ of
$\xi$. For an auxiliary prime $p$ inert in $K$ and a
prime $\mfp$ of $E$ above $p$, reduction modulo
$\mfp_c := \mfp \cap H_c$ orients a $\xi$-optimal maximal
order $\calo_\mfp$ in the definite quaternion algebra
of reduced discriminant $p$. The oriented optimal
embedding $\calo_c\longhookrightarrow\calo_\mfp$ gives
vectors $[\one]_\mfp$ and $[\xi]_\mfp$ and a standard
Hecke-series theta lift $\Theta_\mfp^{\std}$. After
scaling by a comparison constant $c_{\chi,\mfp}$, set
$\Theta_\mfp^f := c_{\chi,\mfp}\Theta_\mfp^{\std}$.
Our comparison theorem gives
$f^\opt(z, pz) = \Theta_\mfp^f\paren{
[\one]_\mfp \otimes [\xi]_\mfp}$. In the unramified,
odd-discriminant setting, the same standard kernel
recovers the full two-character trace identity of
\cite[Theorem~2.2]{dhrv}, while this formula is its
dual-pair specialization. For arbitrary conductor,
the transfer theorem of Part~I carries the
Harris--Venkatesh period of the optimal form to the
newform tensor.

\begin{theorem}
\label{thm:opt-unique}
  Let $f$ be an imaginary dihedral modular form of
  weight $1$ and level $N$. There is a two-variable
  modular form $f^\opt(z_1, z_2)$ of finite full level
  such that, for every prime $p \geq 5$ coprime to
  $N$ and inert in $K$, and every prime $\mfp$ of $E$
  above $p$,
  \[
    f^\opt(z, pz)
    =
    \Theta_\mfp^f
    \paren{
      \sbrac{\one}_\mfp
      \otimes
      \sbrac{\xi}_\mfp}.
  \]
  Moreover, let $g$ be another two-variable modular
  form of finite full level satisfying these identities
  for all such primes $p$ and $\mfp$.
  If $N'$ is a positive integer for which both $f^\opt$
  and $g$ descend to $X(N') \times X(N')$, and if
  $X(N')^\circ$ denotes the connected component of
  $X(N')_\C$ containing the standard cusp, then
  $g = f^\opt$ on
  $X(N')^\circ \times X(N')^\circ$.
\end{theorem}

In Section~\ref{sec:proof-opt}, we combine the theta
identity of Theorem~\ref{thm:opt-unique} with the higher
Eisenstein and elliptic-unit calculations of
Darmon--Harris--Rotger--Venkatesh~\cite{dhrv}. The
comparison scalar $c_{\chi,\mfp}$ depends on the oriented
prime $\mfp$, but its $\mfp$-dependence agrees with that
of the regulator of a single element
$u_f^\opt\in\calu(\Ad(\rho_f))
\otimes_{R_f}\Q(\chi)$. Thus $u_f^\opt$ is independent
of the auxiliary primes.

For a prime $\mfp$ above $p$, let $\Reg_\mfp$ denote
the $\F_p^\times$-regulator defined in the
introduction using evaluation at $x_\mfp$ and reduction
modulo $\mfp$. Define
\[
  m(\xi)
  :=
  \begin{cases}
    v
      & \text{if $\verts{\Im(\xi)}$ is a power
        of a prime $v$}, \\
    1
      & \text{otherwise}.
  \end{cases}
\]

\begin{theorem}
\label{thm:opt}
\label{THM:OPT}
  Let $f$ be an imaginary dihedral modular form of
  weight $1$ and level $N$, and let $f^\opt$ be its
  optimal form. There is an element
  \[
    u_f^\opt
    \in
    \calu\paren{\Ad(\rho_f)}
      \otimes_{R_f}\Q(\chi),
  \]
  independent of $p$, $\mfp$, and $\ell$,
  such that, for every prime $p \geq 5$ coprime to
  $N$ and inert in $K$, every prime $\mfp$ of $E$
  above $p$, and every prime $\ell \geq 5$ dividing
  $p - 1$ and coprime to $N$,
  \[
    \log_\ell\mfs_p\paren{f^\opt(z, pz)}
    =
    -\frac{\sbrac{H_c : K}}{12m(\xi)}
    \log_\ell\Reg_\mfp\paren{u_f^\opt}.
  \]
\end{theorem}

%%%%%%%%%%%%%%%%%%%%%%%%%%%%%%%

\section*{Outline}
The contents of this paper are divided into two parts.

Part~\ref{part:hv} develops the theory of
the Harris--Venkatesh period $\calp_\hv$ and proves the equivalence
between the Harris--Venkatesh conjecture for $f$ and the
Harris--Venkatesh conjecture for any
$\varphi \in \pi_f \otimes \pi_{f^*}$.
Section~\ref{sec:hv} recalls necessary background results on
modular curves and automorphic forms on $\gl_2$.
Section~\ref{sec:hv-period} defines
the Harris--Venkatesh period $\calp_\hv$
on the cohomology of modular curves
using the Shimura class.
Section~\ref{sec:comparison}
gives a multiplicity-one argument
comparing pairings on
the tensor product $\pi \otimes \wt{\pi}$
for $\pi$ an irreducible representation
of $\gl_2$ over a $p$-adic field.
Section~\ref{sec:proof-hv-reform} proves Theorem~\ref{thm:hv-reform}
using the comparison between the Harris--Venkatesh pairing
and canonical pairing on $\pi_f \otimes \pi_{f^*}$
on general forms and newforms.

Part~\ref{part:opt} develops the theory of optimal forms
and proves the Harris--Venkatesh conjecture for optimal forms.
Section~\ref{sec:global-theta} describes global theta liftings
for $G = \gl_2 \times_{\G_m} \go(V)$, especially
for $V$ a quadratic number field or a quaternion algebra.
Section \ref{sec:opt} defines the optimal form as
an automorphic form $\varphi^\opt$ in
the representation $\pi_f \otimes \pi_{f^*}$ using theta series
and as a modular form $f^\opt$ of finite full level.
Section~\ref{sec:opt-unique} constructs the standard
oriented theta lift $\Theta_\mfp^{\mathrm{std}}$, defines
its $f$-normalization $\Theta_\mfp^f$, and proves
Theorem~\ref{thm:opt-unique}. In particular,
$f^\opt(z,pz)$ is the value of $\Theta_\mfp^f$ at
$[\one]_\mfp\otimes[\xi]_\mfp$.
Section~\ref{sec:proof-opt} makes the oriented-order
description precise. It recovers the trace identity of
\cite[Theorem~2.2]{dhrv}, constructs an auxiliary-free
elliptic unit and its fixed adjoint lift, and proves
Theorem~\ref{thm:opt}.
Section~\ref{sec:proof-hv}
combines this result with the comparison theorem of
Part~\ref{part:hv} and treats split auxiliary primes
by simultaneous vanishing.

\section*{AI Usage}

The paper was proofread and copyedited with assistance from GPT 5.6 Pro.

\section*{Acknowledgements}

I am deeply thankful to
Michael Harris, who supervised my thesis of which
this article is roughly the first half,
for suggesting this area of research,
for sharing his deep and broad insight,
and for his continued guidance throughout my doctoral studies.
I am grateful to Henri Darmon for discussions about
quaternion algebras and periods on Shimura curves,
to Wee-Teck Gan for discussions regarding
explicit local theta liftings,
to Dick Gross for observations
about Stark units and the Stark conjectures,
to Will Sawin for a helpful observation
about the compatibility of pullbacks of the
Harris--Venkatesh period,
and to Emmanuel Lecouturier, Lo\"{i}c Merel,
and Akshay Venkatesh
for illuminating conversations.
I thank the anonymous referee for a careful reading and
for many helpful comments.

This material is based on work
supported by the National Science Foundation
under Grants DGE-1644869, DGE-2036197, and DMS-2303280.

%%%%%%%%%%%%%%%%%%%%%%%%%%%%%%%%%%%%%%%%%%%%%%%%%%%%%%%%%%%%%%%%%%%%%%%

\numberwithin{equation}{section}
\numberwithin{theorem}{section}

\newpage

\part{Harris--Venkatesh periods}
\label{part:hv}

%%%%%%%%%%%%%%%%%%%%%%%%%%%%%%%

\section{Automorphic forms on \texorpdfstring{$\GL_2$}{GL(2)}}
\label{sec:hv}

This section collects the background on automorphic forms,
modular forms, and modular curves used later. For
automorphic forms on $\GL_2$, we follow
\cite{jacquet-langlands,jacquet,casselman,gelbart,bump,
zhang-asian}; for modular forms and modular curves, we use
\cite{deligne-rapoport,katz-mazur,shimura}. No result in
this section is new. Readers familiar with these
conventions may proceed directly to
Section~\ref{sec:hv-period}.

\subsection{Automorphic forms}
\label{sec:automorphic-forms}

We start with some notation:
\[
	\wh \Z := \varprojlim _n \Z/n\Z=\prod_p\Z_p,
		\qquad \wh \Q := \Q\otimes_\Z \wh\Z.
\]
Let $\A := \wh \Q \times \R$ denote the ring of ad\`{e}les of $\Q$.
For any abelian group $M$, let $\wh M := M \otimes_\Z \wh \Z$.
For any $\Q$-vector space $V$, let $V_\A := V \otimes_\Q \A$.
We use superscripts to remove specified places and
subscripts to restrict to specified places.
For example, $\A^\infty = \wh \Q$
and $\A_\infty = \Q_\infty = \R$.

Let $\psi: \Q \bs \A \lra \C^\times$ be the standard additive character:
\[
	\psi (x) = \prod_{p \leq \infty} \psi_p(x_p)
		= e^{2\pi i x_\infty}
			e^{-2\pi i (x^\infty \, \mathrm{mod} \, \wh \Z)},
\]
where we have used the identity $\Q / \Z \iso \wh \Q / \wh \Z$ by 
the Chinese remainder theorem. 

For an algebraic group $G$ over $\Q$, we denote its center by
$Z_G$, denote the quotient $Z_G \bs G$ by $PG$,
and denote the set $G(\Q)\bs G(\A)$ by $[G]$. For example,
$[\PGL_2] = Z(\A) \bs \GL_2(\Q) \bs \GL_2(\A)$.
 
For a reductive group $G$ over $\Q$, let $\cala([G])$ denote its space
of automorphic functions. These are smooth functions with 
some growth conditions on $G(\Q)\bs G(\A)$. Let $\cala_0([G])$ denote
the subspace of cusp forms: the space of automorphic functions that
vanish at cusps. 

\subsubsection*{Global Whittaker models and Kirillov models}
Let $N \subset \GL_2$ denote the standard upper-triangular
unipotent subgroup, and write
$[N] := N(\Q) \bs N(\A)$.
We will mainly focus on $\GL_2$,
so $\cala_0([\GL_2])$ denotes the
space of smooth functions $\varphi: \GL_2(\A)\lra \C$ invariant
under left translation by $\GL_2(\Q)$ that vanish at cusps:
\[
	\int_{[N]} \varphi(ng) dn = 0,
\]
where $dn$ is a Haar measure on $N(\A)$ such that the volume of
$N(\Q)\bs N(\A)$ is $1$. 
For such a function, we can define its Whittaker function:
\[
	W_\varphi(g) := \int_{[N]} \varphi(ng) \psi_N^{-1}(n) dn,
\]
where $\psi_N$ is the character on $N(\A)$ via the canonical
isomorphism
\begin{align*}
  n: \A &\iso N(\A) \\
  x &\longmapsto \begin{pmatrix} 1 & x \\ 0 & 1 \end{pmatrix}.
\end{align*}
Then we have the Fourier expansion 
\[
	\varphi (g) = \sum_{a \in \Q^\times} W_\varphi
		\paren{\begin{pmatrix} a & \\ & 1 \end{pmatrix}g}.
\]
Therefore the Fourier transform induces an embedding of representations of
$\GL_2(\A)$:
\[
	W: \cala_0\paren{\sbrac{\GL_2}} \lra \calw(\psi) :=
		\Ind_{N(\A)}^{\GL_2(\A)}\paren{\psi_N}.
\]
 
By a cuspidal automorphic representation,
we mean a subrepresentation $\pi \subset \cala_0([\GL_2])$.
We can embed $\pi$ into $\calw(\psi)$, 
\[
	\pi \longhookrightarrow \calw(\psi).
\]
More precisely, there is a unique embedding
of $\pi$ into $\calw(\psi)$; its image is denoted by
its Whittaker model $\calw(\pi, \psi)$
and has a decomposition
\[
	\calw(\pi, \psi) = \bigotimes_p \calw(\pi_p, \psi_p).
\]
Each element in $\calw(\pi, \psi)$ can be written as a finite linear
combination of pure tensors $\bigotimes W_p$
such that for all but finitely many $p$,
$W_p \in \calw(\pi_p, \psi_p)$ is the
normalized spherical element in the sense that
$W_p$ is invariant under $\GL_2(\Z_p)$ and $W_p(e) = 1$.

Let $\omega$ be the central character of $\pi$. Then define the
Kirillov representation $\calk(\omega, \psi)$ of
the Borel subgroup $B(\A)$ on
$C^\infty(\A^\times)$ by the formula
\[
	\begin{pmatrix} a & b \\ 0 & d \end{pmatrix} f(x)
		= \omega(d) \psi\paren{\frac{bx}{d}}
			f\paren{\frac{ax}{d}}.
\]
Following Jacquet--Langlands~\cite{jacquet-langlands},
the restriction 
\begin{align*}
  \calw(\pi, \psi) &\lra \calk(\omega, \psi), \\
    W &\longmapsto \kappa_W(x) :=
      W\paren{\begin{psmallmatrix} x & 0 \\ 0 & 1 \end{psmallmatrix}},
\end{align*}
is injective. Let $\calk(\pi, \psi)$ denote the image
of this map.

\subsubsection*{Newforms}
Each irreducible cuspidal representation $\pi$ has associated data
(weight, level, central character, and newforms), defined
as follows (cf. \cite[Section~2.3]{zhang-asian}).
\begin{enumerate}
  \item The \textit{weight} $w$ of $\pi$ is the least
    non-negative integer for which $\pi_\infty$
    contains a non-zero vector $v$ satisfying
    \[
      \begin{pmatrix}
        \cos \theta & \sin \theta \\
        -\sin \theta & \cos \theta
      \end{pmatrix}
      v
      =
      e^{iw\theta}v
      \qquad
      \text{for every $\theta \in \R$.}
    \]
    For the corresponding discrete-series
    representation of $\GL_2(\R)$, the remaining
    $\SO_2(\R)$-types have weights
    $\pm(w + 2k)$ for $k \in \Z_{\geq 0}$; see
    \cite[Theorem~2.5.4(ii)]{bump}.
  \item The \textit{level} $N$ of $\pi$ is the least
    positive integer for which
    $(\pi^\infty)^{U_1(N)}$ is non-zero, where
    \[
      U_1(N)
      :=
      \set{
        \begin{pmatrix}a & b \\ c & d\end{pmatrix}
        \in \GL_2(\wh\Z)
        \Mid
        (c,d) \equiv (0,1)
        \pmod{N\wh\Z}}.
    \]
    The space $(\pi^\infty)^{U_1(N)}$ is then
    one-dimensional.
  \item The \textit{central character} $\omega$ of
    $\pi$ is the character of
    $[Z] \iso \Q^\times \bs \A^\times$ acting on
    $\pi$.
  \item For every finite prime $p$, let
    $W_p^\new \in \calw(\pi_p,\psi_p)$ be the
    normalized local new vector, so that
    $W_p^\new(e) = 1$. The \textit{global new vector}
    $\varphi^\new \in \pi$ is the vector whose
    Whittaker function factors as
    \[
      W_{\varphi^\new}
      =
      W_w
      \otimes
      \bigotimes_{p < \infty}W_p^\new,
    \]
    where $W_w$ is the standard holomorphic
    lowest-weight Whittaker function defined below.
\end{enumerate}

For a positive integer $k$, let $W_k$ denote the
standard holomorphic lowest-weight Whittaker function
on $\GL_2(\R)$ of weight $k$. For $y \in \R^\times$, put
$a(y) := \begin{psmallmatrix}y & 0 \\ 0 & 1\end{psmallmatrix}$.
With our choice of additive character, $W_k$ is
supported on $\GL_2(\R)_+$ and satisfies
\begin{equation}
\label{eq:Whittaker-GL2R}
  W_k\paren{a(y)}
  =
  \begin{cases}
    y^{k/2}e^{-2\pi y} & \text{if $y > 0$,} \\
    0 & \text{if $y < 0$.}
  \end{cases}
\end{equation}
This $W_k$ is the classical Whittaker function associated to
holomorphic discrete series
(originally introduced by Whittaker \cite{whittaker},
cf. $W_{\frac{k}{2}, s - \frac{1}{2}}$ in \cite[Section 2.8]{bump}).
Another example is the Whittaker model of
an unramified principal series $\pi(\chi_1, \chi_2)$.
Again we need only consider the value of the
(spherical) Whittaker function at
$a(p^n)=\begin{psmallmatrix}p^n & 0\\ 0 & 1\end{psmallmatrix}$,
which was calculated independently
by Shintani \cite{shintani}, Kato \cite{kato},
Kazhdan (unpublished), and Casselman--Shalika \cite{casselman-shalika}
(see \cite[Theorem 4.6.5]{bump}):

\[
	W\paren{a\paren{p^n}} = p^{\frac{-n}{2}}
		\frac{\alpha_1^{n + 1} - \alpha_2^{n + 1}}{\alpha_1 - \alpha_2},
\]
where $\alpha_1 = \chi_1(p)$ and $\alpha_2 = \chi_2(p)$.

\subsection{Modular forms in finite level}
We recall some background on the theory
of modular forms in finite level, much of
which can be found in the classical
texts of Shimura~\cite[Chapter~6]{shimura},
Deligne--Rapoport \cite[Chapters~IV, VII]{deligne-rapoport},
and Katz--Mazur \cite[Chapters~3-4]{katz-mazur}
(also cf. \cite[Section~1]{katz-1973}
and \cite[Sections~7-9]{diamond-im}).

Let
\[
  \calh := \set{z \in \C \Mid \operatorname{Im}(z) > 0}
\]
denote the complex upper half-plane.
For each positive integer $n$, we have a
modular curve $Y(n)$ over $\Z[1/n]$
parametrizing isomorphism classes of pairs
$(E,\phi)$, where
\[
  \phi: (\Z/n\Z)^2 \iso E[n]
\]
is a full level-$n$ structure.
Over the complex numbers, there is a uniformization
\[
  Y(n)(\C)
  =
  \SL_2(\Z) \bs
  \paren{\calh \times \GL_2(\Z/n\Z)},
\]
where the action of $\gamma \in \SL_2(\Z)$
on $\calh \times \GL_2(\Z/n\Z)$ is given by
\[
  \gamma \cdot (z,g)
  :=
  \paren{\gamma z,(\gamma \bmod n)g}.
\]
If
\[
  g =
  \begin{psmallmatrix}
    a & b\\
    c & d
  \end{psmallmatrix}
  \in \GL_2(\Z/n\Z),
\]
and $a,b,c,d$ are integral lifts of its entries, then the
class of $(z,g)$ corresponds to the pair $(E_z,\phi_g)$, where
\begin{align*}
  E_z
  &:= \C/(\Z z+\Z),\\
  \phi_g(m_1,m_2)
  &:=
  \frac{m_1(dz-b)+m_2(-cz+a)}{n}
  \bmod (\Z z+\Z)
  \in E_z[n].
\end{align*}
The set of connected components is isomorphic to
$(\Z/n\Z)^\times$, and the decomposition is given by
\[
	Y(n)(\C) = \Gamma(n) \bs \calh\times
		\set{\begin{pmatrix} a & 0 \\ 0 & 1 \end{pmatrix}
		\Mid a \in (\Z/n\Z)^\times}.
\]
In fact, each of these connected components is defined over
$\Z[1/n, \zeta_n]$, where $\zeta_n$ is a primitive $n$-th root of
unity,
and the corresponding component
for each $a \in (\Z/n\Z)^\times$
parametrizes pairs $(E, \phi)$ such that the
Weil pairing $\pair{\phi(1, 0), \phi(0, 1)}$ is equal to $\zeta_n^a$.

When $n \geq 3$, there is a universal elliptic curve
\[
  \pi: \cale(n) \lra Y(n),
\]
which can be constructed analytically as
\[
  \cale(n)
  :=
  \paren{\Z^2 \rtimes \SL_2(\Z)}
  \bs
  \paren{\calh \times \C \times \GL_2(\Z/n\Z)}.
\]
Here
$(m_1,m_2,\gamma) \in \Z^2 \rtimes \SL_2(\Z)$
acts by
\[
  (z,u,g)
  \longmapsto
  \paren{
    \gamma z,\,
    j(\gamma,z)^{-1}(u+m_1z+m_2),\,
    (\gamma \bmod n)g
  }.
\]
For
$\gamma =
  \begin{psmallmatrix}
    a & b\\
    c & d
  \end{psmallmatrix}
  \in \GL_2(\R)_+$,
we use the automorphy factor
\[
  j(\gamma,z)
  :=
  \verts{\det\gamma}^{-\frac{1}{2}}(cz+d).
\]

Let
\[
  \omega
  :=
  \pi_*\Omega^1_{\cale(n)/Y(n)}
\]
denote the Hodge bundle, namely the pushforward of the
sheaf of relative invariant differentials.
For any integer $k$, the space
$H^0(Y(n),\omega^k)$
is the space of weakly holomorphic modular forms of weight $k$.
Over the complex numbers, every such form can be written as a
function $f(z, g) du^k$ on $\calh \times \GL_2(\Z/n\Z)$
that is holomorphic in $z$ and invariant under the action by every
$\gamma \in \SL_2(\Z)$:
\[
	f(\gamma z, \gamma g) d(\gamma u)^k = f(z, g) du^k.
\]
Notice that for
$\gamma = \begin{psmallmatrix}a & b \\ c & d \end{psmallmatrix}$, 
\[
	d(\gamma u) = j(\gamma, z)^{-1} du,
\]
so then for all $\gamma \in \SL_2(\Z)$,
\[
	f(\gamma z, \gamma g) = f(z, g) j(\gamma, z)^k.
\]

The modular curve $Y(n)$ admits a compactification $X(n)$
obtained by adjoining the cuspidal divisor $C(n)$, and
the universal elliptic curve extends to a generalized elliptic
curve over $X(n)$.
The Hodge bundle $\omega$ extends canonically to $X(n)$.
We call
\[
  H^0(X(n),\omega^k)
  \qquad\text{and}\qquad
  H^0\paren{X(n),\omega^k(-C(n))}
\]
the spaces of modular forms and cusp forms of weight $k$,
respectively.

Over the complex numbers, we have
\[
  X(n)(\C)
  =
  \SL_2(\Z) \bs
  \paren{\wh\calh \times \GL_2(\Z/n\Z)},
\]
where $\wh\calh := \calh \cup \P^1(\Q)$
is the extended upper half-plane.
The cuspidal divisor is described as
\[
  C(n)
  =
  \SL_2(\Z) \bs
  \paren{\P^1(\Q) \times \GL_2(\Z/n\Z)}.
\]
At each cusp $c$, there is a well-defined holomorphic coordinate
$q_c$ on $X(n)$.
If $c$ is represented by $(a, g) \in \P^1(\Q) \times \GL_2(\Z/n\Z)$
with stabilizer $N_{a, g}$ in $\SL_2(\Z)$,
then $q_c$ is represented by a generator of holomorphic functions
invariant under $N_{a, g}$.
For example, if $c$ is represented by $\infty \times e
\in \P^1 \times \GL_2(\Z/n\Z)$, which has stabilizer $N(n\Z)$,
then we take $q_c = e^{\frac{2 \pi i z}{n}}$.
A modular form $f(z) du^k$, or rather $f(z)$, is then a 
weakly holomorphic modular form with Taylor expansion at each cusp, 
\[
	f(z) = \sum_{i \geq 0} a_{c, i} q_c^i.
\]
Such a form is a cusp form if it vanishes at cusps,
i.e. $a_{c, 0} = 0$.

For a $\Z[1/n]$-algebra $R$, we respectively
define the space of modular forms and
the space of cusp forms,
\begin{align*}
	M_k \big(U(n), R\big) &= H^0\paren{X(n)_R, \omega^k}, \\
	S_k \big(U(n), R\big) &= H^0\paren{X(n)_R, \omega^k\big(-C(n)\big)}.
\end{align*}
These modular forms have Taylor expansions at each cusp with
$R[\zeta_n]$-coefficients.
Now we have the following $q$-expansion principle
(cf. \cite[Section~1.6]{katz-1973}).

\begin{prop}[$q$-expansion principle]
	\label{prop:q-series}
	Let $f$ be a modular form and $c$ a cusp of $X(n)$.
	Assume that the $q$-series of $f$ at $c$ vanishes, then
	$f$ vanishes on the connected component of $X(n)$ containing
	$c$.
\end{prop}

In practice, we only consider the cusp $c_u = (i \infty, a(u))$ for
$u \in (\Z/n\Z)^\times$ and write 
$a_f (\frac{i}{n}, u)$ for $a_{i, c_u}$. Then we have a $q$-expansion
at $c_u$ by 
\begin{equation}
	\label{eq:q-xn}
	f(z) = \sum_{i \geq 0} a_f\paren{\frac{i}{n}, u} q^{\frac{i}{n}}.
\end{equation}
The advantage of this expression is the invariance under pullback by
$X(n') \lra X(n)$ for any multiple $n'$ of $n$.

\subsection{Modular forms in infinite level}
The transition maps in the towers
$\{Y(n)_\Q\}_n$ and $\{X(n)_\Q\}_n$,
indexed by positive integers ordered by divisibility,
are finite and hence affine.
We take the inverse limits below in the category of schemes
over $\Q$, and the direct limit defining $S_k$ in the category
of $\Q$-vector spaces, with transition maps given by pullback:
\begin{align*}
  Y &:= \varprojlim_n Y(n)_\Q, \\
  X &:= \varprojlim_n X(n)_\Q, \\
  S_k &:= \varinjlim_n S_k\big(U(n),\Q\big).
\end{align*}

Over the complex numbers, we have
\begin{align*}
	Y(\C) &= \SL_2(\Z) \bs \calh \times \GL_2(\wh \Z), \\
	X(\C) &= \SL_2(\Z) \bs \wh\calh \times \GL_2(\wh \Z).
\end{align*}
The determinant map identifies the sets of connected components:
\begin{equation}
  \label{eq:pi_0}
  \pi_0\paren{Y(\C)}
  \iso
  \pi_0\paren{X(\C)}
  \iso
  \Q_+^\times \bs \wh\Q^\times
  \cong
  \wh\Z^\times.
\end{equation}
Indeed, the first identification follows by passing to the
inverse limit of the finite-level decompositions above
(equivalently, from strong approximation for $\SL_2$),
and the second follows from
\[
  \wh\Q^\times
  =
  \Q_+^\times \wh\Z^\times.
\]

Let
\[
  \calh^\pm
  :=
  \set{z \in \C \Mid \operatorname{Im}(z) \neq 0}
\]
denote the union of the upper and lower half-planes, and set
\[
  \wh\calh^\pm
  :=
  \calh^\pm \cup \P^1(\Q).
\]
Using the identity
\[
	\GL_2(\wh \Q) = \GL_2(\Q)_+ \GL_2(\wh \Z),
\]
we can write the above as
\begin{align*}
	Y(\C) &= \GL_2(\Q) \bs \calh^\pm \times \GL_2(\wh \Q), \\
	X(\C) &= \GL_2(\Q) \bs \wh\calh^\pm \times \GL_2(\wh \Q).
\end{align*}
In this terminology, a cusp form $f \in S_k(\C)$ is a function $f$ on
$\calh^\pm \times \GL_2(\wh\Q)$ such that,
\begin{enumerate}
	\item $f$ is invariant under right translation of some open subgroup
		$U$ of $\GL_2(\wh \Q)$;
	\item $f$ is holomorphic in $z$;
	\item for any $\gamma \in \GL_2(\Q)$,
		\[
			f(\gamma z, \gamma g) = j(\gamma, z)^k f(z, g);
		\]
	\item $f$ vanishes at each cusp.
\end{enumerate}

\subsubsection*{$q$-expansions}
Notice that the set of cusps in level $U$ is given by 
\[
	C_U(\C) = \GL_2(\Q)_+ \bs \P^1(\Q) \times \GL_2(\wh \Q) / U \iso
		\set{i \infty} \times B(\Q)_+ \bs \GL_2(\wh \Q) / U.
\]
Notice that $N(\Q)$ is dense in $N(\wh \Q)$.
So taking a limit gives, 
\[
	C(\C) \iso \set{i \infty} \times N(\wh \Q) M(\Q)_+
		\bs \GL_2(\wh \Q),
\]
where $M(\Q)_+$ denotes the group of diagonal matrices in $\Q$ with
positive determinant.
Then for the last condition, we need only consider the cusp
represented by $i \infty \times \GL_2(\wh\Q) / U$.
Assume that $c$ is represented by $(i \infty, gU)$, then we have the
stabilizer
\[
	N_c := N(\Q) \cap gUg^{-1} = N(m\Z)
\]
for some $m \in \N$. Then we have the $q$-expansion
\[
	f(z, g) = \sum_{n \geq 1} A_f\paren{\frac{n}{m}, g} q^{\frac{n}{m}}.
\]
So we defined a function
$A_f: \Q^\times_+ \times \GL_2(\wh \Q) \lra \C$,
which does not depend on the choice of $U$.

\subsubsection*{Relation to automorphic forms}

Now we describe how to consider modular forms as automorphic forms
(cf. \cite[Section~3]{casselman}, \cite[Section~3]{gelbart},
\cite[Section~3.6]{bump}).
For a modular form $f$ of weight $k$,
first define a function 
on $\GL_2(\R)_+ \times \GL_2(\wh \Q)$,
\begin{equation}
	\label{eq:varphi-f}
	\varphi(g) := f(g_\infty(i), g^\infty) j(g_\infty, i)^{-k},
\end{equation}
where $g = (g_\infty, g^\infty)$.
Then for $\gamma \in \GL_2(\Q)_+$,
\begin{align*}
	\varphi(\gamma g)
		&= f\paren{\gamma g_\infty(i), \gamma g^\infty}
			j\paren{\gamma g_\infty, i}^{-k} \\
		&= f\paren{\gamma g_\infty(i), \gamma g^\infty}
			j\paren{\gamma, g_\infty(i)}^{-k} j\paren{g_\infty, i}^{-k} \\
		&= f\paren{g_\infty(i), g^\infty}
			j\paren{g_\infty, i}^{-k} \\
		&= \varphi(g).
\end{align*}
From the construction, it is clear that $\varphi$ is invariant under 
$\R_+$ and has weight $k$ under $\SO_2$, as
\[
	j\paren{
		\begin{pmatrix}
			\cos\theta & \sin\theta \\
			-\sin\theta & \cos\theta
		\end{pmatrix}, i}
		= -(\sin\theta) i + \cos \theta
		= e^{-i\theta}.
\]
So $\varphi$ is invariant under $\GL_2(\Q)_+$ on both factors.
This function can be uniquely extended to a function on 
$\GL_2(\A)$ that is invariant under left translation by 
defining $\varphi(g) := \varphi(h g)$ for
$h \in \GL_2(\Q)_{-}$ and
$g \in \GL_2(\A) \setminus (\GL_2(\R)_+ \times \GL_2(\wh \Q))$.
For details, see \cite[Section~3]{casselman} and \cite[Section~3.6]{bump}.

With $z = x + yi$, we can also recover $f$ from $\varphi$ by,
\begin{equation}
	\label{eq:f-varphi}
	f\paren{z, g^\infty} = y^{-k/2}
		\varphi\paren{\begin{pmatrix} y & x \\ 0 & 1 \end{pmatrix}, g^\infty}.
\end{equation}
Comparing the Fourier expansions of both sides of
Equation~\ref{eq:f-varphi}, we obtain
\[
  \sum_{r \in \Q^\times_+}
    A_f\paren{r, g^\infty} e^{2\pi i r z}
  =
  y^{-k/2}
  \sum_{r \in \Q^\times}
    W_\varphi\paren{
      \begin{pmatrix} r & 0 \\ 0 & 1 \end{pmatrix}
      \begin{pmatrix} y & x \\ 0 & 1 \end{pmatrix},
      \begin{pmatrix} r^\infty & 0 \\ 0 & 1 \end{pmatrix}
      g^\infty
    }.
\]
Here $W_\varphi(g)$ is the Whittaker function of $\varphi$:
\[
  W_\varphi(g)
  :=
  \int_{[N]}
    \varphi(ng)\psi_N^{-1}(n)\,dn.
\]
Since
\[
  \begin{pmatrix} r & 0 \\ 0 & 1 \end{pmatrix}
  \begin{pmatrix} y & x \\ 0 & 1 \end{pmatrix}
  =
  \begin{pmatrix} ry & rx \\ 0 & 1 \end{pmatrix},
\]
we obtain
\begin{align*}
  \sum_{r \in \Q^\times_+}
    A_f\paren{r, g^\infty}e^{2\pi i r z}
  &=
  y^{-k/2}
  \sum_{r \in \Q^\times}
    W_\varphi\paren{
      \begin{pmatrix} ry & rx \\ 0 & 1 \end{pmatrix},
      \begin{pmatrix} r^\infty & 0 \\ 0 & 1 \end{pmatrix}
      g^\infty
    } \\
  &=
  y^{-k/2}
  \sum_{r \in \Q^\times}
    e^{2\pi i r x}
    W_\varphi\paren{
      \begin{pmatrix} ry & 0 \\ 0 & 1 \end{pmatrix},
      \begin{pmatrix} r^\infty & 0 \\ 0 & 1 \end{pmatrix}
      g^\infty
    }.
\end{align*}
Comparing the coefficients of $e^{2\pi i r x}$, we obtain
an equality for each $r$ (see \cite[Lemma 3.6]{gelbart}):
\[
  W_\varphi\paren{
    \begin{pmatrix} ry & 0 \\ 0 & 1 \end{pmatrix},
    \begin{pmatrix} r^\infty & 0 \\ 0 & 1 \end{pmatrix}
    g^\infty
  }
  =
  \begin{cases}
    A_f\paren{r, g^\infty}
      y^{k/2}e^{-2\pi r y}
      & \text{if $r>0$,} \\
    0 & \text{if $r<0$.}
  \end{cases}
\]

Equivalently, for $x \in \R$ and $y > 0$, one has
\begin{equation}
\label{eq:def-wk}
  W_k\begin{pmatrix}y & x \\ 0 & 1\end{pmatrix}
  :=
  y^{k/2}e^{2\pi i(x + iy)}.
\end{equation}
At $x = 0$, this agrees with
Equation~\ref{eq:Whittaker-GL2R}.
Then the adelization of classical modular forms can be
characterized by the above comparison of Fourier coefficients
(the following is essentially \cite[Proposition 3.1 and Lemma 3.6]{gelbart}).
\begin{prop}
\label{prop:f-phi}
	The Whittaker function $W_\varphi$ of $\varphi$ factors as
	\[
		W_\varphi(g) = W_k\paren{g_\infty} W_\varphi\paren{g^\infty}
	\]
	Moreover, Equations~\ref{eq:varphi-f} and \ref{eq:f-varphi} give
	a one-to-one correspondence between holomorphic forms of weight $k$
	and automorphic forms of weight $k$ with the following compatibility
	of Fourier coefficients for any $r \in \Q_+^\times$:
	\[
		A_f\paren{r, g^\infty} =
			r^{k/2} W_\varphi \paren{
				\begin{pmatrix}r^\infty & 0 \\ 0 & 1 \end{pmatrix} g^\infty},
	\]
\end{prop}

In the correspondence of Proposition~\ref{prop:f-phi},
we call $f$ the {\em modular avatar} of $\varphi$, and
$\varphi$ the {\em automorphic avatar} of $f$.

\subsection{Galois action and
	\texorpdfstring{$q$}{q}-expansion principle}
Recall that
\[
  S_k := \varinjlim_n S_k\big(U(n),\Q\big),
\]
and consider the space of weight-$k$
cusp forms defined over $\C$
\[
  S_k(\C) := S_k \otimes_\Q \C.
\]
Let $\cala_{0,k}([\GL_2])$ denote the space of cuspidal
automorphic forms $\varphi$ of weight $k$.
Both have an action by $\GL_2(\wh \Q)$.
Equations~\ref{eq:varphi-f} and \ref{eq:f-varphi} define an
isomorphism between these representations of $\GL_2(\wh \Q)$:
\[
	S_k(\C) \iso \cala_{0, k} \paren{\sbrac{\GL_2}}.
\]

\subsubsection*{Galois action}
$\Aut(\C/\Q)$ acts on $S_k(\C)$
and we can recover $S_k$ as the invariants of this action.
This induces an action on $\cala_{0, k} \paren{\sbrac{\GL_2}}$
via the isomorphism from $S_k(\C)$.
In the following, we write down the corresponding formula for
the induced Galois action on $\cala_{0, k} \paren{\sbrac{\GL_2}}$.
First, we describe the Galois action on cusp forms 
in terms of $q$-expansions.

Notice that the set of cusps is defined over $\Q^\ab$,
the maximal abelian extension of $\Q$.
The action of $\gal(\Q^\ab/\Q)$ on this set is given by the
left action by an element
$a_\sigma :=
\begin{psmallmatrix}\lambda_\sigma & 0 \\ 0 & 1 \end{psmallmatrix}$
where $\lambda: \gal(\Q^\ab/\Q) \iso \wh\Z^\times$ is the
isomorphism induced from action on roots of unity and
$\sigma \in \Aut(\C)$.
From the action of $\Aut(\C)$ on the space of modular forms
(via its action on cusps and coefficients),
$f^\sigma(z,g)$ has Fourier coefficients given by
\[
	A_{f^\sigma}(r, g) = A_f\paren{r, a_\sigma^{-1} g}^\sigma.
\]

By Proposition \ref{prop:f-phi},
\[
	W_\varphi\paren{g^\infty} = A_f\paren{1, g^\infty}.
\]
Thus an action of $\Aut(\C)$
can be defined on the space
$\calw(\psi^\infty)$ of the Whittaker function on $\GL_2(\A^\infty)$ by
\begin{align}
	\label{eq:W-galois}
	W &\lra W^\sigma \\
	g^\infty &\longmapsto W\paren{a_\sigma^{-1} g^\infty}^\sigma. \nonumber
\end{align}

\subsubsection*{\texorpdfstring{$q$}{q}-expansion principle}
By the $q$-expansion principle in finite level
of Proposition~\ref{prop:q-series}, 
we have a $q$-expansion principle in infinite level as well:
a modular form $f$ vanishes on $X$ if and only if its
$q$-expansion vanishes on at least one cusp for each connected
component of $X$.

By Equation~\ref{eq:pi_0},
the set of connected components of $X$ is given by $X_u$ for
$u \in \wh \Z^\times$, the connected component of $X$ containing 
the image of $(z, a(u)) \in \wh\calh \times \GL_2(\A^\infty)$.
We can define the standard cusp on $X_u$ by
$c_u = (i\infty, a(u))$. Then by the $q$-expansion principle,
$f$ is determined by its $q$-expansion at the cusp $c_u$.
Let $\kappa_\varphi$ be the Kirillov model of the
automorphic avatar $\varphi$ of $f$.

\begin{proposition}
\label{prop:a-kappa}
  Let $\varphi \in \cala_{0, k} \paren{\sbrac{\GL_2}}$
  be the automorphic avatar of $f \in S_k(\C)$.
  Then for each $u \in \wh \Z^\times$,
  $f$ has the $q$-expansion at $c_u$
  \[
    f(q, u) = \sum_{r \in \Q_+^\times } a_f(r, u) q^r
  \]
  such that
  \[
    a_f(r, u) = r^{\frac{k}{2}} \kappa_{\varphi}(ru).
  \]
\end{proposition}

\begin{proof}
  Denote for each $u \in \wh \Z^\times$
  (note the abuse of notation),
  \begin{align*}
  	f(z, u) &= f\big(z, a(u)\big), \\
  	a_f(r, u) :&= A_f\big(r, a(u)\big).
  \end{align*}
  Then the $q$-expansion of $f$ at $c_u$ is given by 
  \[
  	f(q, u) = \sum_{r \in \Q_+^\times } a_f(r, u) q^r.
  \]

  The Galois action of $\Aut(\C)$ on $q$-expansions can then be
  written for each $u \in \wh \Z^\times$ as
  \begin{equation}
  	\label{eq:f-sigma}
  	f^\sigma(q, u) = \sum_{r \in \Q_+^\times}
  		a_f\paren{r, \lambda _\sigma ^{-1}u}^\sigma q^r.
  \end{equation}

  Recall that the Kirillov model for
  a cuspidal automorphic form
  $\varphi \in \cala_{0, k}([\GL_2])$ is
  \[
  	\kappa_\varphi(x) = W_\varphi \big(a(x)\big).
  \]
  By Proposition \ref{prop:f-phi},
  \[
    a_f(r, u) = r^{\frac{k}{2}} \kappa_{\varphi}(ru)
  \]
  as desired.
\end{proof}

Due to the decomposition $\A^{\infty, \times} =
\Q_+^\times\times \wh \Z^\times$, Proposition~\ref{prop:a-kappa}
allows one to recover $a_f$ and $\kappa_\varphi$ from each other.

\subsection{Newforms}

\subsubsection*{Decompositions}
We study the decomposition of $S_k$ into the direct sum of
irreducible representations of $\GL_2(\wh \Q)$. We may do this by
first working on $\C$ and then studying the action 
by $\Aut(\C)$ later. Over the complex numbers, via
Equations~\ref{eq:f-varphi} and \ref{eq:varphi-f},
there is an isomorphism:
\begin{equation}
	\label{eq:SA}
	S_k(\C) \iso \cala_{0, k}\paren{\sbrac{\GL_2}}.
\end{equation}
We know that the right-hand side is a subspace $\cala_{0}([\GL_2])$ of
cusp forms which can be decomposed into irreducible representations:
\[
	\cala_{0}\paren{\sbrac{\GL_2}} = \bigoplus_\pi \pi,
\]
where $\pi$ ranges over all irreducible cuspidal representations of
$\GL_2(\A)$.
Notice that by their Whittaker functions, each $\pi$ has a
further decomposition into irreducible representations
of $\GL_2(\Q_p)$:
\[
	\pi = \bigotimes_{p \leq \infty} \pi_p.
\]

We define $\calw_k \subset \calw(\psi_\infty)$ to be the
representation of $\GL_2(\R)$ generated by the weight-$k$
holomorphic Whittaker function $W_k$
defined in Equation~\ref{eq:def-wk}.
Then we have
\[
  \cala_{0,k}\paren{\sbrac{\GL_2}}
  =
  \bigoplus_{\pi_\infty = \pi_k}
    \calw_k \otimes
    \calw\paren{\pi^\infty,\psi^\infty}
  \iso
  \bigoplus_{\pi_\infty = \pi_k}
    \calw\paren{\pi^\infty,\psi^\infty}.
\]
 
Combining this with the isomorphism from Equation~\ref{eq:SA},
we obtain a decomposition of $S_k(\C)$ into the direct sum of
irreducible representations:
\begin{equation}
	\label{eq:S-dec}
	S_k(\C) \iso
		\bigoplus_{\pi_\infty = \pi_k} \calw\paren{\pi^\infty, \psi^\infty}.
\end{equation}

\subsubsection*{Definition of newforms}
For each irreducible representation $\pi$ contributing to
the right-hand side of Equation~\ref{eq:S-dec}, let
$\pi^\infty$ denote its finite part. The Atkin--Lehner
theory of newforms \cite{atkin-lehner} gives a least
positive integer $N$ such that
$\paren{\pi^\infty}^{U_1(N)} \neq 0$,
where
\[
  U_1(N)
  :=
  \set{
    u \in \GL_2(\wh \Z)
    \Mid
    u \equiv
    \begin{pmatrix}* & * \\ 0 & 1\end{pmatrix}
    \pmod{N}}.
\]
The space $\paren{\pi^\infty}^{U_1(N)}$ is
one-dimensional. Let
\[
  W^{\new,\infty}
  \in
  \calw\paren{\pi^\infty, \psi^\infty}^{U_1(N)}
\]
be its normalized Whittaker vector, so that
$W^{\new,\infty}(e) = 1$, and set
\[
  W^\new
  :=
  W_k \otimes W^{\new,\infty}.
\]
The corresponding global newform is
\[
  \varphi^\new(g)
  :=
  \sum_{a \in \Q^\times}
  W^\new\paren{
    \begin{pmatrix}a & 0 \\ 0 & 1\end{pmatrix}g}.
\]
Let $\omega$ be the central character of $\pi$.
Then $\varphi^\new$ has character $\omega$ under the
action of the larger group
\[
  U_0(N)
  :=
  \set{
    u \in \GL_2(\wh \Z)
    \Mid
    u \equiv
    \begin{pmatrix}* & * \\ 0 & *\end{pmatrix}
    \pmod{N}}
  =
  \wh\Z^\times \cdot U_1(N).
\]
 
The isomorphism in Equation~\ref{eq:S-dec} gives a corresponding
weight-$k$ cusp form $f^\new \in S_k(\C)$. 
We show that $f^\new$ is the classical newform for
$\Gamma_1(N)$ with nebentypus $\omega^\infty$,
so we may equivalently consider the corresponding
$\pi$ and $\varphi^\new$ to be ``newforms''
(for details, see \cite[Theorem 5.19]{gelbart}).

\begin{prop}
\label{prop:new-new}
	There are natural one-to-one correspondences between the following
	objects:
	\begin{enumerate}
		\item irreducible subrepresentations $\pi$ of $\cala_{0, k}([\GL_2])$
			under $\GL_2(\A)$;
		\item newforms $\varphi^\new$ in $\cala_{0, k}([\GL_2])$;
		\item newforms $f^\new$ in $S_k(\C)$;
		\item irreducible subrepresentations $\pi^\infty$ of $S_k(\C)$
			under $\GL_2(\wh \Q)$.
	\end{enumerate}
\end{prop}

\begin{proof}[Sketch of proof]
	First, since $f^\new$ is invariant under $U_1(N)$, we see that 
	$f^\new$ is a modular form on the modular curve
	\[
		X_{U_1(N), \C} :=
			\GL_2(\Q)_+ \bs \calh \times \GL_2(\wh\Q) / U_1(N).
	\]
	Use the decomposition $\GL_2(\wh \Q) = \GL_2(\Q)_+ \cdot U_1(N)$ to
	see that this modular curve is actually the classical modular curve
	\[
		X_1(N) = \Gamma_1(N) \bs \calh,
	\]
	where 
	\[
		\Gamma_1(N) := \set{\gamma \in \SL_2(\Z) \Mid
			\gamma \equiv \begin{pmatrix}* & * \\ 0 & 1 \end{pmatrix}
			\pmod {N}}
		= U_1(N) \cap \GL_2(\Q)_+.
	\]
	Therefore, $f^\new$ is a classical modular form for
	$\Gamma_1(N)$.
	 
	Second, since $\varphi^\new$ has character $\omega$ under $U_0(N)$,
	$f^\new$ has character $\omega$ under
	\[
		\Gamma_0(N) := \set{\gamma \in \SL_2(\Z) \Mid \gamma \equiv
			\begin{pmatrix}* & * \\ 0 & * \end{pmatrix} \pmod {N}}
			= U_0(N) \cap \GL_2(\Q)_+,
	\]
	as follows:
	\[
		\omega \paren{\begin{pmatrix}a & b \\ 0 & d \end{pmatrix}} =
			\omega^\infty(d).
	\]
	Therefore $f^\new$ has nebentypus $\omega^\infty$, i.e.
	$f^\new \in S_k \paren{\Gamma_1(N), \omega^\infty}$.
	 
	Third, $W^{\new,\infty}$ lies in the
  one-dimensional space
  $\paren{\pi^\infty}^{U_1(N)}$, which is an
  eigenspace for the Hecke algebra
  \[
    \T_1(N)
    :=
    \C\sbrac{
      U_1(N) \bs \GL_2(\wh \Q) / U_1(N)}.
  \]
  Thus $\varphi^\new$, and hence $f^\new$, is an
  eigenform for
  \[
    \C\sbrac{
      \Gamma_1(N) \bs \GL_2(\Q) / \Gamma_1(N)}
    \iso
    \T_1(N).
  \]
  In particular,
  $f^\new \in S_k(\Gamma_1(N), \omega^\infty)$ is an
  eigenform.

  Fourth, since $N$ is the least integer for which
  $\paren{\pi^\infty}^{U_1(N)} \neq 0$, it is also the
  minimal level of $f^\new$. Thus $f^\new$ is a
  newform of level $N$.
   
  Finally, since $\begin{psmallmatrix}1 & \Z \\ 0 & 1 \end{psmallmatrix}
  \subset \Gamma_1(N)$, $f^\new(z)$ has a $q$-expansion:
  \[
    f^\new = \sum_{n \geq 1} a_n q^n.
  \]
  By Proposition~\ref{prop:a-kappa},
  $a_n = n^{\frac{k}{2}} \cdot \kappa\paren{n^\infty}$.
  This shows that $a_1 = 1$
  (recall that $W^{\new,\infty}(e) = 1$);
  combined with the previous steps,
  $f^\new$ is a normalized newform in
  $S_k(\Gamma_1(N), \omega^\infty)$.
  Furthermore, $f^\new$ generates the irreducible
  subrepresentation $\pi(f^\new)$ of $S_k(\C)$ that
  corresponds to $\pi^\infty$ under the isomorphism in
  Equation~\ref{eq:S-dec}.

  Conversely, starting with a newform $f^\new$, we
  can reverse the above procedure to construct a
  newform $\varphi^\new \in \cala_{0, k}([\GL_2])$
  whose finite Whittaker function has minimal level
  $U_1(N)$ and is normalized by $\kappa(1) = 1$.
  It is well-known that such an automorphic form
  generates an irreducible subrepresentation
  $\pi(\varphi^\new)$ of $\cala_{0, k}([\GL_2])$.
\end{proof}

\subsubsection*{Rationality and integrality}
The correspondences in
Proposition~\ref{prop:new-new} are $\Aut(\C)$-equivariant, with
the action on the automorphic side
given by Equation~\ref{eq:W-galois}.
Then the objects in Proposition~\ref{prop:new-new} also have the same
field of definition. Such a field is largely
easy to describe in terms of newforms (as modular forms):
if $f$ is a newform
with $q$-expansion $\sum_n a_n q^n$, then for any
$\sigma \in \Aut(\C)$, the form $f^\sigma$ will have $q$-expansion 
$\sum_n a_n^\sigma q^n$. In fact, if the $q$-expansion of 
$f^\sigma$ is $\sum_n b_n q^n$, then 
\begin{align*}
	a_n &= a_f(n, 1), \\
	b_n &= a_{f^\sigma}(n, 1).
\end{align*}
By Equation~\ref{eq:f-sigma},
\[
	b_n = a_{f^\sigma}(n, 1)
		= a_f\paren{n, \lambda_\sigma^{-1}}^\sigma
		= a_f(n, 1)^\sigma
		= a_n^\sigma
\]
where in the last step, we use the fact that $f$ on
$\calh \times \GL_2(\wh \Q)$ is invariant under $U_1(N)$.
Therefore the field of definition of $f$ is the subfield
$\Q(f)$ of $\C/\Q$ generated by $\set{a_n}$.

Another advantage of using $f$ for this description is that the
coefficients $a_n$ are always algebraic integers. They all come from Galois
representations.
\begin{thm}[{Eichler--Shimura~\cite{eichler-1957, shimura}
    for $k = 2$, Deligne~\cite{deligne} for $k > 2$,
    and Deligne--Serre~\cite{deligne-serre} for $k = 1$}]
  Let $f$ be a newform of level $N$ with $q$-expansion
  $\sum_n a_n q^n$. For every
  prime $\lambda$ of $\Q(f)$ above a rational prime $\ell$,
  there is a continuous Galois representation
  \[
    \rho_{f,\lambda}:
    \gal\paren{\overline{\Q}/\Q}
    \longrightarrow
    \GL_2\paren{\Q(f)_{\lambda}}
  \]
  such that, for every prime $p \nmid \ell N$,
  \[
    a_p
    =
    \Tr\paren{\rho_{f,\lambda}(\Frob_p)}.
  \]
\end{thm}

Then we can define the ring $\calo(f)$ as the ring of integers of
$\Q(f)$. Then all of the objects in Proposition \ref{prop:new-new}
have an integral model defined over $\calo(f)$.

\section{The Harris--Venkatesh period}
\label{sec:hv-period}

We extend the constructions of
Harris--Venkatesh~\cite{hv} and
Darmon--Harris--Rotger--Venkatesh~\cite{dhrv} by defining
the Harris--Venkatesh period $\calp_\hv$
on the tensor product of two spaces of weight-one cusp forms.
It recovers the pairing with
the Shimura class $\mfs_p$ used in \cite{hv, dhrv}.

\subsection{Modular curves}

We recall the basics about modular curves,
first for a positive integer $N$ and then
for a fixed finite set of primes $\Sigma$.
Recall that, for each positive integer $N$,
$X(N)$ denotes the modular curve over $\Z[1/N]$.
Its ring of global functions is
\[
  H^0\paren{X(N), \calo_{X(N)}}
  =
  \Z\sbrac{1/N, \zeta_N}.
\]
Thus the field of constants of its generic fiber is
$\Q(\zeta_N)$.

If $N \geq 3$, then there is a bundle $\omega$ of weight-one forms on
$X(N)$ and a bundle $\Omega = \Omega_{X(N)/\Z[1/N]}$ of relative
differentials with the Kodaira--Spencer map
\[
	\ks: \omega \otimes \omega \iso \Omega (C(N)),
\]
where $C(N)$ is the cuspidal divisor on $X(N)$.
For each positive integer $N$,
the Kodaira--Spencer map induces a pairing
\begin{align}
\label{eq:KS-N}
  H^0\paren{X(N),\omega(-C(N))}
  \otimes_{\Z[1/N]}
  H^0\paren{X(N),\omega(-C(N))}
  \longrightarrow
  H^0\paren{X(N),\Omega(-C(N))}.
\end{align}
This is the product of two weight-one cusp forms,
followed by the Kodaira--Spencer isomorphism.
Serre duality defines another pairing for each positive integer $N$,
\begin{align}
\label{eq:SD-N}
	H^0\paren{X(N), \Omega_{X(N)}} \otimes_{\Z[1/N]} H^1\paren{X(N), \calo_{X(N)}}
		\longrightarrow H^1\paren{X(N), \Omega}
		\iso \Z \sbrac{1/N, \zeta_N},
\end{align}
where the last isomorphism comes from the trace map.
Both of these pairings are compatible with pullback maps and
the action by $\GL_2(\Z/N\Z)$.

Now fix a finite set $\Sigma$ of primes and consider the projective
system of smooth curves over $\Z[1/\Sigma]$ indexed by positive
integers $N$, with prime factors in $\Sigma$:
\[
	X_\Sigma(N) := X(N) \otimes_{\Z[1/N]} \Z[1/\Sigma].
\]
Let
\[
  X_\Sigma := \varprojlim_N X_\Sigma(N),
\]
where the limit is taken in the category of schemes over
$\Z[1/\Sigma]$, and $N$ ranges over the positive integers
supported on $\Sigma$, ordered by divisibility. Set
\begin{align*}
  \Q_\Sigma &:= \prod_{q \in \Sigma} \Q_q,
  \\
  \Z_\Sigma &:= \prod_{q \in \Sigma} \Z_q.
\end{align*}
The subgroup $\GL_2(\Z_\Sigma)$ acts on $X_\Sigma$ by
changing the full level structure. More generally, an element
$g \in \GL_2(\Q_\Sigma)$ sends the standard Tate-module
lattice to a commensurable lattice. The associated
$\Sigma$-primary quasi-isogeny induces compatible maps
between sufficiently deep levels, since every prime in
$\Sigma$ is inverted on the base. Thus $X_\Sigma$ carries
a right action of $\GL_2(\Q_\Sigma)$.
Taking limits, we obtain the following two pairings.
\begin{equation}
\label{eq:KS-Sigma}
	H^0\paren{X_\Sigma, \omega(-C_\Sigma)} \otimes_{\Z[1/\Sigma]}
	H^0\paren{X_\Sigma, \omega(-C_\Sigma)}
	\longrightarrow H^0\paren{X_\Sigma, \Omega(-C_\Sigma)}.
\end{equation}
\begin{equation}
\label{eq:SD-Sigma}
	H^0\paren{X_\Sigma, \Omega_{X_\Sigma}} \otimes_{\Z[1/\Sigma]}
	H^1\paren{X_\Sigma, \calo_{X_\Sigma}}
	\longrightarrow H^1\paren{X_\Sigma, \Omega}
	\iso \Z\sbrac{1/\Sigma, \mu_\Sigma},
\end{equation}
where $\mu_\Sigma$ is the group of roots of unity whose order
is divisible only by primes in $\Sigma$.
The last isomorphism in Equation \ref{eq:SD-Sigma}
is deduced from the inductive system
of modified trace maps
\[
  \begin{tikzcd}[row sep = tiny]
      H^1\paren{X(N), \Omega} \arrow[r, "\widesim{}"]
      & \Z\sbrac{1/N, \zeta_N} \\
    x \arrow[r, mapsto]
      & \frac{\verts{\SL_2\paren{\Z/\big(\prod_{q \in \Sigma} q\big)\Z}}}
      {\verts{\SL_2(\Z/N\Z)}} \cdot \Tr(x),
  \end{tikzcd}
\]
where the numerator is included
to ensure that the entire ratio is in $\Z[1/\Sigma]$,
while the denominator is included to account for the degree of a cover.
Both of these pairings are compatible with pullback maps and
the action by $\GL_2(\Z/N\Z)$.
These pairings are compatible with the action by $\GL_2(\Q_\Sigma)$.
Notice that the action of $\GL_2(\Q_\Sigma)$
on $\Z \sbrac{\mu_\Sigma, 1/\Sigma}$ is given by 
\[
	\GL_2\paren{\Q_\Sigma} \xrightarrow{\det} \Q_\Sigma^\times
	\longrightarrow \Z_\Sigma^\times
	= \Aut\paren{\mu_\Sigma},
\]
where the second step is taking the standard projection for each
factor $\Q_p^\times$:
\[
	\Q_p^\times = p^\Z\times \Z_p^\times
		\longrightarrow \Z_p^\times.
\]

For each positive integer $N$ supported on $\Sigma$, define
\[
  \tau_N:
  \Z\sbrac{1/\Sigma, \mu_N}
  \longrightarrow
  \Z\sbrac{1/\Sigma}
\]
using the trace map $\Tr_{\Q(\mu_N)/\Q}: \Q(\mu_N) \rightarrow \Q$:
\[
  \tau_N(x)
  :=
  \frac{\prod_{p \in \Sigma}(1 - p)}{\phi(N)}
  \Tr_{\Q(\mu_N)/\Q}(x).
\]
If $N \mid N'$ are supported on $\Sigma$, then trace
transitivity and
\[
  \sbrac{\Q(\mu_{N'}):\Q(\mu_N)}
  =
  \frac{\phi(N')}{\phi(N)}
\]
show that $\tau_{N'}$ restricts to $\tau_N$. These maps
therefore define a $\GL_2(\Q_\Sigma)$-invariant map
\begin{equation}
\label{eq:tau}
  \tau:
  \Z\sbrac{1/\Sigma, \mu_\Sigma}
  \longrightarrow
  \Z\sbrac{1/\Sigma}.
\end{equation}
Concretely, if a root of unity $\zeta \in \mu_N$ has order
$N = \prod_{p \in \Sigma} p^{\alpha_p}$, then
\[
  \tau(\zeta)
  =
  \begin{cases}
    \prod_{\alpha_p = 0}(1 - p)
      & \text{if every $\alpha_p \leq 1$,} \\
    0 & \text{otherwise.}
  \end{cases}
\]
Then we get a new pairing,
\begin{equation}
	H^0 \paren{X_\Sigma, \Omega_{X_\Sigma}}
		\otimes_{\Z[1/\Sigma]} H^1\paren{X_\Sigma, \calo_{X_\Sigma}}
		\longrightarrow \Z[1/\Sigma],
\end{equation}
compatible with the action by $\GL_2(\Q_\Sigma)$.

\subsection{The Harris--Venkatesh period in finite level}

Fix a prime $p \geq 5$. Retain the notation
\[
  \Delta_p
  =
  \F_p^\times \otimes_\Z \Z\sbrac{1/6}
\]
from the introduction. In this subsection and the next
two subsections, we use the same notation for a modular
curve and its base change to $\Z\sbrac{1/6}$.

For every positive integer $N$ coprime to $p$, let
$X_0(p,N)$ be the modular curve associated with
\[
  U_0(p,N) := U_0(p)\cap U(N).
\]
Set
\[
  R_{p,N}
  :=
  \Z\sbrac{1/(6N)}/(p - 1).
\]
Equivalently, $R_{p,N}$ retains the primary
components of $p - 1$ at primes coprime to $6N$.

For each prime $\ell \geq 5$, the construction of
Harris--Venkatesh~\cite[Section~3.2]{hv} gives the
pullback of the $\ell$-primary Shimura class.
Taking their product gives
\[
  \mfs_p(N)
  \in
  H^1\paren{
    X_0(p,N)_{R_{p,N}},
    \calo}
  \otimes_{\Z\sbrac{1/6}}\Delta_p,
\]
which satisfies the following properties:
\begin{enumerate}
	\item $\mfs_p(N)$ is invariant under $\GL_2(\Z/N\Z)$;
	\item for any projection $\pi: X_0(p, N_2) \rightarrow X_0(p, N_1)$
		with $N_1 \mid N_2$, we have that $\pi^* \mfs_p(N_1) = \mfs_p(N_2)$.
\end{enumerate}
Notice that $U = U_0(p)$ has two embeddings into the maximal
subgroup $U(1)$: the trivial embedding $i_1$ and the embedding
$i_2$ obtained via conjugation by $\begin{psmallmatrix} 
p & \\ & 1 \end{psmallmatrix}$. This induces two projections
\[
	\pi_1, \pi_2: X_0(p, N) \longrightarrow X(N).
\]
The pullback on $\omega$ yields a pairing
\[
  \begin{tikzcd}[row sep = tiny]
  H^0\paren{X(N),\omega(-C(N))}
    \otimes
    H^0\paren{X(N),\omega(-C(N))}
    \arrow[r] &
    H^0\paren{X_0(p,N),\Omega(-C(N))}, \\
  \alpha\otimes\beta
    \arrow[r, mapsto] &
    \pi_1^*\alpha\cdot\pi_2^*\beta.
  \end{tikzcd}
\]
Composition with the pairing $\abrac{-, \mfs_p}$
(from Equation~\ref{eq:SD-N}) then gives a
$\GL_2\paren{\Z/N\Z}$-equivariant pairing on
cuspidal one-forms, compatible with pullbacks,
with values in
$\Z\sbrac{\zeta_N, 1/(6N)}
  \otimes_{\Z\sbrac{1/6}}\Delta_p$:
\[
  \begin{tikzcd}[row sep = tiny]
  \calp_\hv:
  H^0\paren{X(N),\omega(-C(N))}
  \otimes
  H^0\paren{X(N),\omega(-C(N))}
    \arrow[r] &
    \Z\sbrac{\zeta_N,1/(6N)}
    \otimes_{\Z\sbrac{1/6}}\Delta_p, \\
  \alpha\otimes\beta
    \arrow[r, mapsto] &
    \mfs_p\paren{\pi_1^*\alpha\cdot\pi_2^*\beta}.
  \end{tikzcd}
\]
We call $\calp_\hv$
the \emph{Harris--Venkatesh period}
(or \emph{Harris--Venkatesh pairing})
\emph{in level $N$}.

\subsection{The Harris--Venkatesh period in infinite level}
Fix a finite set $\Sigma$ of primes not containing $p$,
and consider the projective system $\{X_0(p, N)\}_{N}$
over all positive integers $N$ whose prime divisors are in $\Sigma$.
Let $X_{0, \Sigma}(p)$ denote its limit, which has a natural action by
$\GL_2(\Q_\Sigma)$.
Then the Harris--Venkatesh period in level $N$ induces a pairing at
infinite level:
\begin{equation}
\label{eq:HV-Sigma}
	H^0\paren{X_{\Sigma}, \omega(-C_\Sigma)} \otimes
		H^0\paren{X_{\Sigma}, \omega(-C_{\Sigma})}
		\longrightarrow \Z\sbrac{1/6, \mu_\Sigma, 1/\Sigma}
    \otimes_{\Z\sbrac{1/6}} \Delta_p.
\end{equation}
Composition of the pairing of Equation~\ref{eq:HV-Sigma} with the map
$\tau$ defined in Equation~\ref{eq:tau} yields
the \emph{Harris--Venkatesh period in infinite level},
\begin{equation}
\label{eq:HV-infinite}
  \calp_{\hv}:
    H^0\paren{X_{\Sigma}, \omega(-C_\Sigma)}
      \otimes
    H^0\paren{X_{\Sigma}, \omega(-C_\Sigma)}
    \longrightarrow
    \Z\sbrac{1/6, 1/\Sigma}
    \otimes_{\Z\sbrac{1/6}}\Delta_p.
\end{equation}
Note that the left-hand side of Equation~\ref{eq:HV-infinite}
consists of the spaces of cuspidal
modular forms of weight $1$ unramified outside of $\Sigma$.

\subsection{The period in the original setting of Harris--Venkatesh}
The period $\calp_\hv(\alpha \otimes \beta)$ is
related to the setting of
Harris--Venkatesh~\cite{hv} in the following way.
Let $f$ be a cuspidal newform of weight $1$ and
level $\Gamma_1(N)$ with coefficients generating
a subfield $\Q(f) \subset \C$ with ring of integers
$\calo(f)$. Assume that $N$ is coprime to $p$ and
set
\[
  R := \calo(f)\sbrac{1/(6N)}/(p - 1).
\]
Then we may regard $f$ and its dual $f^*$ as
elements of
\[
  H^0\paren{X_1(N)_R,
    \omega\paren{-C(N)}}.
\]
Define
\[
  G :=
    \Tr_{\Gamma_0(Np) / \Gamma_0(p)}
      \paren{f(z)f^*(pz)}
    \in H^0\paren{X_0(p)_R, \Omega}.
\]
Pairing with the Shimura class $\mfs_p$ on
$X_0(p)_R$ yields a value $\mfs_p(G)$ that is
related to the Harris--Venkatesh period by
\begin{align*}
	\mfs_p(G) &= \sum_{\gamma \in U_0(p) / U_0(Np)}
		\calp_\hv(\gamma f \otimes \gamma f^*) \\
		&= \sbrac{U(1) : U_0(N)} \cdot \calp_\hv(f \otimes f^*).
\end{align*}
Since $[U(1) : U_0(N)]$ is not necessarily invertible in
$R$ (its order at prime $\ell$ is $\sum_{q \mid N} \ord_\ell(q + 1)$),
the value $\calp_\hv(f \otimes  f^*)$ is more primitive than
$\abrac{G, \mfs_p}$.

%%%%%%%%%%%%%%%%%%%%%%%%%%%%%%%

\section{Liftings of pairings}
\label{sec:comparison}

Let $F$ be a non-archimedean local field of residue
characteristic $r$, let $\calo_F$ be its ring of integers,
fix a uniformizer $\varpi$, and let $q$ be the cardinality
of its residue field. Let $A \subset \C$ be a principal
ideal domain in which $r$ is invertible, set
$G := \GL_2(F)$, and let $\pi$ be a smooth,
infinite-dimensional irreducible representation of $G$
over $\C$.

We say that $\pi$ has an ${A}$\emph{-model}
$\pi_A \subset \pi$ if $\pi_A$ is an $A$-module such that
$\pi_A \otimes_A \C \iso \pi$, the module $\pi_A$ is
$G$-stable, and $\pi_A^H$ is free of finite type for every
compact open subgroup $H \leq G$
(cf. \cite{vigneras-1989}).

For each integer $c \geq 0$, set
\[
  U_1\paren{\varpi^c}
  :=
  \set{\gamma \in \GL_2\paren{\calo_F} \Mid
    \gamma \equiv
    \begin{pmatrix}
      * & * \\
      0 & 1
    \end{pmatrix}
    \pmod{\varpi^c}}.
\]
The conductor exponent $c(\pi)$ is the least $c$ such that
$\pi^{U_1(\varpi^c)} \neq 0$. By the newvector
theorem~\cite{casselman}, the space
\[
  \pi^\new := \pi^{U_1(\varpi^{c(\pi)})}
\]
is one-dimensional. If $\pi_A$ is an $A$-model of $\pi$,
then $\pi_A^\new := \pi_A \cap \pi^\new$ is free of rank
$1$ over $A$. In this section, we study pairings of
$A$-models generated by chosen newvectors.
	
\begin{prop}
\label{prop:varphi-A}
  Let $\pi_1$ and $\pi_2$ be smooth,
  infinite-dimensional irreducible representations of
  $\GL_2(F)$ that are dual in the sense that
  \[
    \Hom_{\C\sbrac{\GL_2(F)}}
    \paren{\pi_1 \otimes \pi_2, \C}
    \neq 0.
  \]
  Assume that their common conductor exponent
  \[
    c := c(\pi_1) = c(\pi_2)
  \]
  is positive. Let $\pi_{1,A}$ and $\pi_{2,A}$ be
  $A$-models generated by newvectors
  $v_{1,A}^\new$ and $v_{2,A}^\new$, respectively.
  There is a unique element
  \[
    \calp_0 \in
    \Hom_{A\sbrac{\GL_2(F)}}
    \paren{\pi_{1,A} \otimes \pi_{2,A}, A}
  \]
  such that
  \[
    \calp_0\paren{v_{1,A}^\new, v_{2,A}^\new}
    = q - 1.
  \]
\end{prop}

The main tool is a suitably normalized Haar measure on
$U_1\paren{\varpi^c}$.

\begin{lem}
\label{lem:haar-u0}
  There is a Haar measure $dh$ on
  $U_1\paren{\varpi^c}$ with values in
  $\Z\sbrac{1/r}$ and total volume $q - 1$.
\end{lem}

\begin{proof}
  Let
  \[
    H :=
    \set{\gamma \in U_1\paren{\varpi^c} \Mid
      \gamma \equiv
      \begin{pmatrix}
        1 & * \\
        0 & 1
      \end{pmatrix}
      \pmod{\varpi}}.
  \]
  Since $c \geq 1$, the group $H$ is pro-$r$, and
  upper-left reduction modulo $\varpi$ induces an
  isomorphism
  \[
    U_1\paren{\varpi^c}/H
    \iso
    \paren{\calo_F/\varpi\calo_F}^{\times}.
  \]
  Normalize the $\Z\sbrac{1/r}$-valued Haar measure on
  $H$ by $\vol(H) = 1$, and extend it to
  $U_1\paren{\varpi^c}$. Then
  \[
    \vol\paren{U_1\paren{\varpi^c}}
    =
    \verts{U_1\paren{\varpi^c}/H}
    =
    q - 1.
  \]
\end{proof}
	
\begin{proof}[Proof of Proposition~\ref{prop:varphi-A}]
  Since $\pi_1$ and $\pi_2$ are dual, choose a nonzero
  pairing
  \[
    \calp \in
    \Hom_{\C\sbrac{\GL_2(F)}}
    \paren{\pi_1 \otimes \pi_2, \C}.
  \]
  For $g_1,g_2 \in \GL_2(F)$, invariance gives
  \[
    \calp\paren{g_1v_{1,A}^\new,g_2v_{2,A}^\new}
    =
    \calp\paren{v_{1,A}^\new,
      g_1^{-1}g_2v_{2,A}^\new}.
  \]
  Averaging over $U_1\paren{\varpi^c}$ and using
  Lemma~\ref{lem:haar-u0}, we obtain
  \begin{align*}
    &(q - 1)
    \calp\paren{v_{1,A}^\new,
      g_1^{-1}g_2v_{2,A}^\new} \\
    &\qquad =
    \int_{U_1\paren{\varpi^c}}
      \calp\paren{hv_{1,A}^\new,
        hg_1^{-1}g_2v_{2,A}^\new}\,dh \\
    &\qquad =
    \calp\paren{v_{1,A}^\new,
      \int_{U_1\paren{\varpi^c}}
        hg_1^{-1}g_2v_{2,A}^\new\,dh}.
  \end{align*}
  The last integral belongs to
  $\pi_{2,A}^{U_1(\varpi^c)}
  = A v_{2,A}^\new$, so it is equal to
  $\lambda v_{2,A}^\new$ for some $\lambda \in A$.
  Consequently,
  \[
    (q - 1)
    \calp\paren{g_1v_{1,A}^\new,g_2v_{2,A}^\new}
    =
    \lambda
    \calp\paren{v_{1,A}^\new,v_{2,A}^\new}.
  \]
  If
  $\calp\paren{v_{1,A}^\new,v_{2,A}^\new}=0$,
  then the preceding identity shows that $\calp$ vanishes
  on the $G$-translates of the chosen newvectors. Since
  these generate $\pi_{1,A}$ and $\pi_{2,A}$, this
  contradicts the choice of $\calp$. Hence
  \[
    \calp\paren{v_{1,A}^\new,v_{2,A}^\new}
    \neq 0.
  \]
  Define
  \[
    \calp_0 :=
    \frac{q - 1}
    {\calp\paren{v_{1,A}^\new,v_{2,A}^\new}}
    \calp.
  \]
  The preceding identity shows that $\calp_0$ is
  $A$-valued. By Schur's lemma,
  \[
    \dim_{\C}
    \Hom_{\C\sbrac{\GL_2(F)}}
    \paren{\pi_1 \otimes \pi_2, \C}
    = 1,
  \]
  so the prescribed value on the newvectors determines
  $\calp_0$ uniquely.
\end{proof}
	
The main result of this section is the following multiplicity-one
type statement.
\begin{prop} 
\label{prop:AB-comparison}
	Let $\pi_1, \pi_2$ be two infinite-dimensional irreducible
	representations of $\GL_2(F)$ that are
	dual to each other in the sense that 
	\[
		\Hom_{\C\sbrac{\GL_2(F)}} \paren{\pi_1 \otimes \pi_2, \C} \neq 0,
	\]
  and assume that their common conductor exponent
  \[
    c := c(\pi_1) = c(\pi_2)
  \]
  is positive.
	Let $\pi_{1,A}$ and $\pi_{2,A}$ be $A$-models
  of $\pi_1$ and $\pi_2$, respectively, generated by
  newvectors $v_{1,A}^\new$ and $v_{2,A}^\new$.
  Let $A \twoheadrightarrow B$ be a surjective
  homomorphism of rings, and set
  \[
    \pi_{i,B} := \pi_{i,A} \otimes_A B,
    \qquad
    v_{i,B}^\new := v_{i,A}^\new \otimes 1.
  \]
	Then the cokernel of the homomorphism,
	\begin{align*}
		 \Hom_{A\sbrac{\GL_2(F)}} \paren{\pi_{1, A} \otimes \pi_{2, A}, A}
			&\lra \Hom_{B\sbrac{\GL_2(F)}}
			\paren{\pi_{1, B} \otimes \pi_{2, B}, B}, \\
		\calp	&\longmapsto \calp \otimes B
	\end{align*}
	is annihilated by $(q - 1)^2$.
	More precisely, for any $\calp_B \in
	\Hom_{B[\GL_2(F)]}\paren{\pi_{1, B}\otimes \pi_{2, B}, B}$, we have 
	\[
		(q - 1)^2 \cdot \calp_B = 
			(q - 1) \cdot \calp_B\paren{v_{1, B}^\new, v_{2, B}^\new}
			\cdot \calp_0 \otimes B,
	\]
	where $\calp_0$ is defined in Proposition \ref{prop:varphi-A}.
\end{prop}

We first need the following vanishing lemma.
\begin{lem}
	\label{lem:varphi-B}
	Let $\calp \in
	\Hom_{B[\GL_2(F)]} \paren{\pi_{1, B} \otimes \pi_{2, B}, B}$
	such that 
	\[
		\calp \paren{v_{1, B}^\new, v_{2, B}^\new} = 0.
	\]
	Then $(q - 1)\calp = 0$.
\end{lem}

\begin{proof}[Proof of Lemma~\ref{lem:varphi-B}]
  For any $g_1,g_2 \in \GL_2(F)$, the same averaging
  argument as in the proof of
  Proposition~\ref{prop:varphi-A} gives
  \begin{align*}
    &(q - 1)
    \calp\paren{v_{1,B}^\new,
      g_1^{-1}g_2v_{2,B}^\new} \\
    &\qquad =
    \calp\paren{v_{1,B}^\new,
      \int_{U_1\paren{\varpi^c}}
        hg_1^{-1}g_2v_{2,B}^\new\,dh}.
  \end{align*}
  The integral on the right is the image in
  $\pi_{2,B}$ of an element of
  $\pi_{2,A}^{U_1(\varpi^c)}=Av_{2,A}^\new$.
  It is therefore equal to
  $\lambda v_{2,B}^\new$ for some $\lambda \in B$.
  Hence
  \[
    (q - 1)
    \calp\paren{v_{1,B}^\new,
      g_1^{-1}g_2v_{2,B}^\new}
    =
    \lambda\calp\paren{v_{1,B}^\new,
      v_{2,B}^\new}
    =0.
  \]
  Since the chosen newvectors generate $\pi_{1,B}$ and
  $\pi_{2,B}$, invariance implies $(q - 1)\calp=0$.
\end{proof}
	
\begin{proof}[Proof of Proposition \ref{prop:AB-comparison}]
	Let $a \in A$ be a lift of $b := \calp_B
	(v_{1, B}^\new, v_{2, B}^\new) \in B$.
	Then 
	\[
		\calq := (q - 1) \calp_B - a\calp_0\otimes B \in
			\Hom_{B\sbrac{\GL_2(F)}} \paren{\pi_{1, B}\otimes \pi_{2, B}, B},
	\]
	vanishes at $(v_{1, B}^\new \otimes v_{2, B}^\new)$.
	By Lemma \ref{lem:varphi-B},
	$(q - 1) \calq  = 0$.
\end{proof}

%%%%%%%%%%%%%%%%%%%%%%%%%%%%%%%

\section{Proof of Theorem~\ref{thm:hv-reform}}
\label{sec:proof-hv-reform}
Let $L_f := \Q(f)$,
$F_f := \Q\sbrac{\chi_{\Ad(\rho_f)}}$, and
$\calo_{f,N} := \calo(f)\sbrac{1/N}$.
Let $\Sigma$ be the set of primes dividing $N$.
Let $\pi_{f, \Sigma, \calo_{f,N}}$ denote the
$\calo_{f,N}\sbrac{\GL_2(\Q_\Sigma)}$-subrepresentation of
\[
  H^0\paren{X_{\Sigma,\calo_{f,N}},
    \omega\paren{-C_\Sigma}}
\]
generated by $f$, and define
$\pi_{f^*, \Sigma, \calo_{f,N}}$ similarly. Set
\[
  \Pi_{f, \Sigma, \calo_{f,N}}
  :=
  \pi_{f, \Sigma, \calo_{f,N}}
    \otimes_{\calo_{f,N}}
  \pi_{f^*, \Sigma, \calo_{f,N}}.
\]
The tensor product of the one-variable correspondences in
Proposition~\ref{prop:f-phi} identifies
$\Pi_{f, \Sigma, \calo_{f,N}}$ with an $\calo_{f,N}$-module of
cuspidal two-variable modular forms. We denote the
form corresponding to $\varphi$ by $f_\varphi$; on pure
tensors, this convention is
\[
  f_{\varphi_1 \otimes \varphi_2}(z_1,z_2)
  =
  f_{\varphi_1}(z_1)f_{\varphi_2}(z_2),
\]
and it is extended $\calo_{f,N}$-linearly. We also set
\[
  f^\new := f \otimes f^*,
  \qquad
  f^\new(z_1,z_2) = f(z_1)f^*(z_2).
\]
Choose $\varphi \in \Pi_{f, \Sigma, \calo_{f,N}}$ such that
$\calp(\varphi) \neq 0$.

Fix a prime $p \geq 5$ coprime to $N$, a prime
$\ell \geq 5$ with $\ell^t \Vert p - 1$, and a
discrete logarithm
\[
  \log_\ell:
  \F_p^\times
  \longrightarrow
  \Z/\ell^t\Z.
\]
This map factors uniquely through $\Delta_p$. Moreover,
since $\ell \geq 5$, the natural map
\[
  \calo_{f,N} / \ell^t \calo_{f,N}
  \longrightarrow
  \calo_{f,N}\sbrac{1/6}/
    \ell^t \calo_{f,N}\sbrac{1/6}
\]
is an isomorphism. Consequently, although
$\calp_\hv$ is defined after inverting $6$, the
composition
\[
  \calp_{\hv,\ell}
  :=
  \paren{1 \otimes \log_\ell}\circ\calp_\hv
\]
is a well-defined $\calo_{f,N}/\ell^t\calo_{f,N}$-valued pairing on
$\Pi_{f, \Sigma, \calo_{f,N}}$.

For every maximal ideal $\mfl$ of $\calo_{f,N}$ above $\ell$,
the localization $\calo_{f,N,\mfl}$ is a discrete valuation
ring.
After localization, the integral models decompose as
\begin{align*}
  \pi_{f,\Sigma,\calo_{f,N,\mfl}}
  &\iso
  \bigotimes_{q \mid N}\pi_{f,q,\calo_{f,N,\mfl}},
  \\
  \pi_{f^*,\Sigma,\calo_{f,N,\mfl}}
  &\iso
  \bigotimes_{q \mid N}\pi_{f^*,q,\calo_{f,N,\mfl}},
\end{align*}
where each local model is generated by its newvector.
Write $\calp_{\hv,\ell,\mfl}$ for the base change of
$\calp_{\hv,\ell}$ to $\calo_{f,N,\mfl}/\ell^t\calo_{f,N,\mfl}$. Thus
Propositions~\ref{prop:varphi-A}
and~\ref{prop:AB-comparison} apply over
$\calo_{f,N,\mfl}$ and $\calo_{f,N,\mfl}/\ell^t\calo_{f,N,\mfl}$.

Set
\[
  D := \prod_{q \mid N}(q - 1).
\]
Tensoring the normalized local pairings gives an
$\calo_{f,N,\mfl}$-valued pairing $\calp_{0,\mfl}$ such that
\[
  \calp_{0,\mfl}\paren{f^\new} = D.
\]
Over $L_f$, multiplicity one gives
\[
  \calp_{0,\mfl}
  =
  \frac{D}{\calp\paren{f^\new}}\calp.
\]
Applying Proposition~\ref{prop:AB-comparison}
successively at the primes $q \mid N$ gives
\begin{equation}
\label{eq:comparison-localized}
  D^2\calp_{\hv,\ell,\mfl}
  =
  D\calp_{\hv,\ell,\mfl}\paren{f^\new}
    \calp_{0,\mfl}
\end{equation}
as $\calo_{f,N,\mfl}/\ell^t\calo_{f,N,\mfl}$-valued pairings.

Set
\[
  \alpha_\varphi
  :=
  \frac{\calp(\varphi)}{\calp\paren{f^\new}}
  \in L_f^\times.
\]
After clearing denominators, choose nonzero
$a_\varphi,b_\varphi \in \calo(f)$ such that
$\alpha_\varphi=a_\varphi/b_\varphi$.
Evaluating Equation~\ref{eq:comparison-localized} at
$f_\varphi$ and localizing at every $\mfl \mid \ell$
gives
\begin{equation}
\label{eq:comparison-cross-multiplied}
  b_\varphi D^2
  \calp_{\hv,\ell}\paren{f_\varphi}
  =
  a_\varphi D^2
  \calp_{\hv,\ell}\paren{f^\new}
\end{equation}
in $\calo_{f,N}/\ell^t\calo_{f,N}$.

Extend the regulator $\calo(f)$-linearly to
\[
  \calu_f\paren{\Ad(\rho_f)}
  :=
  \calu\paren{\Ad(\rho_f)}
    \otimes_{R_f}\calo(f).
\]
Choose $a_\varphi^\vee,b_\varphi^\vee \in \calo(f)$
and positive integers $n_a,n_b$ such that
\[
  a_\varphi^\vee a_\varphi=n_a,
  \qquad
  b_\varphi^\vee b_\varphi=n_b.
\]
Such elements exist by applying the adjugate construction
to multiplication by $a_\varphi$ and $b_\varphi$ on the
free $\Z$-module $\calo(f)$.

Suppose first that the Harris--Venkatesh conjecture
for $f^\new$ holds, so that
\[
  m_f\calp_{\hv,\ell}\paren{f^\new}
  =
  \log_\ell\Reg_{\F_p^\times}\paren{u_f}.
\]
Equation~\ref{eq:comparison-cross-multiplied} implies
\[
  n_bD^2m_f
  \calp_{\hv,\ell}\paren{f_\varphi}
  =
  \log_\ell\Reg_{\F_p^\times}
  \paren{b_\varphi^\vee a_\varphi D^2u_f},
\]
with the unit on the right belonging to
$\calu_f\paren{\Ad(\rho_f)}$.

Conversely, suppose that
\[
  m_\varphi\calp_{\hv,\ell}\paren{f_\varphi}
  =
  \log_\ell\Reg_{\F_p^\times}\paren{u_\varphi}
\]
for some
$u_\varphi \in \calu_f\paren{\Ad(\rho_f)}$.
Then Equation~\ref{eq:comparison-cross-multiplied} gives
\[
  n_aD^2m_\varphi
  \calp_{\hv,\ell}\paren{f^\new}
  =
  \log_\ell\Reg_{\F_p^\times}
  \paren{a_\varphi^\vee b_\varphi D^2u_\varphi}.
\]
Let $\calo_{F_f}$ be the ring of integers of $F_f$, and
choose a positive integer $c_f$ such that
$c_f\calo_{F_f}\subset R_f$. Define
\[
  u_f' :=
  c_f\Tr_{L_f/F_f}
  \paren{a_\varphi^\vee b_\varphi D^2u_\varphi}.
\]
Then $u_f' \in \calu\paren{\Ad(\rho_f)}$, and
\[
  c_f\sbrac{L_f:F_f}n_aD^2m_\varphi
  \calp_{\hv,\ell}\paren{f^\new}
  =
  \log_\ell\Reg_{\F_p^\times}\paren{u_f'}.
\]
This proves the Harris--Venkatesh conjecture for
$f^\new$.

%%%%%%%%%%%%%%%%%%%%%%%%%%%%%%%%%%%%%%%%%%%%%%%%%%%%%%%%%%%%%%%%%%%%%%%%%

\newpage

\part{Optimal forms}
\label{part:opt}

\section{Global theta liftings}
\label{sec:global-theta}

\subsection{Weil representations and general theta series}
Let $(V, Q)$ be an orthogonal quadratic space over $\Q$ of even
dimension $m$ with bilinear form
\[
	\pair{x, y} = Q(x + y) - Q(x) - Q(y).
\]
Let $\go(V)$ denote the group of similitudes on $V$ with
similitude character
\[
  \nu: \go(V) \lra \G_m.
\]
Since $m$ is even, the group of proper similitudes is
\[
  \gso(V)
  :=
  \set{h \in \go(V) \Mid
    \det(h) = \nu(h)^{m/2}}.
\]
It is the identity component of $\go(V)$. Let
\[
  G := \GL_2 \times_{\G_m} \go(V)
\]
be the fiber product of $\nu$ and
$\det: \GL_2 \lra \G_m$.
With the special elements of $\GL_2$,
\begin{align*}
  d(a) &:= \begin{pmatrix} 1 & \\ & a \end{pmatrix}, \\
  m(a) &:= \begin{pmatrix} a & \\ & a^{-1} \end{pmatrix}, \\
  n(b) &:= \begin{pmatrix} 1 & b\\ & 1 \end{pmatrix}, \\
  w &:= \begin{pmatrix} & 1 \\ -1 & \end{pmatrix},
\end{align*}
$G$ is generated by the elements $(d(\nu(h)), h)$ for $h\in \go(V)$,
$m(a)$, $n(b)$, and $w$.
Then we can consider $\SL_2$ and $\rmo(V)$ as subgroups of $G$.
Let $\cals(V_\A)$ be the space of Schwartz functions on
$V_\A$, and let
\[
  \psi: \A/\Q \lra \C^\times
\]
be the standard character. Following Harris--Kudla
\cite[Section~5]{harris-kudla-1992} and
\cite[Section~3]{harris-kudla-2004}, we use the
realization of the similitude Weil representation on
$\cals(V_\A)$ for the fiber product $G$.\footnote{Waldspurger \cite[Section~I]{waldspurger} and
Yuan--Zhang--Zhang \cite[Section~2.1]{yzz}
instead use a realization for the full product
$\GL_2(\A) \times \go(V_\A)$ on
$\cals(V_\A \times \A^\times)$.}
Thus $G(\A)$ acts on $\cals(V_\A)$ through a Weil
representation $r = r_\psi$, given by the following
rules:
\begin{enumerate}
	\item For any $h \in \go(V_\A)$, $\Phi \in \cals(V_\A)$,
		\[
			r\big(d(\nu(h)), h\big) \cdot \Phi(x) =
				\verts{\nu(h)}^{\frac{-m}{4}} \Phi\paren{h^{-1}x}.
		\]
	\item For any $a \in \A^\times$,
		\[
			r\big(m(a)\big) \cdot \Phi(x) = \eta_V(a) \verts{a}^{m/2} \Phi(ax),
		\]
    where
    \[
      \eta_V(a)
      =
      \paren{a,(-1)^{m/2}\det(V)}_\A.
    \]
    Equivalently, $\eta_V$ is the quadratic character
    attached to
    \[
      \Q\paren{\sqrt{(-1)^{m/2}\det(V)}}/\Q.
    \]
	\item For any $b \in \A$,
		\[
			r\big(n(b)\big) \cdot \Phi(x) = \psi\big(bQ(x)\big) \Phi(x).
		\]
	\item For $w$ as above,
		\[
			r(w) \cdot \Phi(x) = \gamma \cdot \wh \Phi(x),
		\]
		where $\gamma$ is an $8$-th root of unity and $\wh{\Phi}$ is the
		Fourier transform,
		\[
			\wh{\Phi}(x) = \int_{V_\A} \Phi(y) \psi\paren{\pair{x, y}} dy.
		\]
\end{enumerate}

\noindent
From the definition, $r(z,z)$ acts on $\cals(V_\A)$
through the character $\eta_V$. Indeed,
\[
  r(z,z)\Phi(x)
    = r\paren{d(z^2)m(z),z}\Phi(x)
    = \verts{z}^{-m/2}
      r\paren{m(z)}\Phi\paren{z^{-1}x}
    = \eta_V(z)\Phi(x).
\]

For $R \in \set{\Q, \A}$, set
\[
  \GL_2^+(R)
  :=
  \set{g \in \GL_2(R) \Mid
    \det(g) \in \nu\paren{\go(V_R)}}.
\]
For any $\Phi \in \cals(V_\A)$, define the theta series
(or theta kernel)
\[
  \theta(g, h, \Phi)
    :=
    \sum_{x \in V} r(g, h)\Phi(x)
    \in \cala\paren{G(\A)}.
\]

For an automorphic distribution $\varphi$ on
$[\GO(V)]$, define, for $g \in \GL_2^+(\A)$,
\begin{equation}
\label{eq:theta-phi}
  \theta(g,\varphi,\Phi)
  :=
  \int_{[\rmo(V)]}
    \theta(g,hh_0,\Phi)\varphi(hh_0)\,dh,
\end{equation}
where $h_0 \in \GO(V_\A)$ satisfies
$\nu(h_0)=\det(g)$. Extend this function to
$\GL_2(\A)$ by left $\GL_2(\Q)$-invariance and by
setting it equal to zero outside
$\GL_2(\Q)\GL_2^+(\A)$.

Now suppose that there is a character $\omega$ of $\Q^\times \bs
\A^\times$ such that for $z \in \A^\times, h\in \go(V_\A)$,
\[
	\varphi(zh) = \omega(z) \varphi (h).
\]
Then
\begin{align*}
  \theta(zg,\varphi,\Phi)
  &= \int_{[\rmo(V)]}
    \theta(zg,zhh_0,\Phi)
    \varphi(zhh_0)\,dh \\
  &= \eta_V(z)\omega(z)
  \int_{[\rmo(V)]}
    \theta(g,hh_0,\Phi)
    \varphi(hh_0)\,dh \\
  &= \eta_V(z)\omega(z)
  \theta(g,\varphi,\Phi).
\end{align*}

\subsubsection*{Whittaker functions}

In the following, we compute the Whittaker function of
$\theta (g,\varphi, \Phi)$ when $\varphi$ is an automorphic function
on $[\gso (V)]$:
\[
	W(g, \varphi, \Phi) :=
		\int_{\Q \bs \A} \theta(n(b)g, \varphi, \Phi) \psi(-b) db.
\]

\begin{prop}
\label{prop:whittaker}
	The function $W(g, \varphi, \Phi)$ is supported on 
	\[
		\GL_2(\A)_{Q(V_\A)} :=
			\set{g\in \GL_2(\A) \Mid \det g \in Q (V_\A)}.
	\]
	Moreover for $g \in \GL_2(\A)_{Q(V_\A)}$ with decomposition
	$g = d(Q(v)^{-1}) g_1$, where $v \in V_\A$ and $g_1 \in \SL_2(\A)$,
	we have the following expression:
	\[
		W(g, \varphi, \Phi) = \verts{\det g}^{-\frac{m}{4}}
			\int_{\rmo\paren{V_\A} / \rmo\paren{V_{v,\A}}}
			r(g_1) \Phi (h v)
			\int_{\sbrac{\rmo(V_0)}} \varphi\paren{u h_0^{-1} h^{-1}} du dh,
	\]
	where
	\begin{enumerate}
		\item $V_{v, \A}$ is the orthogonal complement of $v$ in $V_{\A}$;
		\item $h_0 \in \go (V_{\A})$ such that $v_0 := h_0^{-1} v\in V$,
			which induces an isomorphism 
			\begin{align*}
				\rmo(V_\A) / \rmo\paren{V_{v,\A}} &\iso \rmo(V_\A) / \rmo(V_{0,\A}) \\
					h &\longmapsto h_0 h h_0^{-1},
			\end{align*}
			where $V_0$ is the orthogonal complement of $v_0$ in $V$;
		\item $dh$ is a measure induced by the above isomorphism and
			the quotient measure of the measure on $\rmo(V_\A)$
			by the measure on $\rmo(V_{0, \A})$ so that 
			the volume of $[\rmo(V_0)]$ is $1$.
	\end{enumerate}
\end{prop}

\begin{proof}
	It is clear that the function $W(g, \varphi, \Phi)$ is also
	supported on $\GL_2(\Q) \cdot \GL_2^+(\A)$.
	For $g \in \GL_2(\Q) \cdot \GL_2^+(\A)$, we may write
	$g = d(a \nu(h_0)) g_1$ for some $a \in \Q^\times$,
	$h_0 \in \go(V_\A)$, and $g_1 \in \SL_2(\A)$.
	Then we have 
	\begin{align*}
		W(g, \varphi, \Phi) &= \int_{\Q \bs \A}
				\theta\paren{n(b) d\big(a\nu(h_0) \big) g_1, \varphi, \Phi}
				\psi(-b) db \\
			&= \int_{\Q \bs \A}
				\theta\paren{n(ab) d\big(\nu(h_0) \big) g_1, \varphi, \Phi}
				\psi(-b) db \\
			&= \int_{\Q \bs \A}
				\theta\paren{n(b) d\big(\nu(h_0)\big) g_1, \varphi, \Phi}
				\psi(-a^{-1}b) db \\
			&= \int_{[\rmo(V)]} \varphi(hh_0) \int_{\Q \bs \A}
				\theta\paren{n(b)d\big(\nu(h_0)) g_1, hh_0, \Phi}
				\psi(-a^{-1}b) db dh.
	\end{align*}
	The second integral can be computed directly
	(for general $g'$ and $h'$):
	\begin{align*}
		\int_{\Q \bs \A} \theta\paren{n(b) g', h', \Phi} \psi(-a^{-1}b) db
			&= \int_{\Q \bs \A} \sum_{x\in V}
				\psi\paren{\paren{Q(x)-a^{-1}}b} r(g', h') \Phi(x) db \\
			&= \sum_{x\in V_a} r(g', h') \Phi(x),
	\end{align*}
	where $V_a$ denotes the subset of elements $x \in V$ with norm
	$Q(x) = a^{-1}$. Define
	\[
		\theta_a(g', h', \Phi) := \sum_{x\in V_a} r(g', h') \Phi(x).
	\]
	Using $g' = d\big(\nu(h_0)\big) g_1$
	and $h' = h h_0$,
	we have shown that 
	\[
		W(g, \varphi, \Phi) = \int_{[\rmo(V)]}
			\theta_a \paren{d\paren{\nu(h_0)} g_1, hh_0, \Phi}
			\varphi\paren{h h_0} dh.
	\]
	This shows that $W(g, \varphi, \Phi)$ is actually supported on
	$Q(V^\times) \GL_2^+(\A)$, where $V^\times$ is the subset of
	elements in $V$ with non-zero norm. This proves the first part of
	Proposition~\ref{prop:whittaker}. 

	For the second part of the Proposition \ref{prop:whittaker}, we use
	the fact that $V_a$ is an orbit of some $v_0 \in V_a$. 
	Let $V_0$ be the orthogonal complement of $v_0$ in $V$.
	Then we have
	\[
  \theta_a\paren{
    d\paren{\nu(h_0)}g_1,
    hh_0,
    \Phi}
  = \sum_{\gamma \in \rmo(V_0)\bs\rmo(V)}
    \verts{\det g}^{-\frac{m}{4}}
    r(g_1)
    \Phi\paren{h_0^{-1}h^{-1}\gamma^{-1}v_0}.
\]
	It follows that
	\begin{align*}
		W(g, \varphi, \Phi) &= \verts{\det g}^{-\frac{m}{4}}
				\int_{\rmo\paren{V_{0}} \bs \rmo \paren{V_\A}}
				r\paren{g_1} \Phi\paren{h_0^{-1} h^{-1} v_0}
				\varphi\paren{hh_0} dh \\
			&= \verts{\det g}^{-\frac{m}{4}}
				\int_{\rmo\paren{V_{0, \A}} \bs \rmo(V_\A)}
				r\paren{g_1} \Phi\paren{h_0^{-1}h^{-1}v_0}
				\int_{\sbrac{\rmo\paren{V_0}}} \varphi\paren{uhh_0} du dh.\\
	\end{align*}

	A change of variables $h \mapsto h_0 h^{-1} h_0^{-1}$ yields
	\[
		W(g, \varphi, \Phi) = \verts{\det g}^{-\frac{m}{4}}
			\int_{\rmo\paren{V_\A} / h_0 \rmo\paren{V_{0, \A}} h_0^{-1}}
			r\paren{g_1} \Phi\paren{h h_0^{-1} v_0}
			\int_{\sbrac{\rmo\paren{V_0}}} \varphi\paren{uh_0h^{-1}} du dh.
	\]
	Set $v = h_0^{-1}v_0$. Then $Q(v) = \nu(h_0)^{-1} a^{-1}$.
	Finally, change $h_0$ to $h_0^{-1}$.
\end{proof}
 
It is useful to pass to the Kirillov model,
obtained by restricting Whittaker functions to
$a(x) := \begin{psmallmatrix}x&0\\0&1\end{psmallmatrix}$.
If $x=Q(v)$, then $a(x)=x\,d(x^{-1})$. Applying the
central character $\omega$ gives the following.

\begin{cor}
\label{cor:Kirillov}
	Assume that $\varphi$ has the central character $\omega$.
	Then the Kirillov function $\kappa (x, \varphi, \Phi)$ for the
	theta series $\theta (g, \varphi, \Phi)$
	is supported on $Q(V_\A)$ with the following formula 
	\[
		\kappa(x, \varphi, \Phi) = 
			\eta_V(x) \omega(x) \verts{x}^{\frac{m}{4}}
			\int_{\rmo\paren{V_\A} / \rmo\paren{V_{v, \A}}}
			\Phi(hv) \int_{\sbrac{\rmo(V_0)}}
			\varphi \paren{u h_0^{-1} h^{-1}} du dh,
	\]
	where $v \in V_{\A}$ and $h_0 \in \go(V_\A)$ such that $Q(v)=x$,
	$v_0 = h_0^{-1} v \in V$,
	and $V_0$ is the orthogonal complement of $v_0$.
\end{cor}

\subsection{Theta series for one character}
\label{sec:theta}
Let $K$ be a quadratic field and
$\chi: K^\times \bs K_\A^\times \lra \C^\times$ be a finite character.
Assume the following conditions:
\begin{enumerate}
  \item $\chi$ does not factor through the norm map
    $\rmn_{K/\Q}$ of $K$ over $\Q$; equivalently,
    $\chi$ is not of the form
    $\mu \circ \rmn_{K/\Q}$ for a character
    $\mu: \Q^\times \bs \A^\times \lra \C^\times$;
  \item if $K$ is real, then the two components at the
    archimedean places have different signs.
\end{enumerate}
Then we have an irreducible cuspidal representation $\pi(\chi)$ of
$\GL_2$ of weight $1$.
In the following, we want to construct newforms in $\pi(\chi)$ and
optimal forms in $\pi(\chi) \otimes \pi(\chi^{-1})$ using theta
liftings.
  
We first take the norm quadratic space
$V = (K,\rmn_{K/\Q})$.
Let $c_K$ denote the nontrivial element of
$\Gal(K/\Q)$. For a Hecke character $\chi$ of $K$,
write $\chi^{c_K}:=\chi\circ c_K$.
Multiplication identifies $\gso(V)$ with
$\operatorname{Res}_{K/\Q}\G_m$, and the similitude
character is $\rmn_{K/\Q}$. For each
$\Phi\in\cals(K_\A)$, the theta lift of $\chi^{c_K}$ is an
automorphic form
\[
  \theta\paren{g,\chi^{c_K},\Phi}
  \in
  \cala\paren{\GL_2(\Q)\bs\GL_2(\A)}.
\]
Its Whittaker function is supported on $\GL_2^+(\A)$.
For $g \in \GL_2^+(\A)$, choose
$h_0 \in K_\A^\times$ and $g_1 \in \SL_2(\A)$, set
$x := \rmn_{K/\Q}(h_0)$, and write
$g = d(x^{-1}) g_1$.
Proposition~\ref{prop:whittaker} gives
\begin{equation}
\label{eq:w-chi-phi}
  W\paren{g,\chi^{c_K},\Phi}
  =
  \verts{\det g}^{-\frac{1}{2}}
  \int_{K_\A^1}
    r\paren{g_1}\Phi\paren{hh_0}
    \chi^{c_K}\paren{h_0^{-1}h^{-1}}\,dh,
\end{equation}
where
$K_\A^1:=\ker\paren{\rmn_{K/\Q}:
K_\A^\times\longrightarrow\A^\times}$.

We now explain the corresponding Kirillov formula. The
restriction of $\chi^{c_K}$ to $\A^\times$ is the same as
that of $\chi$, so the central character of the theta
series is $\eta_V\restr{\chi}{\A^\times}$.
For $x = \rmn_{K/\Q}(h_0)$,
we have $\eta_V(x) = 1$ and $a(x) = x\,d(x^{-1})$.
Moreover, for $h\in K_\A^1$,
\[
  \chi(x)\chi^{c_K}\paren{h_0^{-1}h^{-1}}
  =
  \chi\paren{hh_0}.
\]
Indeed, $x=h_0\overline{h}_0$ and
$\overline{h} = h^{-1}$. Equation~\ref{eq:w-chi-phi}
therefore gives
\begin{equation}
\label{eq:kappa-chi-phi}
  \kappa\paren{x,\chi^{c_K},\Phi}
  =
  \verts{x}^{\frac{1}{2}}
  \int_{K_\A^1}
    \Phi\paren{h_0h}\chi\paren{hh_0}\,dh.
\end{equation}

Choose the Haar measures on $K_\A^\times$,
$K_\A^1$, and $\rmn_{K/\Q}(K_\A^\times)$ compatibly.
Unfolding Equation~\ref{eq:kappa-chi-phi} gives the Tate
zeta integral
\begin{align*}
  Z\paren{\theta(g,\chi^{c_K},\Phi),s}
  &:=
  \int_{\Q^\times\bs\A^\times}
    \theta\paren{a(x),\chi^{c_K},\Phi}
    \verts{x}^{s-\frac{1}{2}}\,d^\times x \\
  &=
  \int_{\A^\times}
    \kappa\paren{x,\chi^{c_K},\Phi}
    \verts{x}^{s-\frac{1}{2}}\,d^\times x \\
  &=
  \int_{K_\A^\times}
    \Phi(t)\chi(t)
    \verts{\rmn_{K/\Q}(t)}^s\,d^\times t \\
  &=:
  Z\paren{\chi,\Phi,s}.
\end{align*}

As $\Phi$ varies, the theta series
$\theta(g,\chi^{c_K},\Phi)$ generate an irreducible
representation denoted by $\pi(\chi)$. This
representation decomposes as
(cf. \cite[Equation~5.1]{shimizu})
\[
  \pi(\chi)\iso
  \bigotimes_{v\leq\infty}'\pi(\chi_v).
\]
For each place $v$, put
$K_v^1 := \ker\paren{\rmn_{K_v/\Q_v}:
K_v^\times\longrightarrow\Q_v^\times}$. Choose
$h_0 \in K_v^\times$, set
$x := \rmn_{K_v/\Q_v}(h_0)$, and write
$g = d(x^{-1}) g_1$ with $g_1 \in \SL_2(\Q_v)$.
The Whittaker and Kirillov models of $\pi(\chi_v)$ are
generated by the following functions
(cf. \cite[Equation~5.2]{shimizu}):
\begin{align}
  W\paren{g,\chi_v^{c_K},\Phi_v}
  &=
  \verts{\det g}^{-\frac{1}{2}}
  \int_{K_v^1}
    r\paren{g_1}\Phi_v\paren{hh_0}
    \chi_v^{c_K}\paren{h_0^{-1}h^{-1}}\,dh,
  \label{eq:wp-chi-phi} \\
  \kappa\paren{x,\chi_v^{c_K},\Phi_v}
  &=
  \verts{x}^{\frac{1}{2}}
  \int_{K_v^1}
    \Phi_v\paren{h_0h}\chi_v\paren{hh_0}\,dh.
  \label{eq:kappap-chi-phi}
\end{align}

\subsubsection*{Newforms}
Now assume $V = (K, \rmn)$, where $\rmn = \rmn_{K/\Q}$ is the norm of $K$
over $\Q$.
We construct a newform $\varphi^\new \in \pi(\chi)$
by picking a ``standard'' Schwartz function for $\chi$,
\[
	\Phi_\chi = \bigotimes_v \Phi_{\chi_v} \in \cals(\A_K),
\]
where the tensor product is over places of $K$. We pick
$\Phi_{\chi_v}$ as follows (cf. \cite[Section~2.1]{zhang-asian}):
\begin{enumerate}
	\item If $v$ is complex, $K_v \iso \C$, and $\chi_v$ is trivial,
		take
		\[
			\Phi_{\chi_v}(x + yi) = e^{-2\pi\paren{x^2 + y^2}}.
		\]
	\item If $v$ is real, $K_v = \R$, and $\chi_v(x) = \sgn(x)^m$ with
		$m = 0, 1$, take
		\[
			\Phi_{\chi_v}(x) = x^m e^{-\pi x^2}.
		\]
	\item If $v$ is finite and $\chi_v$ is unramified, take
		\[
			\Phi_{\chi_v} = \restr{\one}{\calo_{K_v}}.
		\]
	\item If $v$ is finite and $\chi_v$ is ramified, take 
		\[
			\Phi_{\chi_v} = \restr{\chi_v^{-1}}{\calo_{K_v}^\times}.
		\]
\end{enumerate}

For each place $p$ of $\Q$, let
\[
	\Phi_{\chi_p} = \bigotimes_{v \, \mid \, p} \Phi_{\chi_v}
		\in \cals(K_p).
\]
In the real case, to avoid the need to remember signs,
we use the following function instead:
\[
	\Phi_\infty(x, y) = \frac{1}{2}(x + y)e^{-\pi\paren{x^2 + y^2}}.
\]
Using the ``standard'' $\Phi_\chi$ as defined above,
we have the theta series $\theta(g, \chi^{c_K}, \Phi_\chi)$
and the following local description of the Whittaker function
$W(g, \chi^{c_K}, \Phi_{\chi})$ at each place of $\Q$
(cf. \cite[Section~2.3]{zhang-asian}).
\begin{enumerate}
   \item If $p = \infty$, then
    $W(g,\chi_\infty^{c_K},\Phi_{\chi_\infty})$ is the
    standard holomorphic Whittaker function of weight
    $1$. Writing
    $k_\theta
      :=
      \begin{psmallmatrix}
        \cos\theta & \sin\theta \\
        -\sin\theta & \cos\theta
      \end{psmallmatrix}$,
    one has
    \[
      W\paren{
        r
        \begin{pmatrix}y & x \\ 0 & 1\end{pmatrix}
        k_\theta,
        \chi_\infty^{c_K},
        \Phi_{\chi_\infty}}
      =
      \begin{cases}
        \sgn(r)^m y^{1/2}
        e^{2\pi i(x + iy)}e^{i\theta}
          & \text{if $y > 0$,} \\
        0 & \text{if $y < 0$,}
      \end{cases}
    \]
    for $r \in \R^\times$ and $x,\theta \in \R$,
    where $m = 0$ if $K$ is imaginary and $m = 1$
    if $K$ is real.
  \item If $p$ is not ramified in $K$, then
    $W(g, \chi_p^{c_K}, \Phi_{\chi_p})$ is the
    unique Whittaker newform
    $W^\new_\chi$ in $\pi(\chi_p)$ that
    is invariant under $U_1\paren{p^{c(\pi(\chi_p))}}$
    and takes value $1$ at the identity $e$.
  \item If $p$ is ramified in $K$, then
    $W(g, \chi_p^{c_K}, \Phi_{\chi_p})$ is the
    restriction of a newform
    $W_\chi^\new $ on $\GL_2(\Q_p)^+$.
    One can also recover the newform by taking
    \[
      W^\new_\chi(g)
      :=
      \sum_{
        \beta_p \in
        \Z_p^\times/
        \rmn(\calo_{K_p}^\times)}
      W\paren{
        g a(\beta_p),
        \chi_p^{c_K},
        \Phi_{\chi_p}}.
    \]
\end{enumerate}
The newform is given by
\begin{equation}
  \label{eq:theta}
  \varphi^\new_\chi(g) =
    \sum_{\beta \in \wh\Z^\times / \rmn(\wh\calo_K^\times)}
    \theta\paren{g a(\beta), \chi^{c_K}, \Phi_\chi}.
\end{equation}
From the local description above,
the multiplicative local factor at archimedean and unramified finite
places is a single term with $\beta_p = 1$,
while the multiplicative local factor at ramified finite places
is the sum of two terms
since $\#(\Z_p^\times / \rmn(\calo_{K_p}^\times)) = 2$ if $p$ is ramified.
Hence the sum in the right-hand side has $2^n$ terms,
where $n$ is the number of ramified primes in $K$.

\subsubsection*{Comparison of models}

Let $V = (Ke, Q)$ and $V' = (Ke', Q')$ be
quadratic spaces with an action of $K$, and let
\[
  \iota: V'_\A \iso V_\A
\]
be an isomorphism of $K_\A$-modules. There is a unique
$a_\iota \in K_\A^\times$ such that
$\iota(e') = a_\iota e$.
Set
\[
  Q(\iota)
  :=
  \frac{Q\paren{\iota(e')}}{Q'(e')}
  \in \A^\times.
\]
Pullback along $\iota$ defines an isomorphism
\begin{align*}
  \iota^*: \cals(V_\A)
    &\iso \cals(V'_\A), \\
  \Phi
    &\longmapsto \Phi \circ \iota.
\end{align*}

We compare the corresponding Kirillov functions. Let
\[
  x = Q(t_0e),
  \qquad
  t_0' := t_0a_\iota^{-1}.
\]
Then
\[
  Q'(t_0'e') = xQ(\iota)^{-1},
\]
and
\[
  \iota^*\Phi(t_0'te') = \Phi(t_0te).
\]
Equation~\ref{eq:kappa-chi-phi} therefore gives
\begin{align*}
  &\kappa_{V'}\paren{
    xQ(\iota)^{-1},
    \chi^{c_K},
    \iota^*\Phi} \\
  &\qquad =
  \verts{xQ(\iota)^{-1}}^{\frac{1}{2}}
  \int_{K_\A^1}
    \Phi(t_0te)
    \chi\paren{tt_0a_\iota^{-1}}\,dt \\
  &\qquad =
  \chi\paren{a_\iota}^{-1}
  \verts{Q(\iota)}^{-\frac{1}{2}}
  \kappa_V\paren{x, \chi^{c_K}, \Phi}.
\end{align*}
Thus
\begin{equation}
\label{eq:model-comparison-kirillov}
  \kappa_V\paren{x, \chi^{c_K}, \Phi}
  =
  \chi\paren{a_\iota}
  \verts{Q(\iota)}^{\frac{1}{2}}
  \kappa_{V'}\paren{
    xQ(\iota)^{-1},
    \chi^{c_K},
    \iota^*\Phi}.
\end{equation}

Although $a_\iota$ is defined using the chosen
generators of the one-dimensional $K$-spaces $V$ and
$V'$, the scalar $\chi(a_\iota)$ and the idele
$Q(\iota)$ are independent of those choices. Indeed,
changing either generator multiplies $a_\iota$ by an
element of $K^\times$, on which $\chi$ is trivial.
Since Kirillov functions determine automorphic forms, we
obtain the following comparison.

\begin{prop}
\label{prop:theta-comparison}
  Let $V = (Ke,Q)$ and $V' = (Ke',Q')$ be
  two-dimensional quadratic spaces with an action of
  $K$, and let
  \[
    \iota:V'_\A \iso V_\A
  \]
  be an isomorphism of $K_\A$-modules. Define
  $a_\iota \in K_\A^\times$ and
  $Q(\iota) \in \A^\times$ by
  \[
    \iota(e') = a_\iota e,
    \qquad
    Q(\iota)
    =
    \frac{Q\paren{\iota(e')}}{Q'(e')}.
  \]
  Then, for every $\Phi \in \cals(V_\A)$,
  \[
    \theta_V\paren{g,\chi^{c_K},\Phi}
    =
    \chi\paren{a_\iota}
    \verts{Q(\iota)}^{\frac{1}{2}}
    \theta_{V'}\paren{
      g a\paren{Q(\iota)^{-1}},
      \chi^{c_K},
      \iota^*\Phi}.
  \]
\end{prop}

For example, if we compare the theta functions defined by two
opposite spaces $V_{\pm} := (V, \pm Q)$, then 
for each $\Phi \in \cals(V_\A)$, we get two theta series:
$\theta_\pm (g, \chi^{c_K}, \Phi)$.
We use the identity map $\iota: V_\A \lra V_\A $ for the quadratic
space, so $Q(\iota) = -1$.
Then
\[
  \theta_{-}\paren{g,\chi^{c_K},\Phi}
  =
  \theta_{+}\paren{ga(-1),\chi^{c_K},\Phi}.
\]

In the case that $V^\pm = (K, \pm\rmn)$ and
$\Phi = \Phi_\chi$, Proposition~\ref{prop:theta-comparison}
shows that the Whittaker function
$W_-(g, \chi^{c_K}, \Phi_\chi)$ of
$\theta_-(g, \chi^{c_K}, \Phi_\chi)$ is still new at every
finite place, but has weight $-1$ at infinity. Its value
on diagonal elements is
\[
  W_-\paren{a(y)}
    =
    \begin{cases}
      \verts{y}^{1/2}e^{2\pi y}
        & \text{if $y<0$,} \\
      0 & \text{if $y>0$.}
    \end{cases}
\]
Let $(-1)_\infty \in \A^\times$ be the idele with
archimedean component $-1$ and finite components $1$.
Then
\[
  \theta_{-}\paren{g,\chi^{c_K},\Phi}
  =
  \theta_{+}\paren{
    ga\paren{(-1)_\infty},
    \chi^{c_K},
    \Phi}.
\]

\subsection{Theta series for two characters}
\label{sec:theta-two}
\subsubsection*{Theta series for automorphic forms}
Now we consider the theta lifting for $V = (B, \rmn)$, with $B$ a
quaternion algebra over $\Q$ and with norm $\rmn$ given by the reduced
norm on $B$. Then
\begin{align*}
	\GSO(V) &\cong (B^\times\times B^\times)/\Delta\Q^\times, \\
  \GO(V) &= \abrac{\GSO(V),\iota}.
\end{align*}
where $(b_1, b_2) \in B^\times \times B^\times$ brings
$x \in V$ to $b_1 x b_2^{-1}$, and $\iota(x) = \overline{x}$.
Let $G$ denote the group over $\Q$ defined by
\[
	G := \GL_2 \times_{\G_m} \gso (V).
\]
Then we have a Weil representation of $G(\A)$ on $\cals(V(\A))$.
For each $\Phi \in \cals (V(\A))$, we have a theta series
\[
	\theta(g, h, \Phi) = \sum_{x \in V} r(g, h) \Phi(x).
\]
Also, for each automorphic form, or more generally each
automorphic distribution, $\varphi$ on $\go(V)$, we
obtain a form on $\GL_2^+(\A)$:
\[
  \theta(g, \varphi, \Phi)
  =
  \int_{[O(V)]}
    \theta\paren{g, h h_0, \Phi}
    \varphi\paren{h h_0}\,dh.
\]
We want to interpret these theta liftings as Hecke
operators. For any $g \in B_\A^\times$, define the
right-translation operator $\rho(g)$ on
$\cala([B^\times])$ by
\[
  \rho(g)\varphi(x) := \varphi(xg).
\]
For any $x \in \rmn(B_\A^\times)$
and $\Phi \in \cals(B_\A)$, define
\begin{align}
\label{eq:hecke}
  \rmt_\Phi(x)
  &:=
  \int_{B_\A^1}
    \Phi\paren{b_0b}
    \rho\paren{b_0b}\,db, \\
  \rmt_\Phi^*(x)
  &:=
  \int_{B_\A^1}
    \Phi\paren{b^{-1}b_0^{-1}}
    \rho\paren{b_0b}\,db,
\end{align}
where $b_0 \in B_\A^\times$ satisfies
$\rmn(b_0) = x$.

\begin{prop}
\label{prop:hecke-auto}
	Let $\varphi = \varphi_1 \otimes \varphi_2$ with $\varphi_i$
	automorphic forms on $[B^\times]$ with central characters $\omega$
	and $\omega^{-1}$.
	Then the Kirillov function
  $\kappa(x, \varphi, \Phi)$ is supported on
  $\rmn(B_\A^\times)$ and is given by
  \[
    \kappa(x, \varphi, \Phi)
    = \omega(x)\verts{x}
      \pair{\varphi_1, \rmt_\Phi(x)\varphi_2}
    = \omega(x)\verts{x}
      \pair{\rmt_\Phi^*(x)\varphi_1, \varphi_2},
  \]
  where $\pair{\cdot, \cdot}$ denotes the bilinear form
  \[
    \pair{\varphi_1, \varphi_2}
    = \int_{\sbrac{B^\times / \Q^\times}}
      \varphi_1(u)\varphi_2(u)\,du.
  \]
\end{prop}
	
\begin{proof}
	By Corollary~\ref{cor:Kirillov}, take
  $x = \rmn(b_0)$
  for some $b_0 \in B_\A^\times$,
  together with
  $h_0 = (1, b_0^{-1})$ and $v_0 = e$.
  Then
	\[
		\kappa (x, \varphi, \Phi) = \omega(x) \verts{x}
			\int_{\rmo(V_\A) / \rmo(V_{b_0, \A})} \Phi\paren{h b_0}
			\int_{\sbrac{\rmo(V_0)}}
			\varphi\paren{u \cdot (1, b_0) \cdot h^{-1}} du dh.
	\]
	Here,
	$\rmo(V) = B^\times \times_{\Q^\times} B^\times / \Delta \Q^\times$,
	and $\rmo(V_{b_0, \A})$ consists of 
	elements of the form $(b_0 b b_0^{-1}, b)$ for all
	$b \in B_\A^\times$.
	So we can use elements $(1, b^{-1})$ for $b \in B_\A^1 $
	to represent quotient elements.
	Then the above integral becomes:
	\begin{align*}
		\kappa(x, \varphi, \Phi)
			&= \omega(x) \, \verts{x} \int_{B_\A^1} \Phi\paren{b_0 b}
				\int_{\sbrac{B^\times / \Q^\times}}
				\varphi\paren{u, u b_0 b} du db \\
			&= \omega(x) \, \verts{x} \int_{B_\A^1} \Phi\paren{b_0 b}
				\int_{\sbrac{B^\times / \Q^\times}}
				\varphi_1(u) \varphi_2\paren{u b_0 b} du db \\
			&= \omega(x) \, \verts{x} \int_{B_\A^1} \Phi\paren{b_0 b}
				\pair{\varphi_1, \rho\paren{b_0 b} \varphi_2} db\\
			&= \omega(x) \, \verts{x} \int_{B_\A^1} \Phi\paren{b_0 b}
				\pair{\rho \paren{b^{-1}b_0^{-1}} \varphi_1, \varphi_2} db\\
	\end{align*}
	The proposition follows from the last two identities. 
\end{proof}

\subsubsection*{Theta series for two characters}
Now let $K$ be a quadratic field embedded into $B$. Then we have a
decomposition $B = K + Kj$, which gives an orthogonal decomposition
$V = V_1 + V_2$ for $V = (B, \rmn)$. Then we have an embedding 
\[
	\go(V_1) \times_{\G_m} \go(V_2) \subset \go(V).
\]
The restriction to the connected component can be described as:
\[
	\begin{tikzcd}
		T := K^\times \times_{\Q^\times} K^\times / \Delta\paren{\Q^\times}
				\arrow[r, leftarrow, "\widesim{}"]
			& K^\times \times K^\times/\Delta (\Q^\times)
				\arrow[r, hookrightarrow]
			& B^\times \times B^\times / \Delta\paren{\Q^\times}.
	\end{tikzcd}
\]
The first map is induced by
  $(t_1, t_2)
    \mapsto
    \paren{
      t_1/t_2,
      t_1/\overline{t}_2}$.
It is surjective by Hilbert's theorem 90. Indeed, if
$(a, b) \in T$, then
$\rmn_{K/\Q}(b/a) = 1$,
so there is a $t_2 \in K^\times$ such that
$\frac{b}{a} = t_2/\overline{t}_2$.
Taking $t_1 := at_2$ gives
$a = t_1/t_2$ and $b = t_1/\overline{t}_2$.
The same argument place by place gives the corresponding
surjectivity on adelic points
(cf. Serre~\cite[Chapter~X,
Proposition~3]{serre-local-fields}).

There are two equivalent descriptions of an automorphic
character of $T$. One may use two automorphic characters
$\xi_1$ and $\xi_2$ of $\A_K^\times$ whose restrictions
to $\A^\times$ are inverse to one another, or the
restriction to
$[K^\times \times_{\Q^\times} K^\times]$
of a character $\chi_1 \otimes \chi_2$ on
$[K^\times \times K^\times]$.
Using the conjugation notation fixed above, the two
descriptions are related by
\begin{align*}
  \chi_1\paren{t_1/t_2}
  \chi_2\paren{t_1/\overline{t}_2}
  &=
  \xi_1(t_1)\xi_2(t_2), \\
  \xi_1
  &=
  \chi_1\chi_2, \\
  \xi_2
  &=
  \chi_1^{-1}\paren{\chi_2^{c_K}}^{-1}.
\end{align*}

Let
$T^1 := \ker\paren{\nu: T \longrightarrow \G_m}$,
and write
$[T^1] := T^1(\Q)\bs T^1(\A)$.
For an automorphic character
$\xi = \xi_1 \otimes \xi_2$ and a function
$\Phi \in \cals(B_\A)$, define
\begin{equation}
\label{eq:theta-xi}
  \theta(g, \xi, \Phi)
  :=
  \int_{[T^1]}
    \theta\paren{g, t_0t, \Phi}
    \xi\paren{t_0t}\,dt,
\end{equation}
where $t_0 \in T(\A)$ satisfies
$\nu(t_0) = \det(g)$.
This is the specialization of the theta lifting in
Equation~\ref{eq:theta-phi} to the automorphic
distribution on $[\gso(V)]$ induced from $T$.

Assume that $\Phi = \Phi_1 \otimes \Phi_2 \in \cals(V_{\A}) =
\cals(V_{1, \A}) \otimes \cals(V_{2, \A})$ is a decomposable function.
Then for $h = (h_1, h_2) \in \go(V_1) \times_{\G_m} \go(V_2)$,
we have
\[
	\theta(g, h, \Phi) = \theta\paren{g, h_1, \Phi_1} \cdot
		\theta\paren{g, h_2, \Phi_2}.
\]
Thus, if $\xi_1 \otimes \xi_2$ is the restriction of
$\chi_1 \otimes \chi_2$, then
\begin{equation}
	\label{eq:xi-chi}
	\theta\paren{g, \xi_1\otimes \xi_2, \Phi} =
		\theta\paren{g, \chi_1, \Phi_1} \cdot
		\theta\paren{g, \chi_2, \Phi_2}.
\end{equation}

We now specialize to the case when $K$ is imaginary.

\begin{definition}
\label{def:standard-Eichler}
  Let $K$ be an imaginary quadratic field and let $B$
  be a definite quaternion algebra. Let $\calo$ be an
  Eichler order of $B$. Write $D_B := \disc(B)$,
  and let $N_\calo$ be the Eichler level of $\calo$.
  Thus the reduced discriminant of $\calo$ is
  $M := D_B N_\calo$;
  see \cite[Chapter~23]{voight}. Define the standard
  Schwartz function
  \[
    \Phi_\calo
    :=
    \Phi_{\calo,\infty}
    \otimes
    \Phi_\calo^\infty
    \in \cals(B_\A)
  \]
  by
  \begin{enumerate}
    \item
      $\Phi_{\calo,\infty}(x)
      := e^{-2\pi\verts{x}^2}$;
    \item
      $\Phi_\calo^\infty
      := \one_{\wh\calo}$.
  \end{enumerate}
\end{definition}

In the remainder of this subsection, assume that $M$
is odd; this is the case in all applications below.
Let $\xi_1$ and $\xi_2$ be finite characters with
opposite restrictions to $\A^\times$, and set
\[
  \omega
  :=
  \restr{\xi_1}{\A^\times}
  =
  \restr{\xi_2^{-1}}{\A^\times}.
\]

\begin{prop}
\label{prop:theta-xi}
  The theta series
  \[
    \theta\paren{g, \xi_1 \otimes \xi_2,
      \Phi_\calo}
  \]
  is holomorphic of weight $2$, is right-invariant
  under $U_1(M)$, and has central character $\omega$.
\end{prop}

\begin{proof}
  The archimedean Gaussian factors under the
  orthogonal decomposition $V = V_1 \oplus V_2$.
  At the finite adeles, the space
  $\cals(V(\A^\infty))$ is spanned by decomposable
  Schwartz functions with respect to
  $V=V_1\oplus V_2$. Hence $\Phi_\calo$ is a finite
  sum of decomposable Schwartz functions.
  is an isomorphism. Hence $\Phi_\calo$ is a finite
  sum of decomposable Schwartz functions. Applying
  Equation~\ref{eq:xi-chi} term by term shows that
  the theta series is holomorphic of weight $2$.

  We next check the level. For a finite prime $q$, set
  $m_q := \ord_q(M)$,
  and let $\calo_q^\vee$ be the dual of $\calo_q$
  for the pairing
  $(x, y) \mapsto
  \operatorname{trd} \paren{x\overline{y}}$.
  In the standard split Eichler model, and in the
  maximal-order model when $q \mid D_B$, one computes
  $q^{m_q}\rmn\paren{\calo_q^\vee} \subset \Z_q$.
  The standard local lattice calculation in the
  similitude Weil representation therefore shows that
  the local theta vector attached to
  $\one_{\calo_q}$ is fixed by
  $U_1\paren{q^{m_q}}$
  (cf. \cite[Lemmas~7 and~10]{shimizu}).
  Thus the global
  theta series is right-invariant under $U_1(M)$.

  Since $\eta_V$ is trivial for the reduced norm form
  on $B$, the central-character calculation preceding
  Proposition~\ref{prop:whittaker} gives $\omega$.
\end{proof}

By Proposition~\ref{prop:theta-xi}, the theta series
$\theta(g,\xi,\Phi_\calo)$ is right-invariant under
$U_1(M)$. Its archimedean Whittaker function is the
fixed weight-two function $W_2$, so the theta series is
determined by its finite Kirillov function. We would
like to apply Proposition~\ref{prop:hecke-auto} to
express this function in terms of quaternionic Hecke
operators. We first need to replace the characters
$\xi_1$ and $\xi_2$ by automorphic forms on
$[B^\times]$.

\subsection{Hecke operators}
\label{subsec:hecke-operators}
Generalizing \cite[Section~2.2]{dhrv} from characters
of the ideal class group of an imaginary quadratic $K$
to general Hecke characters, we define a projection map
$[\dotspace]$ from characters of $[K^\times]$ to
quaternionic automorphic forms.

We first specify the quaternionic level. For every prime
$q \nmid D_B$, choose an isomorphism
$\iota_q: B_q \iso M_2(\Q_q)$
that carries $\calo_q$ to the standard Eichler order of
level $q^{\ord_q(N_\calo)}$. Set
\[
  U_{1,q}^B(M)
  :=
  \begin{cases}
    \iota_q^{-1}
    \paren{U_1\paren{q^{\ord_q(N_\calo)}}}
      & \text{if $q \nmid D_B$}, \\
    \calo_q^\times
      & \text{if $q \mid D_B$},
  \end{cases}
\]
and
\begin{align*}
  U_1^B(M)
    &:=\prod_{q < \infty} U_{1,q}^B(M), \\
  U_0^B(M)
    &:= \wh\calo^\times.
\end{align*}
At a prime $q \mid D_B$, the order $\calo_q$ is the
unique maximal order of $B_q$, so the second line means
invariance under its full unit group.

Let $B_\infty^1$ be the norm-one subgroup of
$B_\infty^\times$. For a finite automorphic character
$\omega$ of $\Q^\times \bs \A^\times$, let
$\cala^\pm(\omega)$ be the space of automorphic forms
$\varphi$ on $[B^\times]$ such that
\[
  \varphi\paren{z g u_\infty u}
  =
  \omega^{\pm 1}(z)\varphi(g)
\]
for all
$z \in \A^\times$,
$u_\infty \in B_\infty^1$,
and $u \in U_1^B(M)$.
Since $B$ is definite, these spaces are
finite-dimensional, and the pairing
\[
  \pair{\varphi_+, \varphi_-}
  :=
  \int_{\sbrac{B^\times / \Q^\times}}
    \varphi_+(x)\varphi_-(x)\,dx
\]
identifies $\cala^+(\omega)$ with the dual of
$\cala^-(\omega)$.
Analogously, let $\Xi^\pm(\omega)$ be the space of
characters on $[K^\times]$ that are invariant under
$U_K := U_0^B(M) \cap \wh K^\times$
and whose restrictions to $\A^\times$ are
$\omega^{\pm 1}$.

For every $\xi^\pm \in \Xi^\pm(\omega)$, define
$\sbrac{\xi^\pm} \in \cala^\pm(\omega)$ by the
identity
\begin{equation}
\label{eq:proj-complex}
  \pair{\varphi, \sbrac{\xi^\pm}}
  =
  \int_{\sbrac{K^\times / \Q^\times}}
    \varphi(t)\xi^\pm(t)\,dt
  \qquad
  \paren{\varphi \in \cala^\mp(\omega)}.
\end{equation}
The restrictions of $\varphi$ and $\xi^\pm$ to the
center are inverse, so the integral on the right is
well-defined. Applying Equation~\ref{eq:proj-complex}
in both variables gives, for
$\xi_1 \in \Xi^+(\omega)$ and
$\xi_2 \in \Xi^-(\omega)$,
\begin{equation}
\label{eq:xi-varphi}
  \theta\paren{g, \xi_1 \otimes \xi_2,
    \Phi_\calo}
  =
  \theta\paren{g, \sbrac{\xi_1} \otimes
    \sbrac{\xi_2}, \Phi_\calo}.
\end{equation}

Normalize the Haar measure on $B_\infty^1$ to have
volume $1$. For
$x^\infty \in \rmn\paren{B_{\A^\infty}^\times}$,
define the finite Hecke operator
\begin{equation}
\label{eq:hecke-finite}
  \rmt_{\Phi_\calo^\infty}\paren{x^\infty}
  :=
  \int_{B_{\A^\infty}^1}
    \Phi_\calo^\infty\paren{b_0^\infty b}
    \rho\paren{b_0^\infty b}\,db,
\end{equation}
where $\rmn\paren{b_0^\infty} = x^\infty$.
For a holomorphic weight-two theta series, let
$\kappa^\infty$ denote the finite Kirillov function,
characterized by
\[
  W\paren{a(x), \varphi, \Phi_\calo}
  =
  W_2\paren{a(x_\infty)}
  \kappa^\infty\paren{x^\infty,
    \varphi, \Phi_\calo}
\]
for $x = (x_\infty, x^\infty)$ and $x_\infty > 0$.

\begin{prop}
\label{prop:hecke-complex}
  Let $\omega$ be a finite automorphic character of
  $\Q^\times \bs \A^\times$, and let
  \[
    \varphi_1 \in \cala^+(\omega),
    \qquad
    \varphi_2 \in \cala^-(\omega).
  \]
  Then
  \[
    \theta\paren{g, \varphi_1 \otimes
      \varphi_2, \Phi_\calo}
  \]
  is holomorphic of weight $2$, is right-invariant
  under $U_1(M)$, and has central character $\omega$.
  Moreover,
  \[
    \kappa^\infty\paren{x^\infty,
      \varphi_1 \otimes \varphi_2,
      \Phi_\calo}
    =
    \omega^\infty\paren{x^\infty}
    \verts{x^\infty}
    \pair{\varphi_1,
      \rmt_{\Phi_\calo^\infty}
      \paren{x^\infty}\varphi_2}.
  \]
  In particular, if
  $\xi_1 \in \Xi^+(\omega)$ and
  $\xi_2 \in \Xi^-(\omega)$, then
  \[
    \kappa^\infty\paren{x^\infty,
      \xi_1 \otimes \xi_2,
      \Phi_\calo}
    =
    \omega^\infty\paren{x^\infty}
    \verts{x^\infty}
    \pair{\sbrac{\xi_1},
      \rmt_{\Phi_\calo^\infty}
      \paren{x^\infty}
      \sbrac{\xi_2}}.
  \]
\end{prop}

\begin{proof}
  The archimedean Gaussian has weight $2$, and the
  local calculation in the proof of
  Proposition~\ref{prop:theta-xi} gives right
  invariance under $U_1(M)$. The central-character
  calculation gives $\omega$.

  Proposition~\ref{prop:hecke-auto} factors as the
  product of its archimedean and finite components.
  Because $\varphi_1$ and $\varphi_2$ are invariant
  under $B_\infty^1$, its archimedean factor is
  \[
    \int_{B_\infty^1}
      \Phi_{\calo,\infty}
      \paren{b_{0,\infty}b}\,db
    =
    e^{-2\pi x_\infty}.
  \]
  Thus
  $\verts{x_\infty}e^{-2\pi x_\infty}
  = W_2\paren{a(x_\infty)}$.
  After removing this archimedean Whittaker factor,
  Proposition~\ref{prop:hecke-auto} gives the displayed
  formula for
  $\kappa^\infty(x^\infty,
  \varphi_1 \otimes \varphi_2, \Phi_\calo)$
  with the finite Hecke operator of
  Equation~\ref{eq:hecke-finite}. Substituting
  $\varphi_i = \sbrac{\xi_i}$ and applying
  Equation~\ref{eq:xi-varphi} gives the formula for
  $\xi_1 \otimes \xi_2$.
\end{proof}

%%%%%%%%%%%%%%%%%%%%%%%%%%%%%

\section{Construction of optimal forms}
\label{sec:opt}

\subsection{The automorphic avatar of optimal forms}

Let $K$ be an imaginary quadratic field and
$\chi: K^\times \bs K_\A^\times \rightarrow \C^\times$
a finite character as in Section \ref{sec:theta}.
The tensor product representation
$\pi(\chi)\otimes_{\Q(\chi)}\pi(\chi^{-1})$ has the
distinguished newvector
$\varphi^\new=\varphi_\chi^\new
\otimes\varphi_{\chi^{-1}}^\new$; see
Equation~\ref{eq:theta}.

We construct another element, the \emph{optimal form}
$\varphi^\opt$, which depends only on the antinorm
$\xi:=\chi^{1-c_K}$. This is a ring class character.
Let $c:=c(\xi)$ be its conductor and put
$\calo_c:=\Z+c\calo_K$. Then $\xi$ is trivial on
$\wh\calo_c^\times$, and hence may be viewed as a
character of
\[
  K^\times\bs K_\A^\times/
  K_\infty^+\wh\calo_c^\times
  =
  K^\times\bs\wh K^\times/\wh\calo_c^\times
  =:
  \Pic^+(\calo_c).
\]

\begin{definition}
\label{def:opt}
  Write $\calo_c = \Z + \Z\tau$ and set
  $\delta := \tau - \overline{\tau}$. Then
  $\cald := \delta \calo_c$ is the different ideal, i.e.
  the $\calo_c$-ideal generated by
  $x-\overline{x}$ for $x\in\calo_c$.
  For each $\alpha \in \calo_c/\cald$, define
  \[
    \Phi_\alpha^\opt
      = \Phi_{\alpha,\infty}^\opt\otimes
        \Phi_\alpha^{\opt,\infty}
        \in\cals(K_\A)
  \]
  as follows:
  \begin{enumerate}
    \item
      $\Phi_{\alpha,\infty}^\opt(x)
        = e^{-2\pi\verts{x}^2}$;
    \item
      $\Phi_\alpha^{\opt,\infty}$ is the
      characteristic function of
      $\wh\calo_c+\alpha/\delta$.
  \end{enumerate}
  Let $(-1)^\infty\in\A^\times$ be the idele with
  archimedean component $1$ and finite components $-1$.
  Using the theta series
  $\theta(g,\chi^{c_K},\Phi_\alpha^\opt)$, define
  \begin{equation}
  \label{eq:opt}
    \varphi^\opt(g_1,g_2)
    :=
    \sum_{\alpha\in\calo_c/\cald}
      \theta\paren{
        g_1,
        \chi^{c_K},
        \Phi_\alpha^\opt}
      \theta\paren{
        g_2a\paren{(-1)^\infty},
        \paren{\chi^{-1}}^{c_K},
        \Phi_{-\alpha}^\opt}.
  \end{equation}
\end{definition}

\begin{remark}
\label{rem:optimal-name}
  The adjective ``optimal'' was chosen because of the
  relation with optimal embeddings. For the history of
  this terminology, see
  \cite[Remark~30.3.17]{voight}; for the theory of
  optimal embeddings and orders, see
  \cite[Section~3]{eichler-1955},
  \cite[Sections~1 and~3]{gross-1987}, and
  \cite[Section~30.3]{voight}.
  An embedding $\iota:K\longhookrightarrow B$ is
  \emph{optimal} with respect to orders
  $\calo_c\subset K$ and $\calo\subset B$ if
  $\iota(K)\cap\calo=\iota(\calo_c)$. Thus the
  quadratic order cut out inside $K$ by the quaternionic
  order is exactly $\calo_c$, rather than a larger order.
  This also explains the Schwartz functions in
  Definition~\ref{def:opt}. In the quaternionic
  realization below, an optimal embedding and a
  decomposition $B=K\oplus Kj$ make the local standard
  lattice decomposable away from the different $\cald$.
  At a prime $q\mid\cald$, its characteristic function
  instead decomposes as a finite sum indexed by
  $\alpha\in\calo_c/\cald$; the two factors in the
  $\alpha$-summand are the characteristic functions of
  $\calo_{c,q}+\alpha/\delta$ and
  $\calo_{c,q}-\alpha/\delta$.
  Consequently, Equation~\ref{eq:opt} packages the
  decomposition of the standard quaternionic theta kernel
  attached to an optimal embedding. This realization is
  made precise in Proposition~\ref{prop:theta-opt}. It
  also explains why $\varphi^\opt$ is defined before an
  auxiliary quaternion algebra is chosen.
\end{remark}
  
We can write down the optimal form's Kirillov functions.
By Equation~\ref{eq:kappap-chi-phi},
\begin{equation}
\label{eq:kappa-def}
  \kappa\paren{x,\chi^{c_K},\Phi_\alpha^\opt}
  =
  \verts{x}^{\frac{1}{2}}
  \int_{K_\A^1}
    \Phi_\alpha^\opt(hh_0)\chi(hh_0)\,dh,
\end{equation}
where $x=\rmn(h_0)$. The Kirillov function of
$\varphi^\opt$ is
\begin{equation}
\label{eq:kappa-opt}
  \kappa^\opt(x_1, x_2)
  =
  \sum_{\alpha\in\calo_c/\cald}
    \kappa\paren{x_1, \chi^{c_K},
      \Phi_\alpha^\opt}
    \kappa\paren{x_2(-1)^\infty,
      \paren{\chi^{-1}}^{c_K},
      \Phi_{-\alpha}^\opt}.
\end{equation}

\subsubsection*{Nonvanishing of the optimal period}

Before turning to the theta identity, we establish the
nonvanishing needed later to use the optimal form as a test
vector in Theorem~\ref{thm:hv-reform}. The argument is
local: the finite part of the global optimal vector factors
as a restricted tensor product of local optimal vectors, whose
invariant periods are nonzero by unitarity and positivity.

The following lemma is a small integral refinement of
Hilbert's theorem 90. The usual field-valued statement is
recalled by Serre~\cite[Chapter~X,
Proposition~3]{serre-local-fields};
the additional point is
that the coboundary can be represented by a unit of the
possibly nonmaximal order.

\begin{lem}
\label{lem:integral-hilbert-90}
  Let $F$ be a non-archimedean local field, let $E$ be a
  quadratic \'{e}tale $F$-algebra, and let
  $\calo \subset E$ be an order. Choose
  $\tau \in \calo$ such that
  $\calo = \calo_F + \calo_F \tau$,
  and define $\delta := \tau - \overline{\tau}$.
  If $h \in 1 + \delta\calo$ and
  $\rmn_{E/F}(h) = 1$, then
  $h = t/\overline{t}$
  for some $t \in \calo^\times$.
\end{lem}

\begin{proof}
  Suppose first that $E = F \times F$. Write
  $\tau = (r, s)$, set $d := r - s$, and set
  $\mff := d\calo_F$. Then
  \[
    \calo
    =
    \set{(a, b) \in \calo_F^2 \Mid
      a \equiv b \pmod{\mff}}.
    \]
    Moreover, $\delta = (d, -d)$. If
  $h = (u, u^{-1}) \in 1 + \delta\calo$, then
  $u \equiv 1 \pmod{\mff}$. Hence
  $t := (u, 1) \in \calo^\times$ satisfies
  $h = t/\overline{t}$.

  Suppose now that $E/F$ is a field. Write
  $h = 1 + \delta y$ with $y \in \calo$. Since
  $\overline{\delta} = -\delta$, the norm-one condition
  gives
  $y - \overline{y}
  = \delta y\overline{y}$.
  If $y = 0$, take $t = 1$.
  Otherwise,
  $h = 1 + \delta y = y/\overline{y}$.
  Write $y = a + b\tau$ with
  $a, b \in \calo_F$. Comparing coefficients gives
  $b = y\overline{y} = \rmn_{E/F}(y)$.

  Normalize $v_E$ to have value group $\Z$, and set
  $n := v_E(y) \geq 0$. If $n = 0$, then
  $b \in \calo_F^\times$, and hence
  $y^{-1} = \frac{\overline{y}}{b}
  \in \calo$.
  Thus we may take $t = y$.

  Suppose that $n > 0$. Then
  $v_E(b) = 2n$, while
  $v_E(b\tau) > n$, so $v_E(a) = n$. Since
  $a \in F$, the ramification index $e(E/F)$ divides
  $n$. Choose $s \in F^\times$ with
  $v_E(s) = -n$. Then
  $sa \in \calo_F^\times$
  and $sb \in \calo_F$,
  so $t := sy$ belongs to $\calo$. Moreover,
  $\rmn_{E/F}(t) =
  s^2b \in \calo_F^\times$,
  and therefore $t\in\calo^\times$. Finally,
  $s = \overline{s}$, so
  $t/\overline{t} = y/\overline{y} = h$.
\end{proof}

We apply Lemma~\ref{lem:integral-hilbert-90} to show
that every local optimal vector has nonzero invariant
period. For a finite prime $q$, write
$K_q:=K\otimes_\Q\Q_q$ and
$\calo_{c,q}:=\calo_c\otimes_\Z\Z_q$, let $\delta_q$
be the $q$-component of $\delta$, and put
$\xi_q:=\chi_q/\chi_q^{c_K}$.
For $a \in \calo_{c,q}/\delta_q\calo_{c,q}$, let
$\Phi_{a,q}^\opt := \one_{\calo_{c,q} + a/\delta_q}$ and define
\[
  W_{a,q}^\opt(g)
  := W\paren{g,\chi_q^{c_K},\Phi_{a,q}^\opt}.
\]
Define
\[
  W_q^\opt(g_1,g_2)
  :=
  \sum_{a\in\calo_{c,q}/\delta_q\calo_{c,q}}
    W_{a,q}^\opt(g_1)
    W\paren{
      g_2a(-1),
      \paren{\chi_q^{-1}}^{c_K},
      \Phi_{-a,q}^\opt}.
\]
By local multiplicity one, the space
\[
  \Hom_{\C\sbrac{\GL_2(\Q_q)}}
  \paren{\pi(\chi_q)\otimes\pi(\chi_q^{-1}),\C}
\]
is one-dimensional. Fix a nonzero element and denote the
corresponding pairing on the Whittaker models by
$\calp_{\rs,q}$.

\begin{prop}
\label{prop:local-opt-pairing-nonzero}
  For every finite prime $q$,
  \[
    \calp_{\rs,q}\paren{W_q^\opt} \neq 0.
  \]
\end{prop}

\begin{proof}
  Since $\chi_q$ has finite order, $\pi(\chi_q)$ is
  unitary. Define an anti-linear equivariant isomorphism
  \[
    J_q:\calw\paren{\pi(\chi_q),\psi_q}
    \longrightarrow
    \calw\paren{\pi(\chi_q^{-1}),\psi_q}
  \]
  by
  \[
    (J_qW)(g):=\overline{W\paren{a(-1)g}}.
  \]
  A comparison of Kirillov functions, using
  $\Phi_{-a,q}^\opt(x)=\Phi_{a,q}^\opt(-x)$, gives
  \[
    W_q^\opt
    =
    \chi_q(-1)
    \sum_{a\in\calo_{c,q}/\delta_q\calo_{c,q}}
      W_{a,q}^\opt\otimes J_qW_{a,q}^\opt.
  \]
  Choose a positive-definite invariant Hermitian pairing
  $\pair{-,-}_q$ on $\pi(\chi_q)$. By Schur's lemma,
  there is a scalar $c_q\in\C^\times$ such that
  \[
    \calp_{\rs,q}\paren{W_1\otimes J_qW_2}
    =c_q\pair{W_1,W_2}_q.
  \]
  Consequently,
  \[
    \calp_{\rs,q}\paren{W_q^\opt}
    =
    \chi_q(-1)c_q
    \sum_a
      \pair{W_{a,q}^\opt,W_{a,q}^\opt}_q.
  \]

  It remains to show that at least one summand is nonzero.
  Take $h_0=\delta_q^{-1}$ and
  $x=\rmn_{K_q/\Q_q}(h_0)$. Equation~
  \ref{eq:kappap-chi-phi} gives
  \[
    \kappa\paren{x,\chi_q^{c_K},\Phi_{1,q}^\opt}
    =
    \verts{x}^{\frac{1}{2}}\chi_q(\delta_q)^{-1}
    \int_{\paren{1+\delta_q\calo_{c,q}}^1}
      \chi_q(h)\,dh,
  \]
  where the superscript $1$ denotes the kernel of the
  norm to $\Q_q^\times$. By
  Lemma~\ref{lem:integral-hilbert-90}, write
  $h = t/\overline{t}$ with
  $t \in \calo_{c,q}^\times$. Then
  \[
    \chi_q(h)
    = \frac{\chi_q(t)}
      {\chi_q(\overline{t})}
    = \xi_q(t)
    = 1,
  \]
  since $\xi_q$ is trivial on $\calo_{c,q}^\times$.
  Thus the last integral is the nonzero volume of
  $\paren{1+\delta_q\calo_{c,q}}^1$, and
  $W_{1,q}^\opt \neq 0$. The sum of Hermitian norms above
  is therefore strictly positive, which proves the claim.
\end{proof}

\begin{prop}
\label{prop:opt-pairing-nonzero}
  Let
  \[
    \calp:
    \pi(\chi)^\infty\otimes
    \pi(\chi^{-1})^\infty
    \longrightarrow\C
  \]
  be a nonzero $\GL_2(\A^\infty)$-invariant pairing.
  Then
  \[
    \calp\paren{\varphi^\opt} \neq 0,
  \]
  where the fixed archimedean weight-one vectors are
  suppressed from the notation.
\end{prop}

\begin{proof}
  The Chinese remainder theorem identifies the global
  indexing set with the product of its local analogues:
  \[
    \calo_c/\cald
    \iso
    \prod_{q \mid d(\xi)}
    \calo_{c,q}/\delta_q\calo_{c,q}.
  \]
  Factoring each theta kernel into its local Whittaker
  functions therefore gives
  \[
    W_{\varphi^\opt}^\infty
    =
    \bigotimes_q'W_q^\opt.
  \]
  Normalize the local pairings so that they take value $1$
  on the unramified local optimal vectors for almost all
  $q$. Their restricted tensor product is a nonzero
  $\GL_2(\A^\infty)$-invariant pairing, and
  Proposition~\ref{prop:local-opt-pairing-nonzero} gives
  \[
    \bigotimes_q'\calp_{\rs,q}
    \paren{W_{\varphi^\opt}^\infty}
    =
    \prod_q\calp_{\rs,q}\paren{W_q^\opt}
    \neq 0.
  \]
  The space of invariant pairings on
  $\pi(\chi)^\infty\otimes
  \pi(\chi^{-1})^\infty$ is one-dimensional. Hence
  $\calp$ is a nonzero scalar multiple of this pairing.
\end{proof}

\subsection{A theta identity for the automorphic avatar of optimal forms}
 
We again assume that $K$ is imaginary. Let $B$ be a
definite quaternion algebra over $\Q$ equipped with an
embedding $K \longhookrightarrow B$. Choose
$j \in B^\times$ such that $B = K \oplus Kj$ and
$jx = \overline{x}j$ for every $x \in K$. Set
$V_1 := \paren{K, \rmn_{K/\Q}}$ and
$V_2 := \paren{Kj, -j^2\rmn_{K/\Q}}$, so that
$B = V_1 \oplus V_2$ as an orthogonal direct sum. As in
Section~\ref{sec:theta-two}, take $\xi_1 := \one$ and
$\xi_2 := \xi := \chi^{1 - c_K}$. The character $\xi$
is a ring class character. Let $c := c(\xi)$ be its
conductor, so that $\calo_c = \Z + c\calo_K$ and
$d(\xi) := \disc(\calo_c) = c^2\disc(K)$. Assume that
$\gcd\paren{d(\xi), \disc(B)} = 1$. In particular,
$B_q$ is split for every prime $q \mid d(\xi)$.

\begin{definition}
\label{defn-xi-opt}
  An Eichler order $\calo \subset B$ is
  $\xi$-optimal if
  \begin{enumerate}
    \item
      $K \cap \calo = \calo_c$;
    \item
      for every prime $q \nmid d(\xi)$, there is
      an element $j_q \in B_q^\times$ such that
      $\calo_q
        =
        \calo_{K,q}
        \oplus
        \calo_{K,q}j_q$ and
      $j_qx = \overline{x}j_q$
      for all $x \in K_q$, and the principal ideal
      $j_q^2\Z_q$ is the local reduced
      discriminant ideal of $\calo_q$.
  \end{enumerate}
\end{definition}

The $\xi$-optimal orders needed
in the application to the auxiliary
prime always exist.

\begin{lemma}
\label{lem:xi-optimal-order}
  Let $p \nmid d(\xi)$ be inert in $K$, and let $B$ be
  the definite quaternion algebra of reduced
  discriminant $p$, equipped with an embedding
  $K \longhookrightarrow B$. Then there is a
  $\xi$-optimal maximal order $\calo \subset B$.
  In particular,
  $K \cap \calo = \calo_c$.
\end{lemma}

\begin{proof}
  For every finite prime $q \neq p$, the algebra $B_q$
  is split. Choose an isomorphism
  \[
    B_q \iso \End_{\Q_q}(K_q)
  \]
  that identifies the given copy of $K_q$ with its
  action on itself by multiplication, and set
  \[
    \calo_q
    :=
    \End_{\Z_q}\paren{\calo_{c,q}}.
  \]
  This is a maximal order, and
  \begin{align*}
    K_q \cap \calo_q
    =
    \set{a \in K_q \Mid
      a\calo_{c,q}
      \subseteq \calo_{c,q}}
    =
    \calo_{c,q}.
  \end{align*}
  At $p$, the extension $K_p/\Q_p$ is unramified and
  $\calo_{c,p} = \calo_{K,p}$. The unique maximal
  order $\calo_p$ of the division algebra $B_p$
  therefore satisfies
  \[
    K_p \cap \calo_p = \calo_{K,p}.
  \]

  These local maximal orders agree with a fixed maximal
  order at all but finitely many primes. The
  local--global dictionary for lattices therefore gives
  a maximal order $\calo \subset B$ with localizations
  $\calo_q$; see
  \cite[Theorem~9.4.9]{voight}.
  For $q \neq p$ with $q \nmid d(\xi)$, one has
  $\calo_{c,q} = \calo_{K,q}$, and the regular
  representation gives
  \[
    \End_{\Z_q}\paren{\calo_{K,q}}
    =
    \calo_{K,q}
    \oplus
    \calo_{K,q}j_q,
  \]
  where $j_q$ acts by the nontrivial automorphism of
  $K_q/\Q_q$ and satisfies $j_q^2 = 1$.

  At $p$, choose a quaternionic uniformizer
  $\Pi_p \in B_p$ satisfying
  $\Pi_p x = \overline{x}\Pi_p$ for all $x \in K_p$
  and $\Pi_p^2 = p$.
  Then
  \[
    \calo_p
    =
    \calo_{K,p}
    \oplus
    \calo_{K,p}\Pi_p.
  \]
  These descriptions verify the second condition in
  Definition~\ref{defn-xi-opt}.
\end{proof}
 
For a positive integer $M$, let
$M^\infty \in \A^{\infty,\times}$
denote the image of $M \in \Q^\times$ in the finite
ideles.

The following proposition identifies the optimal form
with the theta lift of the standard Schwartz function
attached to a $\xi$-optimal order. The scalar
$c_{\chi,\calo}$ records the normalization of the
adelic comparison between $V_1$ and $V_2$.

\begin{prop}
\label{prop:theta-opt}
  Let $\calo$ be a $\xi$-optimal Eichler order of $B$
  whose reduced discriminant $M$ is coprime to
  $d(\xi)$. There is a unique scalar
  $c_{\chi,\calo} \in \Q(\chi)^\times$
  such that
  \[
    c_{\chi,\calo}
    \theta\paren{
      g,
      \one \otimes \xi,
      \Phi_\calo}
    =
    M^{-\frac{1}{2}}
    \varphi^\opt\paren{
      g,
      g a\paren{(M^\infty)^{-1}}}.
  \]
\end{prop}

\begin{proof}
  We first decompose the standard Schwartz function
  $\Phi_\calo$ with respect to
  \[
    B_\A = V_{1,\A} \oplus V_{2,\A}.
  \]

  At the archimedean place, choose
  $j_\infty \in K_\infty j$ with
  $j_\infty^2 = -1$. Using the identifications
  \[
    K_\infty \iso \C,
    \qquad
    K_\infty j_\infty \iso \C,
    \quad
    zj_\infty \longleftrightarrow z,
  \]
  write $z_i = x_i + y_i\sqrt{-1}$
  for $x_i, y_i \in \R$.
  Set
  \[
    \Phi_{i,\infty}(z_i)
    :=
    e^{-2\pi(x_i^2 + y_i^2)}
    \qquad (i = 1,2).
  \]
  Then
  \[
    \Phi_{\calo,\infty}
    =
    \Phi_{1,\infty}
    \otimes
    \Phi_{2,\infty}.
  \]

  Let $q$ be a finite prime with
  $q \nmid d(\xi)$. By
  Definition~\ref{defn-xi-opt}, choose $j_q$ such that
  \[
    \calo_q
    =
    \calo_{K,q}
    \oplus
    \calo_{K,q}j_q.
  \]
  Set $\Phi_{1,q} := \one_{\calo_{K,q}}$ and
  $\Phi_{2,q} := \one_{\calo_{K,q}j_q}$.
  Then
  \[
    \Phi_{\calo,q}
    =
    \Phi_{1,q}
    \otimes
    \Phi_{2,q}.
  \]

	Now let $q \mid d(\xi)$. Since
  $\gcd\paren{M, d(\xi)} = 1$,
  the order $\calo_q$ is maximal, and $B_q$ is split.
  The embedding
  $\calo_{c,q} \longhookrightarrow \calo_q$
  is optimal. By the uniqueness of local optimal
  embeddings into a split maximal order
  \cite[Proposition~30.5.3(a)]{voight}, we may identify
  \[
    (B_q, \calo_q, K_q, \calo_{c,q})
    \iso
    \paren{
      \End_{\Q_q}(K_q),
      \End_{\Z_q}\paren{\calo_{c,q}},
      K_q,
      \calo_{c,q}},
  \]
  where $K_q$ acts on itself by multiplication. Let
  $j_q$ be the endomorphism induced by the nontrivial
  automorphism of $K_q/\Q_q$. Then $j_q^2 = 1$ and
  $j_qx = \overline{x}j_q$ for every $x \in K_q$.
  Recall that $\calo_{c,q} = \Z_q + \Z_q\tau$ and
  $\delta = \tau - \overline{\tau}$. For
  $a,b \in K_q$, one has
  \begin{align*}
    a + bj_q
    \in \End_{\Z_q}\paren{\calo_{c,q}}
    \quad\Longleftrightarrow\quad
      &a + b \in \calo_{c,q}, \\
      &a\tau + b\overline{\tau}
        \in \calo_{c,q} \\
    \quad\Longleftrightarrow\quad
      &a + b \in \calo_{c,q}, \\
      &a\delta \in \calo_{c,q}.
  \end{align*}
  It follows that
  \[
    \calo_q
    =
    \coprod_{\alpha_q \in
      \calo_{c,q}/\delta\calo_{c,q}}
    \left(
      \calo_{c,q}
      + \frac{\alpha_q}{\delta}
    \right)
    +
    \left(
      \calo_{c,q}
      - \frac{\alpha_q}{\delta}
    \right)j_q.
  \]
  Consequently, under
  $B_q = K_q \oplus K_qj_q$,
  we have
  \[
    \Phi_{\calo,q}
    =
    \sum_{\alpha_q \in
      \calo_{c,q}/\delta\calo_{c,q}}
    \one_{\calo_{c,q}+\alpha_q/\delta}
    \otimes
    \one_{(\calo_{c,q}-\alpha_q/\delta)j_q}.
  \]

  The Chinese remainder theorem identifies
  \[
    \calo_c/\cald
    \iso
    \prod_{q \mid d(\xi)}
    \calo_{c,q}/\delta\calo_{c,q}.
  \]
  For $\alpha \in \calo_c/\cald$, let
  $\Phi_{1,\alpha}$ and $\Phi_{2,\alpha}$ be the
  global Schwartz functions with the local components
  described above, and with
  \begin{align*}
    \Phi_{1,\alpha,q}
    &:=
    \one_{\calo_{c,q}+\alpha/\delta}, \\
    \Phi_{2,\alpha,q}
    &:=
    \one_{(\calo_{c,q}-\alpha/\delta)j_q}
  \end{align*}
  when $q \mid d(\xi)$. Then
  \begin{equation}
  \label{eq:standard-opt-decomposition}
    \Phi_\calo
    =
    \sum_{\alpha \in \calo_c/\cald}
    \Phi_{1,\alpha}
    \otimes
    \Phi_{2,\alpha}.
  \end{equation}
  Moreover, $\Phi_{1,\alpha} = \Phi_\alpha^\opt$.

  Under the character correspondence preceding
  Equation~\ref{eq:xi-chi}, the pair
  $(\chi^{c_K}, (\chi^{-1})^{c_K})$
  corresponds to $\paren{\one,\xi}$. Therefore,
  Equations~\ref{eq:standard-opt-decomposition}
  and~\ref{eq:xi-chi} give
  \begin{equation}
  \label{eq:theta-opt-decomposition}
    \theta\paren{
      g, \one \otimes \xi, \Phi_\calo}
    = \sum_{\alpha \in \calo_c/\cald}
      \theta\paren{
        g, \chi^{c_K}, \Phi_{1,\alpha}}
      \theta\paren{
        g, \paren{\chi^{-1}}^{c_K},
        \Phi_{2,\alpha}}.
  \end{equation}

	For the second factor, define a $K_\A$-linear
  isomorphism
  \[
    \iota_\calo: V_{1,\A} \iso V_{2,\A}
  \]
  place by place. At the archimedean place, set
  $\iota_{\calo,\infty}(x) := xj_\infty$.
  At a finite prime $q \nmid d(\xi)$, the element
  $j_q$ in Definition~\ref{defn-xi-opt} is determined
  up to multiplication by $K_q^\times$. Since $K_q/\Q_q$
  is split or unramified, the norm map on units is
  surjective. We may therefore choose $j_q$ so that
  $j_q^2 = q^{\ord_q(M)}$,
  and set
  $\iota_{\calo,q}(x) := xj_q$.
  At a prime $q \mid d(\xi)$, use the conjugation
  operator $j_q$ constructed above and again set
  $\iota_{\calo,q}(x) := xj_q$.

  Recall that $(-1)^\infty \in \A^\times$ has
  archimedean component $1$ and finite components $-1$.
  The preceding choices give
  \[
    Q\paren{\iota_\calo}
    =
    (-1)^\infty M^\infty,
    \qquad
    \verts{Q\paren{\iota_\calo}}
    =
    M^{-1},
  \]
  and
  \[
    \iota_\calo^*\Phi_{2,\alpha}
    =
    \Phi_{-\alpha}^\opt.
  \]
  Write $\iota_\calo(1) = a_\calo j$ and set
  $c_{\chi,\calo} := \chi\paren{a_\calo}$. We call
  $c_{\chi,\calo}$ the comparison scalar attached
  to these local comparison maps.
  Proposition~\ref{prop:theta-comparison}, applied with
  the character $\chi^{-1}$, gives
  \[
    \theta_{V_2}\paren{
      g,
      \paren{\chi^{-1}}^{c_K},
      \Phi_{2,\alpha}}
    =
    c_{\chi,\calo}^{-1}
    M^{-\frac12}
    \theta_{V_1}\paren{
      g a\paren{(M^\infty)^{-1}(-1)^\infty},
      \paren{\chi^{-1}}^{c_K},
      \Phi_{-\alpha}^\opt}.
  \]
  Combining this with
  Equation~\ref{eq:theta-opt-decomposition} and
  Definition~\ref{def:opt}, we obtain
  \[
    c_{\chi,\calo}
    \theta\paren{g, \one \otimes \xi, \Phi_\calo}
    =
    M^{-\frac{1}{2}}
    \varphi^\opt\paren{
      g,
      g a\paren{(M^\infty)^{-1}}}.
  \]

  Proposition~\ref{prop:opt-pairing-nonzero}
  implies that $\varphi^\opt \neq 0$. The displayed
  identity therefore also proves that
  $c_{\chi,\calo}$ is unique.
\end{proof}

By Proposition~\ref{prop:hecke-complex}, the theta
series in Proposition~\ref{prop:theta-opt} is
holomorphic of weight $2$. Its finite Kirillov function
is
\[
  \kappa^\infty\paren{
    x^\infty,
    \one \otimes \xi,
    \Phi_\calo}
  =
  \verts{x^\infty}
  \pair{\sbrac{\one},
    \rmt_{\Phi_\calo^\infty}
    \paren{x^\infty}
    \sbrac{\xi}}
\]
for every $x^\infty \in \rmn\paren{B_{\A^\infty}^\times}$.

\subsection{The modular avatar of optimal forms}
Let $K/\Q$ be an imaginary quadratic extension, let
$\chi$ be a character of $\gal(\overline{K}/K)$, and let
$\xi := \chi^{1 - c_K}$ be its antinorm of conductor
$c = c(\xi)$. Retain the notation
$\calo_c = \Z + \Z\tau$ from Definition~\ref{def:opt},
put $\delta := \tau - \overline{\tau}$, and let
$\cald := \delta\calo_c$ be the different of $\calo_c$.

In general, a modular form
can be associated to $\chi$ and any locally constant function
with a precise description of its $q$-expansion. 
Let $\wh K^1$ denote the subgroup of $\wh K^\times$ of norm $1$
and choose a measure on $\wh K^1$ such that the volume of its maximal compact
subgroup $\wh K^1 \cap \wh\calo_K^\times$ is $1$.

\begin{lemma}
  \label{lem:achi}
  Let $\Phi_\infty$ be the standard archimedean Schwartz function
  from Definition~\ref{def:standard-Eichler}.
  \begin{enumerate}[(a)]
    \item For any locally constant function
      $\Phi^\infty: \wh K \lra \C$ with compact support,
      there is a weight-$1$ modular form
      $f_{\chi,\Phi^\infty}$ whose automorphic avatar is
      \[
        \theta\paren{
          g,
          \chi^{c_K},
          \Phi^\infty \otimes \Phi_\infty}.
      \]
    \item For $u \in \wh\Z^\times$, the form
      $f_{\chi,\Phi^\infty}$ has the
      $q$-expansion
      \[
        f_{\chi,\Phi^\infty}(q, u)
        =
        \sum_{r \in \Q_+}
          a_{\chi,\Phi^\infty}(r, u)q^r.
      \]
      The coefficient
      $a_{\chi,\Phi^\infty}(r, u)$ vanishes
      unless $ru = \rmn(t_0)$ for some
      $t_0 \in \wh K^\times$. In that case,
      \[
        a_{\chi,\Phi^\infty}(r, u)
        =
        \int_{\wh K^1}
          \Phi^\infty(tt_0)
          \chi(tt_0)\,dt.
      \]
    \end{enumerate}
\end{lemma}

\begin{proof}
  Let
  $\varphi :=
  \theta\paren{g, \chi^{c_K},
  \Phi^\infty \otimes \Phi_\infty}$.
  The existence of its modular avatar $f_{\chi, \Phi^\infty}(z, u)$
  follows directly from Proposition \ref{prop:f-phi}.

  By Proposition~\ref{prop:a-kappa},
  the coefficient of $q^r$ is
  \[
    a_{\chi,\Phi^\infty}(r, u)
    =
    r^{\frac{1}{2}}\kappa_\varphi(ru).
  \]
  Equation~\ref{eq:kappa-chi-phi} shows that this
  vanishes unless $ru = \rmn(t_0)$ for some
  $t_0 \in \wh K^\times$. In that case,
  $\verts{ru} = r^{-1}$, since
  $u \in \wh\Z^\times$, and therefore
  \[
    r^{\frac{1}{2}}\kappa_\varphi(ru)
    =
    \int_{\wh K^1}
      \Phi^\infty(tt_0)
      \chi(tt_0)\,dt.
  \]
\end{proof}

\begin{example}
  Let $c(\chi) \subset \calo_K$ be the conductor ideal of $\chi$.
  Classically, there is a newform $f_\chi$ associated to 
  the induced Galois representation $\rho = \Ind_K^\Q(\chi)$,
  whose $q$-expansion $f_\chi = \sum_n a_\chi(n) q^n$
  is determined by the equality of $L$-functions:
  \[
    \sum a_\chi(n) n^{-s}
      = L(f, s)
      = L(\rho, s)
      = L(\chi, s)
      = \prod_{\wp \nmid c(\chi)}
        \paren{1 - \chi(\wp) \rmn(\wp)^{-s}}^{-1}.
  \]
  Let $\varphi_\chi$ be the automorphic avatar of $f_\chi$.
  Then $\varphi_\chi$ can also be defined as a theta lifting as in
  Equation~\ref{eq:theta}.
  Let $\Phi^\infty = \bigotimes_{v\nmid \infty} \Phi_{\chi_v}$
  as defined in Section~\ref{sec:theta}.
  Then $f_{\chi, \Phi^\infty} = f_\chi$.
\end{example}

Let $\pi^\infty(\chi)$ be the irreducible representation of
$\GL_2(\wh \Q)$ of modular forms generated by $f_\chi$.  
Consider the tensor product of representations
generated by $f_\chi$ and $f_{\chi^{-1}}$:
\[
	\pi^\infty(\chi) \otimes \pi^\infty(\chi^{-1}).
\]
Notice that since $\chi$ is unitary,
$f_{\chi^{-1}}$ can be obtained from $f_\chi$ by complex conjugation
on the coefficients in its $q$-expansion.
Recall from Definition \ref{def:opt} that
there is an automorphic optimal form
$\varphi^\opt \in \pi(\chi) \otimes \pi(\chi^{-1})$.
We define the modular optimal form $f^\opt$
through a $q$-expansion like in Equation~\ref{eq:q-xn}.
\begin{proposition}
\label{prop:q-opt}
  Set $D_c := \verts{\disc(\calo_c)} = -\delta^2$.
  The modular optimal form $f^\opt$ has the
  componentwise $q$-expansion
  \begin{align*}
    f^\opt\paren{q_1, q_2, u_1, u_2}
    &=
    \sum_{r_1, r_2 \in \Q_+}
      a^\opt\paren{r_1, r_2, u_1, u_2}
      q_1^{r_1} q_2^{r_2},
    \\
    \intertext{where
      $u_1, u_2 \in \wh\Z^\times$ and}
    a^\opt\paren{r_1, r_2, u_1, u_2}
    &:=
    \sum_{\alpha \in \calo_c/\cald}
      a_{\chi, \Phi_\alpha^{\opt,\infty}}
      \paren{r_1, u_1}
      a_{\chi^{-1},
        \Phi_{-\alpha}^{\opt,\infty}}
      \paren{r_2, -u_2}.
  \end{align*}
  In particular, all exponents in both variables lie in
  $D_c^{-1}\Z$.
\end{proposition}

\begin{proof}
  For each $\alpha \in \calo_c/\cald$, apply
  Proposition~\ref{prop:a-kappa} separately to the two
  one-variable theta series in Equation~\ref{eq:opt}.
  For any weight-one automorphic form $\varphi$ with
  modular avatar $f_\varphi$, the identity
  $a(ru)a\paren{(-1)^\infty}=a(-ru)$ and
  Proposition~\ref{prop:a-kappa} give
  \[
    a_{f_{\varphi(\,\cdot\,a((-1)^\infty))}}(r, u)
    =
    r^{\frac{1}{2}} \kappa_\varphi(-ru)
    =
    a_{f_\varphi}(r, -u).
  \]
  Thus right translation by $a((-1)^\infty)$ changes
  the component index, rather than the sign of the
  positive Fourier exponent. Summing over $\alpha$ and
  using Equation~\ref{eq:kappa-opt} gives the displayed
  two-variable $q$-expansion.

  For a fixed $\alpha$, Lemma~\ref{lem:achi} shows
  that $a_{\chi, \Phi_\alpha^{\opt,\infty}}(r, u)$
  can be nonzero only if
  $ru = \rmn(t_0)$
  for some
  $t_0 \in \wh\calo_c + \frac{\alpha}{\delta}$.
  As $\alpha$ varies, the Chinese remainder theorem
  gives
  \[
    \bigcup_{\alpha \in \calo_c/\cald}
    \paren{
      \wh\calo_c + \frac{\alpha}{\delta}}
    =
    \frac{1}{\delta} \wh\calo_c.
  \]
  Since
  $\rmn_{K/\Q}(\delta)  = -\delta^2 = D_c$,
  it follows that
  $\rmn(t_0) \in D_c^{-1} \wh\Z$.
  Because $u \in \wh\Z^\times$ and $r \in \Q$, this
  implies
  $r \in D_c^{-1}\Z$.
  The same argument applies to the second theta factor.
  Thus the exponents in both variables lie in
  $D_c^{-1}\Z$, as claimed.
\end{proof}

%%%%%%%%%%%%%%%%%%%%%%%%%%%%%%%%%%%%%%%%%%%%%%%%%%%%%%%%

\section{Proof of Theorem~\ref{thm:opt-unique}}
\label{sec:opt-unique}

Enlarge the splitting field $E$, if necessary, so that
$H_c \subset E$ and $\chi$ factors through
$\gal(E/K)$. This does not change $\rho_f$. The
corresponding identification of adjoint-unit groups will
be recorded at the beginning of
Section~\ref{sec:proof-opt}.
Fix a lift
\[
  \widetilde c_K
  \in
  \gal(E/\Q) \setminus \gal(E/K)
\]
of $c_K$.
Let $p \geq 5$ be coprime to $N$ and inert in $K$, and
choose a prime $\mfp$ of $E$ above $p$. There is a unique
$h_\mfp \in \gal(E/K)$ such that
\begin{equation}
\label{eq:frob-cK-h}
  \Frob_\mfp = \widetilde c_K h_\mfp.
\end{equation}
Let $\mfp_c := \mfp \cap H_c$.

Let $B_p$ be the definite quaternion algebra of reduced
discriminant $p$. By
Lemma~\ref{lem:xi-optimal-order}, there is a
$\xi$-optimal maximal order
$\calo_\mfp \subset B_p$. We equip it with the
orientation induced by reduction modulo $\mfp_c$.
Applying the projection construction of
Section~\ref{subsec:hecke-operators} to this oriented
embedding gives the quaternionic vectors
$\sbrac{\one}_\mfp$ and $\sbrac{\xi}_\mfp$.

\subsection{A theta identity for the modular avatar of optimal forms}

Let $\cala_\mfp^\pm$ be the quaternionic automorphic
spaces attached to $\calo_\mfp$ as in
Section~\ref{subsec:hecke-operators}. Since $\one$ and
$\xi$ are trivial on $\A^\times$, these are the spaces
with trivial central character.

For an element $\varphi$ of either space, set
\[
  \operatorname{mass}_\mfp(\varphi)
  :=
  \int_{\sbrac{B_p^\times/\Q^\times}}
    \varphi(x)\,dx.
\]
For $n \geq 1$, let
$n^\infty \in \A^{\infty,\times}$ be the finite idele
induced by $n$, and set
\[
  \rmt_n
  :=
  \rmt_{\Phi_{\calo_\mfp}^\infty}
  \paren{n^\infty}.
\]

Define
\[
  \Theta_\mfp^{\std}:
  \cala_\mfp^+ \otimes \cala_\mfp^-
  \longrightarrow M_2\paren{\Gamma_0(p)}
\]
by requiring
$\Theta_\mfp^{\std}
(\varphi_1 \otimes \varphi_2)$ to be the modular
avatar of
\[
  \theta\paren{
    g,
    \varphi_1 \otimes \varphi_2,
    \Phi_{\calo_\mfp}}.
\]
Propositions~\ref{prop:hecke-complex}
and~\ref{prop:a-kappa} identify its nonconstant
Fourier coefficients. Together with the standard
constant-term calculation, they give
\begin{equation}
\label{eq:Theta-standard}
\begin{split}
  \Theta_\mfp^{\std}
  \paren{\varphi_1 \otimes \varphi_2}
  &=
  \frac{1}{2}
  \operatorname{mass}_\mfp(\varphi_1)
  \operatorname{mass}_\mfp(\varphi_2)
  + \sum_{n \geq 1}
    \pair{\varphi_1,\rmt_n\varphi_2}q^n.
\end{split}
\end{equation}
This is the adelic form of the classical
Gross--Emerton normalization (cf.
\cite[Proposition~5.6]{gross-1987},
\cite{emerton}, and
\cite[Equation~16]{dhrv}).
Set $c_{\chi,\mfp} := c_{\chi,\calo_\mfp}$.
We spell out the passage to modular avatars. Let
$\varphi$ be a weight-one automorphic form with modular
avatar $f_\varphi$, and set
\[
  \varphi^{(p)}(g)
  :=
  \varphi\paren{
    g a\paren{(p^\infty)^{-1}}}.
\]
By left invariance under $a(p) \in \GL_2(\Q)$ and
Equation~\ref{eq:f-varphi}, the modular avatar of
$\varphi^{(p)}$ is
\[
  f_{\varphi^{(p)}}(z) = p^{\frac{1}{2}}f_\varphi(pz).
\]
It follows that the modular avatar of
\[
  p^{-\frac{1}{2}}
    \varphi^\opt\paren{g, g a\paren{(p^\infty)^{-1}}}
\]
is $f^\opt(z,pz)$. Hence
Proposition~\ref{prop:theta-opt}, with $M = p$, gives
\begin{equation}
\label{eq:theta-opt-standard}
  f^\opt(z, pz)
    = c_{\chi,\mfp} \Theta_\mfp^{\std}
    \paren{\sbrac{\one}_\mfp \otimes \sbrac{\xi}_\mfp}.
\end{equation}
Define the $f$-normalized theta lift by
\begin{equation}
\label{eq:Theta-normalized}
  \Theta_\mfp^f := c_{\chi, \mfp} \Theta_\mfp^{\std}.
\end{equation}
Then
\begin{equation}
\label{eq:theta-opt-normalized}
  f^\opt(z, pz)
  =
  \Theta_\mfp^f
  \paren{
    \sbrac{\one}_\mfp
    \otimes
    \sbrac{\xi}_\mfp}.
\end{equation}

\subsection{Uniqueness on the standard component}

Let $g$ be another two-variable modular form of finite
full level satisfying
Equation~\ref{eq:theta-opt-normalized} with $f^\opt$
replaced by $g$.
Choose a positive integer $N'$ divisible by $N$
such that both $f^\opt$ and $g$ descend to
$X(N') \times X(N')$.
Let $X(N')^\circ$ denote the
connected component of $X(N')_\C$ containing the
standard cusp. At the pair of standard cusps of
$X(N')^\circ \times X(N')^\circ$, set
\[
  Q_i := e^{2\pi i z_i/N'}
  \qquad (i = 1, 2),
  \qquad
  Q := e^{2\pi i z/N'}.
\]
We use the theta identity only for sufficiently large
primes $p$ inert in $K$ and coprime to $NN'$. Write
\begin{align*}
  f^\opt(z_1, z_2)
  &=
  \sum_{m, n \geq 0}
    a_{m,n} Q_1^m Q_2^n, \\
  g(z_1, z_2)
  &=
  \sum_{m, n \geq 0}
    b_{m,n} Q_1^m Q_2^n.
\end{align*}
Under the specialization $(z_1, z_2) = (z, pz)$, one has
$Q_1 = Q$ and $Q_2 = Q^p$. Hence comparison of the
coefficient of $Q^k$ gives
\[
  \sum_{m + pn = k} a_{m,n}
  =
  \sum_{m + pn = k} b_{m,n}
\]
for every $k \geq 0$. There are arbitrarily large
primes inert in $K$ and coprime to $NN'$. Taking $p > k$
gives
\[
  a_{k,0} = b_{k,0}.
\]
Now proceed by induction on $n$. Suppose that
$a_{m,n} = b_{m,n}$ for every $m \geq 0$ and
$0 \leq n \leq n_0$. Given $m_0 \geq 0$, choose an
inert prime $p > m_0$ coprime to $NN'$, and set
\[
  k := m_0 + p(n_0 + 1).
\]
Every solution of $m + pn = k$ has
$n \leq n_0 + 1$. Cancelling the terms covered by the
induction hypothesis gives
\[
  a_{m_0,n_0 + 1} = b_{m_0,n_0 + 1}.
\]
Thus the $Q_1,Q_2$-expansions of $g$ and $f^\opt$
agree at the pair of standard cusps. Their difference
therefore vanishes on an analytic neighborhood of that
pair. Since $X(N')^\circ \times X(N')^\circ$ is connected,
analytic continuation gives
\[
  g = f^\opt
  \qquad\text{on}\qquad
  X(N')^\circ \times X(N')^\circ.
\]
This proves the uniqueness assertion in
Theorem~\ref{thm:opt-unique}.

%%%%%%%%%%%%%%%%%%%%%%%%%%%%%%%

\section{Proof of Theorem~\ref{thm:opt}}
\label{sec:proof-opt}

Retain the notation of Section~\ref{sec:opt-unique}.
Let $E_0$ be the original splitting field of
$\rho_f$, and set $H := \gal(E/E_0)$. Since $H$ acts
trivially on $M_{\Ad}$, equivariance forces the image
of every homomorphism defining
$\calu\paren{\Ad(\rho_f)}$ to lie in
\[
  \paren{\calo_E^\times}^H
  =
  \calo_{E_0}^\times.
\]
Thus enlarging $E_0$ to $E$ does not change the
adjoint-unit group.
We have $p \geq 5$ is coprime to $N$ and inert in $K$,
$\mfp$ is a prime of $E$ above $p$, and
\[
  \Frob_\mfp = \widetilde c_K h_\mfp.
\]
Let $\mfp_c := \mfp \cap H_c$. We write $B_p$ for the
definite quaternion algebra of reduced discriminant $p$
and $\calo_\mfp$ for the $\mfp_c$-oriented
$\xi$-optimal maximal order chosen in
Section~\ref{sec:opt-unique}.

The proof combines Equation~\ref{eq:theta-opt-standard}
with two results of
Darmon--Harris--Rotger--Venkatesh~\cite{dhrv}: the
calculation of the adjoint theta lift of the Shimura
class and the relation between a higher Eisenstein
element and an elliptic unit. We keep the dependence on
$\mfp$ explicit and show that it is absorbed by one fixed
element of the adjoint-unit line.

\subsection{Oriented maximal orders}
\label{subsec:oriented-maximal-orders}

An orientation on a maximal order
$\calo' \subset B_p$ is a surjective ring
homomorphism
\[
  \sigma':\calo' \lra \F_{p^2}.
\]
Two oriented maximal orders
$(\calo_1,\sigma_1)$ and $(\calo_2,\sigma_2)$ are
equivalent if there is a $b \in B_p^\times$ such that
\[
  b\calo_1b^{-1} = \calo_2,
  \qquad
  \sigma_1
  =
  \sigma_2 \circ \Ad(b).
\]
Let $\Pic^{\mathrm{or}}(B_p)$ denote the set of
equivalence classes. Deuring's correspondence
identifies this set with the isomorphism classes of
supersingular elliptic curves over
$\overline{\F}_p$; see
\cite[Section~2.1]{dhrv} and
\cite[Corollary~42.3.7]{voight}.

The $\mfp_c$-orientation on $\calo_\mfp$ is the map
\[
  \sigma_\mfp:
  \calo_\mfp \lra \F_{p^2}
\]
whose restriction to $\calo_c$ is reduction modulo
$\mfp_c$.

Consider the adelic class set
\begin{equation}
\label{eq:pic-bp}
  \Pic(B_p)
  :=
  B_p^\times
  \bs
  \wh B_p^\times
  /
  \paren{
    \wh\Q^\times
    \cdot
    \wh\calo_\mfp^\times}.
\end{equation}

\begin{remark}
\label{rem:qhat-pic}
  The factor $\wh\Q^\times$ in the right-hand side of
  Equation \ref{eq:pic-bp} may be omitted. Indeed,
  \[
    \wh\Q^\times
    =
    \Q^\times \cdot \wh\Z^\times,
  \]
  because $\Q$ has class number $1$. The first factor
  is absorbed by $B_p^\times$, and the second is
  contained in $\wh\calo_\mfp^\times$.
\end{remark}

For a finite prime $q$, write
\[
  B_{p,q} := B_p \otimes_\Q \Q_q,
  \qquad
  \calo_{\mfp,q}
  :=
  \calo_\mfp \otimes_\Z \Z_q.
\]

\begin{lemma}
\label{lem:local-oriented-maximal}
  If $q \neq p$, every maximal order in $B_{p,q}$ is
  $B_{p,q}^\times$-conjugate to
  $\calo_{\mfp,q}$, and
  \[
    N_{B_{p,q}^\times}
    \paren{\calo_{\mfp,q}}
    =
    \Q_q^\times
    \calo_{\mfp,q}^\times.
  \]
  At $q = p$, every oriented maximal order is
  $B_{p,p}^\times$-conjugate to
  $(\calo_{\mfp,p},\sigma_\mfp)$, and its
  orientation-preserving stabilizer is
  \[
    \Q_p^\times
    \calo_{\mfp,p}^\times.
  \]
\end{lemma}

\begin{proof}
  If $q \neq p$, then
  $B_{p,q} \iso M_2(\Q_q)$. The assertions follow
  from the conjugacy of lattices in $\Q_q^2$ and
  \[
    N_{\GL_2(\Q_q)}
    \paren{M_2(\Z_q)}
    =
    \Q_q^\times\GL_2(\Z_q).
  \]

  At $p$, the division algebra $B_{p,p}$ has a unique
  maximal order. Conjugation by a quaternionic
  uniformizer induces the nontrivial automorphism of
  its residue field $\F_{p^2}$. Hence an element
  preserves the orientation precisely when its
  quaternionic valuation is even. The subgroup of
  elements of even valuation is
  $\Q_p^\times \calo_{\mfp,p}^\times$.
\end{proof}

For $g \in \wh B_p^\times$, set
\begin{align*}
  \calo_g
  &:= B_p \cap
    g\wh\calo_\mfp g^{-1},
  \intertext{and for $x \in \calo_g$, define}
  \sigma_g(x)
    &:=
      \sigma_\mfp
      \paren{g_p^{-1}xg_p}.
\end{align*}
\begin{prop}
\label{prop:oriented-maximal}
  There is a natural bijection
  \[
    \begin{tikzcd}[row sep = tiny]
    \Pic(B_p)
      \arrow[r, "\widesim{}"] &
    \Pic^{\mathrm{or}}(B_p), \\
    \sbrac{g}
      \arrow[r, mapsto] &
      \sbrac{\paren{\calo_g,\sigma_g}},
    \end{tikzcd}
  \]
\end{prop}

\begin{proof}
  Lemma~\ref{lem:local-oriented-maximal} gives
  surjectivity place by place. At $p$, the local
  conjugating element may be multiplied by a
  quaternionic uniformizer when necessary to match the
  chosen orientation.

  For injectivity, suppose that the oriented orders
  attached to $g_1$ and $g_2$ are equivalent through
  $b \in B_p^\times$. Then
  \[
    g_2^{-1}bg_1
  \]
  stabilizes
  $(\wh\calo_\mfp,\sigma_\mfp)$ at every finite
  place. Lemma~\ref{lem:local-oriented-maximal}
  therefore gives
  \[
    g_2^{-1}bg_1
    \in
    \wh\Q^\times
    \wh\calo_\mfp^\times.
  \]
  Hence $g_1$ and $g_2$ represent the same adelic
  double coset.
\end{proof}

We use
Proposition~\ref{prop:oriented-maximal} to identify
these two descriptions of $\Pic(B_p)$. The inclusion
\[
  \wh K^\times
  \longhookrightarrow
  \wh B_p^\times
\]
and the optimal embedding
$\calo_c \longhookrightarrow \calo_\mfp$ induce
\[
  \iota_\mfp:
  \Pic(\calo_c)
  \longrightarrow
  \Pic(B_p).
\]
Applying the projection construction of
Section~\ref{subsec:hecke-operators} using this
oriented embedding gives
$\sbrac{\one}_\mfp$
and
$\sbrac{\xi}_\mfp$.

\subsection{Recovery of the trace identity of
Darmon--Harris--Rotger--Venkatesh}
\label{subsec:dhrv-recovery}

Before returning to elliptic units, we recover
\cite[Theorem~2.2]{dhrv} in our notation and
normalizations. In the unramified setting below, this
theorem compares, through the Shimura-class pairing, the
trace from level $Dp$ to level $p$ of
$\theta_{\psi_1}(pz)\theta_{\psi_2}(z)$ with four times
the quaternionic theta lift attached to
$\psi_1\psi_2$ and $\psi_1\psi_2^{c_K}$. When $D$ is
prime, the comparison is an equality of
characteristic-zero modular forms. The precise statement
is Proposition~\ref{prop:dhrv-trace}.

The argument has two stages. We first identify the
standard theta kernel of Section~\ref{sec:theta-two}
with the quaternionic correspondence used in
\cite{dhrv} and recover the full two-character trace
identity. We then specialize to the dual pair relevant
to the Harris--Venkatesh conjecture and use
Proposition~\ref{prop:theta-opt} to identify the
quaternionic theta lift with the optimal form.

\subsubsection*{The DHRV setup}

For this subsection only, put
$D:=\verts{\disc(K)}$ and $C:=\Pic(\calo_K)$, and
assume that $D$ is odd. Let $p \geq 5$ be such that $-p$
is a square in $(\Z/D\Z)^\times$. Thus $p$ is inert in
$K$; for composite $D$, this is stronger than inertness.
We use $p$ for the prime denoted by $N$ in
\cite[Section~2]{dhrv}.

Choose an odd prime $q_0\equiv-p\pmod D$ such that
$q_0\calo_K=\mfq_0\overline{\mfq}_0$ and both primes on
the right are principal. Such primes exist by the
Chebotarev density theorem; see
\cite[Section~2.1]{dhrv}.
Let $B_p$ be the definite quaternion algebra of reduced
discriminant $p$, written as $B_p = K \oplus Kj$, where
$j^2 = -pq_0$ and $jx = \overline{x}j$ for $x \in K$.
Put $\calo(q_0):=\calo_K\oplus\calo_Kj$, and fix an
oriented maximal order $\calm\supset\calo(q_0)$.
For $I\in C$, write
$\iota_{\calm}(I):=I^{-1}\calm I$ for its
image in the oriented class set, with orientation
transported from $\calm$; this is
\cite[Equation~13]{dhrv}. To compare this construction
with Proposition~\ref{prop:theta-opt}, we first verify
that $\calm$ satisfies the optimality condition used
there.

\begin{lem}
\label{lem:dhrv-order-optimal}
  For every class-group character $\xi$ of $K$, the
  order $\calm$ is $\xi$-optimal in the sense of
  Definition~\ref{defn-xi-opt}.
\end{lem}

\begin{proof}
  The character $\xi$ is unramified, so
  $c(\xi) = 1$ and $\calo_{c(\xi)} = \calo_K$.
  Since $\calo_K \subset \calm$ and $\calo_K$ is
  maximal in $K$, we have that $K \cap \calm = \calo_K$.

  Let $r \nmid D$. If $r \neq p$, identify
  $B_{p,r} \iso \End_{\Q_r}(K_r)$
  so that $K_r$ acts by multiplication. There is a
  $\Z_r$-lattice $\Lambda_r \subset K_r$ such that
  $\calm_r = \End_{\Z_r}(\Lambda_r)$.
  Because $\calm_r$ contains $\calo_{K,r}$, the
  lattice $\Lambda_r$ is a fractional
  $\calo_{K,r}$-ideal. The quadratic
  $\Z_r$-algebra $\calo_{K,r}$ is \'{e}tale, so every such
  ideal is principal. Write
  $\Lambda_r = a_r\calo_{K,r}$
  with $a_r \in K_r^\times$.
  Let $c_{K,r}$ be the nontrivial automorphism of
  $K_r/\Q_r$, viewed as a $\Q_r$-linear endomorphism of
  $K_r$, and set $j_r := a_r c_{K,r} a_r^{-1}$.
  Conjugating
  $\End_{\Z_r}(\calo_{K,r})$ by $a_r$ gives
  \[
    \calm_r
    =
    \calo_{K,r}
    \oplus
    \calo_{K,r}j_r,
  \]
  with $j_rx = \overline{x}j_r$
  and $j_r^2 = 1$.
  Thus $j_r^2\Z_r = \Z_r$, the local reduced
  discriminant ideal of the split maximal order
  $\calm_r$.

  At $r = p$, the extension $K_p/\Q_p$ is unramified
  and $B_{p,p}$ is the division quaternion algebra. Its
  unique maximal order is
  \[
    \calm_p
    =
    \calo_{K,p}
    \oplus
    \calo_{K,p}\Pi_p,
  \]
  where one may choose $\Pi_p$ so that
  $\Pi_px = \overline{x}\Pi_p$
  and
  $\Pi_p^2 = p$.
  Hence $\Pi_p^2\Z_p = p\Z_p$, the local reduced
  discriminant ideal of $\calm_p$. No local
  condition is required at primes dividing $D$.
\end{proof}

We first specify the Shimura-class functional used
below. Let $L$ contain the values of the characters
under consideration, let $\ell \geq 5$ divide $p - 1$,
write $\ell^t \Vert p - 1$, and set
$R := \calo_L \otimes_\Z \Z_{(\ell)}$.
For $F \in M_2\paren{\Gamma_0(p), R}$, define
\[
  \log_\ell\mfs_p\paren{F}
  :=
  \pair{
    \overline F,
    (1 \otimes \log_\ell)\mfs_p}
  \in
  R/\ell^tR.
\]
Here the bar denotes reduction modulo $\ell^t$,
$(1 \otimes \log_\ell)\mfs_p$ is the coherent Shimura
class after scalar extension, and the pairing is the
Serre-duality pairing of \cite[Equation~11]{dhrv}.
For cuspidal $F$, this agrees with applying
$\log_\ell$ to the $\Delta_p$-valued functional
$\mfs_p$ defined in the introduction.

\subsubsection*{Comparison of theta normalizations}

Let $\cale := \Pic^{\mathrm{or}}(B_p)$.
The proof of Proposition~\ref{prop:oriented-maximal},
with $\calm$ as base point, identifies $\cale$ with the
corresponding adelic class set.
We use the identification
\[
  C
  \iso
  K^\times\bs\wh K^\times/
  \wh\calo_K^\times
\]
that sends the class of $t$ to the ideal class of
$t^{-1}\wh\calo_K\cap K$. For $x \in \cale$,
let $\calm_x$ be the underlying maximal order
and set $w_x := \verts{\calm_x^\times}/2$.
Let $e_x$ be the characteristic function of $x$ in $\Div(\cale)$,
and let $\delta_x$ be the corresponding characteristic
automorphic function.
Set $\wh{\calm}^{\,1} := \wh{\calm}^\times
\cap B_{\A^\infty}^1$. We normalize the measure on the
quaternionic class set by $\vol(\set{x})=1/w_x$, the
Haar measure on $\wh{\calm}^{\,1}$ by
$\vol\paren{\wh{\calm}^{\,1}}=1$, and the toric quotient
in Equation~\ref{eq:proj-complex} so that it induces
counting measure on $C$. With these choices, the adelic
pairing and toric projection agree with the
normalizations in \cite[Equations~(14)--(17)]{dhrv}.

Give $\Div(\cale)$ the DHRV pairing
\[
  \pair{e_x,e_y}_{\mathrm{DHRV}}
  :=
  w_x\delta_{x,y}.
\]
Let $\cala_{\calm}$ be the quaternionic automorphic
space of trivial central character and level
$\wh{\calm}^\times$. For
$\varphi \in \cala_{\calm}$, set
\[
  \operatorname{mass}_{\calm}(\varphi)
  :=
  \int_{\sbrac{B_p^\times/\Q^\times}}
    \varphi(x)\,dx.
\]
Let $\Phi_{\calm}$ be the standard Schwartz
function of Definition~\ref{def:standard-Eichler}.
For $n \geq 1$, let
$n^\infty \in \A^{\infty,\times}$ be the finite idele
induced by $n$, and set
$\rmt_n :=
\rmt_{\Phi_{\calm}^\infty}
\paren{n^\infty}$.
Let $\Theta_{\calm}^{\std}$ be the modular avatar
of
$\theta\paren{g, \varphi_1 \otimes \varphi_2,
\Phi_{\calm}}$.
The calculation preceding
Equation~\ref{eq:Theta-standard} gives
\begin{equation}
\label{eq:dhrv-theta-standard}
\begin{split}
  \Theta_{\calm}^{\std}
  \paren{\varphi_1 \otimes \varphi_2}
  &=
  \frac{1}{2}
  \operatorname{mass}_{\calm}(\varphi_1)
  \operatorname{mass}_{\calm}(\varphi_2)
  + \sum_{n \geq 1}
    \pair{\varphi_1,\rmt_n\varphi_2}q^n.
\end{split}
\end{equation}

\begin{lem}
\label{lem:dhrv-normalizations}
  The map
  \[
    \begin{tikzcd}[row sep = tiny]
      \jmath_{\calm}:
      \Div(\cale)\otimes\C
      \arrow[r] &
      \cala_{\calm}, \\
      e_x \arrow[r, mapsto] & w_x\delta_x,
    \end{tikzcd}
  \]
  is an isometric, Hecke-equivariant isomorphism. It
  sends
  \[
    \Sigma_0
    :=
    \sum_{x \in \cale}\frac{e_x}{w_x}
  \]
  to the constant function $\one$, and it sends the
  vector $\sbrac{\psi}$ of
  \cite[Equation~17]{dhrv} to the vector
  $\sbrac{\psi}_{\calm}$ defined by
  Equation~\ref{eq:proj-complex}. Consequently, the
  standard lift of
  Equation~\ref{eq:dhrv-theta-standard} is
  exactly the correspondence $\Theta$ of
  \cite[Equation~16]{dhrv}.
\end{lem}

\begin{proof}
  The measure normalization gives
  \[
    \pair{\delta_x,\delta_y}
    =
    \frac{\delta_{x,y}}{w_x}.
  \]
  Hence $\jmath_{\calm}$ is an isometry.

  For Hecke equivariance, write
  \[
    T_ne_x
    =
    \sum_y b_{x,y}(n)e_y.
  \]
  Right convolution on characteristic automorphic
  functions is given by the transposed Brandt matrix:
  \[
    \rmt_n\delta_x
    =
    \sum_y b_{y,x}(n)\delta_y.
  \]
  The self-adjointness of $T_n$ for the DHRV pairing
  gives
  \[
    b_{x,y}(n)w_y
    =
    b_{y,x}(n)w_x.
  \]
  Therefore
  \begin{align*}
    \rmt_n\paren{w_x\delta_x}
    &=
    \sum_y w_xb_{y,x}(n)\delta_y \\
    &=
    \sum_y b_{x,y}(n)w_y\delta_y \\
    &=
    \jmath_{\calm}\paren{T_ne_x}.
  \end{align*}
  This includes the operator at $p$; see
  \cite[Equation~14 and Remark~2.1]{dhrv}.

  The assertion concerning $\Sigma_0$ is immediate.
  For the toric vectors, choose a finite idele $t_I$
  representing $I$ under the convention above. Thus
  \[
    t_I\wh\calo_K = \wh I^{-1},
    \qquad
    \psi(t_I) = \psi(I).
  \]
  Then $t_I$ maps to $\iota_{\calm}(I)$. Applying
  Equation~\ref{eq:proj-complex} to $\delta_x$ gives
  \[
    \frac{\sbrac{\psi}_{\calm}(x)}{w_x}
    =
    \sum_{\substack{I \in C\\
      \iota_{\calm}(I) = x}}
    \psi(I).
  \]
  Thus
  \[
    \sbrac{\psi}_{\calm}
    =
    \jmath_{\calm}
    \paren{
      \sum_{I \in C}
      \psi(I)e_{\iota_{\calm}(I)}}.
  \]
  The vector in parentheses is the vector of
  \cite[Equation~17]{dhrv}.

  Finally, for $v \in \Div(\cale)$,
  \[
    \operatorname{mass}_{\calm}
    \paren{\jmath_{\calm}(v)}
    =
    \pair{v,\Sigma_0}_{\mathrm{DHRV}}.
  \]
  The isometry and Hecke equivariance identify the
  constant term and every positive Fourier coefficient
  of Equation~\ref{eq:dhrv-theta-standard} with
  \cite[Equation~16]{dhrv}. This proves the final
  assertion.
\end{proof}

For a full $\Z$-lattice $\Lambda$ in $K$ or $B_p$,
let $\mathfrak n(\Lambda)$ be the positive generator of
the fractional ideal of $\Q$ spanned by
$\set{\rmn(x) \mid x \in \Lambda}$, where $\rmn$ is the
field norm on $K$ and the reduced norm on $B_p$. Put
\[
  \theta_\Lambda(z)
  :=
  \sum_{x \in \Lambda}
    q^{\rmn(x)/\mathfrak n(\Lambda)}.
\]
This is the normalization of
\cite[Equation~21]{dhrv}.
For $I \in C$, set
$\delta_I^{\calm} := \jmath_{\calm}\paren{e_{\iota_{\calm}(I)}}$.
Lemma~\ref{lem:dhrv-normalizations} translates
\cite[Proposition~2.4]{dhrv} into the identity
\begin{equation}
\label{eq:dhrv-basis-theta}
  \theta_{\overline{J}\calm I}
  =
  2\Theta_{\calm}^{\std}
  \paren{
    \delta_I^{\calm}
    \otimes
    \delta_J^{\calm}}
\end{equation}
for all $I,J \in C$. Thus the standard quaternionic theta
kernel used in Proposition~\ref{prop:theta-opt} recovers
the theta-correspondence input of the proof of
\cite[Theorem~2.2]{dhrv}.

\subsubsection*{The trace identity}

For a character
$\psi:C \longrightarrow L^\times$, put
\begin{equation}
\label{eq:dhrv-class-theta}
  \theta_\psi
  :=
  \sum_{I \in C}\psi(I)^{-1}\theta_I.
\end{equation}
This is the convention of
\cite[Equation~22]{dhrv}. For a character $\psi$ of
$C$, write $\psi^{c_K}(I) := \psi(\overline{I})$.
Since $\overline{I}=I^{-1}$ in $C$, one has
$\psi^{c_K}=\psi^{-1}$.

Write
\[
  \Tr_p^{Dp}:
  M_2\paren{\Gamma_0(Dp)}
  \longrightarrow
  M_2\paren{\Gamma_0(p)}
\]
for the usual trace map.
The proof uses three inputs from \cite{dhrv}.
\cite[Proposition~2.3]{dhrv} expresses the trace of a
product of binary theta series as a sum of theta series
of quaternionic lattices. Equation~
\ref{eq:dhrv-basis-theta}, which is
\cite[Proposition~2.4]{dhrv} in our normalization,
identifies each lattice theta series with a value of the
standard quaternionic theta lift. Finally,
\cite[Proposition~2.5]{dhrv} shows that the quadratic
genus twists have the same Shimura-class pairing. We
record the resulting character average in order to make
the factor $4$ explicit.

\begin{prop}[{\cite[Theorem 2.2]{dhrv}}]
\label{prop:dhrv-trace}
  Let $\psi_1, \psi_2$ be characters of $C$, and set
  $\psi_{12} := \psi_1\psi_2$
  and $\psi_{12}' := \psi_1\psi_2^{c_K}$.
  There is a constant
  $\ell_0 = \ell_0(D, \psi_1, \psi_2)$, independent
  of the auxiliary prime $p$, such that, for every prime
  $\ell \geq \ell_0$ dividing $p - 1$, one has
  \begin{equation}
  \label{eq:dhrv-recovered}
    \log_\ell\mfs_p
      \paren{\Tr_p^{Dp}
      \paren{\theta_{\psi_1}(pz) \theta_{\psi_2}(z)}}
    = 4\log_\ell\mfs_p
      \paren{\Theta_{\calm}^{\std}
      \paren{\sbrac{\psi_{12}}_{\calm}
        \otimes
        \sbrac{\psi_{12}'}_{\calm}}}.
  \end{equation}
  If $D$ is prime, the corresponding equality of
  modular forms holds in characteristic zero before
  applying $\log_\ell\mfs_p$; in particular, no
  condition on a coefficient prime $\ell$ is needed.
\end{prop}

\begin{proof}
  Let $a$ be the number of prime divisors of $D$, let
  $C[2] := \set{I \in C \Mid I^2 = 1}$ be the $2$-torsion
  subgroup of $C$, and set
  $\calx := \paren{C/C^2}^\vee$.
  For $\nu \in \calx$, define
  \begin{align*}
    G_\nu
      &:=
      \theta_{\psi_1^{-1}\nu}(pz)
      \theta_{\psi_2^{-1}\nu}(z),
      \\
    G &:= \sum_{\nu \in \calx}G_\nu.
  \end{align*}
  Expanding Equation~\ref{eq:dhrv-class-theta} and
  using character orthogonality gives
  \begin{equation}
  \label{eq:dhrv-average-expansion}
    G =
      \sum_{I_1,I_2 \in C}
      \psi_{12}(I_1)
      \psi_{12}'(\overline{I_2})
      \theta_{I_1I_2}(z)
      \theta_{I_1\overline{I_2}}(pz).
  \end{equation}
  Indeed, the sum of $\nu(IJ)$ over
  $\nu \in \calx$ vanishes unless $IJ \in C^2$, in which
  case it is $\verts{\calx}=\verts{C[2]}$. The map
  \[
    \begin{tikzcd}[row sep = tiny]
      C \times C
        \arrow[r] &
      C \times C, \\
      (I_1, I_2)
        \arrow[r, mapsto] &
      (I_1\overline{I_2}, I_1I_2)
    \end{tikzcd}
  \]
  has kernel naturally identified with $C[2]$. Thus
  every pair $(I,J)$ with $IJ \in C^2$ occurs in exactly
  $\verts{C[2]}$ ways. This gives
  \cite[Equation~27]{dhrv}.

  For each prime $r \mid D$, let $\mathfrak r$ be the
  unique prime of $\calo_K$ above $r$, and retain the
  chosen principal prime $\mfq_0$ above $q_0$. For a
  divisor $d \mid Dq_0$, define
  \[
    \mathfrak d
    :=
    \prod_{\substack{r \mid d \\ r \mid D}}
      \mathfrak r
    \cdot
    \begin{cases}
      \mfq_0 & \text{if $q_0 \mid d$,} \\
      \calo_K & \text{if $q_0 \nmid d$,}
    \end{cases}
  \]
  and set
  $\calm_d :=
  \mfd^{-1} \calm \mfd$.
  These are the ideals and maximal orders of
  \cite[Equation~20]{dhrv}.
  Applying \cite[Proposition~2.3]{dhrv} term by term to
  Equation~\ref{eq:dhrv-average-expansion} gives
  \begin{align}
  \label{eq:dhrv-trace-prop-23}
    \Tr_p^{Dp}(G)
    &=
    \frac{1}{2}
    \sum_{I_1,I_2 \in C}
    \psi_{12}(I_1)
    \psi_{12}'(\overline{I_2})
    \sum_{d \mid Dq_0}
    \theta_{I_1\calm_d I_2}.
  \end{align}

  For a fixed $d$, use
  \[
  I_1\calm_d I_2
  =
  \mathfrak d^{-1}I_1\calm I_2\mathfrak d,
  \]
  set $J_1 := \mathfrak d^{-1}I_1$ and
  $J_2 := I_2\mathfrak d$, and then rename
  $(J_1,J_2)$ as $(I_1, I_2)$.
  If $d \mid D$, then
  $\overline{\mathfrak d} = \mathfrak d$ and the class of
  $\mathfrak d$ has order at most $2$. One has
  $\psi_{12}'(\mathfrak d)
  = \psi_{12}(\mathfrak d)$, and this common value has
  order at most $2$. Hence the character factor
  introduced by the change of variables is
  \[
    \psi_{12}(\mathfrak d)
    \psi_{12}'
    \paren{(\overline{\mathfrak d})^{-1}}
    = \psi_{12}(\mathfrak d)^2
    = 1.
  \]
  The ideal $\mfq_0$ is principal, so its factor is
  also $1$; multiplicativity treats every divisor of
  $Dq_0$. Since $D$ is an odd fundamental
  discriminant, it is squarefree; hence $Dq_0$ has
  $2^{a+1}$ positive divisors. Equation~
  \ref{eq:dhrv-trace-prop-23} becomes
  \begin{equation}
  \label{eq:dhrv-traced-average}
    \Tr_p^{Dp}(G)
    = 2^a
      \sum_{I_1,I_2 \in C}
      \psi_{12}(I_1)
      \psi_{12}'(\overline{I_2})
      \theta_{I_1\calm I_2}.
  \end{equation}

  Apply Equation~\ref{eq:dhrv-basis-theta} with
  $I := I_2$ and $J := \overline{I_1}$. Then
  \[
    \theta_{I_1\calm I_2}
    =
    2\Theta_{\calm}^{\std}
    \paren{
      \delta_{I_2}^{\calm}
      \otimes
      \delta_{\overline{I_1}}^{\calm}}.
  \]
  Moreover,
  \begin{align*}
    \sum_{I_2 \in C}
      \psi_{12}'(\overline{I_2})
      \delta_{I_2}^{\calm}
    &= \sbrac{(\psi_{12}')^{-1}}_{\calm}, \\
      \sum_{I_1 \in C}
        \psi_{12}(I_1)
        \delta_{\overline{I_1}}^{\calm}
    &= \sbrac{\psi_{12}^{-1}}_{\calm}.
  \end{align*}
  Using the symmetry of
  $\Theta_{\calm}^{\std}$ therefore gives
  \begin{equation}
  \label{eq:dhrv-average-theta}
    \Tr_p^{Dp}(G)
    =
    2^{a+1}
    \Theta_{\calm}^{\std}
    \paren{
      \sbrac{\psi_{12}^{-1}}_{\calm}
      \otimes
      \sbrac{(\psi_{12}')^{-1}}_{\calm}}.
  \end{equation}

  By genus theory,
  \[
    \verts{\calx}
    =
    \verts{C[2]}
    =
    2^{a - 1};
  \]
  see \cite[paragraph following
  Equation~29]{dhrv}. By
  \cite[Proposition~2.5]{dhrv}, once $\ell$ is
  sufficiently large relative to $D$, the forms
  $\Tr_p^{Dp}(G_\nu)$ have the same pairing with the
  Shimura class. Since
  $G = \sum_{\nu \in \calx}G_\nu$, it follows that
  \[
    \log_\ell\mfs_p\paren{\Tr_p^{Dp}(G)}
      = 2^{a - 1}
        \log_\ell\mfs_p
        \paren{\Tr_p^{Dp}(G_\one)}.
  \]
  Since $\ell \geq 5$, the scalar $2^{a - 1}$ is
  invertible in $R/\ell^tR$. Combining this equality
  with Equation~\ref{eq:dhrv-average-theta}, dividing
  by $2^{a - 1}$, and replacing
  $(\psi_1, \psi_2)$ by
  $(\psi_1^{-1}, \psi_2^{-1})$ proves
  Equation~\ref{eq:dhrv-recovered}.

  If $D$ is prime, then $a=1$ and
  $\calx=\set{\one}$. Thus $G=G_\one$, and the
  strict identity follows directly from Equation~
  \ref{eq:dhrv-average-theta}, without using
  \cite[Proposition~2.5]{dhrv}.
\end{proof}

After the preceding identifications of notation and
normalization, Proposition~\ref{prop:dhrv-trace} is
\cite[Theorem~2.2]{dhrv}. We now relate its dual-pair
specialization to the optimal form.

\begin{cor}
\label{cor:dhrv-optimal}
  Let $\chi$ be an unramified class-group character
  such that $\theta_\chi$ is cuspidal,
  let $f := \theta_\chi$,
  let $\xi := \chi^{1 - c_K}$,
  and retain the hypotheses of
  Proposition~\ref{prop:dhrv-trace}. Then
  \begin{equation}
  \label{eq:dhrv-optimal-comparison}
    \log_\ell\mfs_p
    \paren{
      \Tr_p^{Dp}
      \paren{f(z)f^*(pz)}}
    = 4c_{\chi,\calm}^{-1}
      \log_\ell\mfs_p
      \paren{f^\opt(z, pz)}.
  \end{equation}
  If $D$ is prime, the corresponding equality of
  modular forms holds in characteristic zero before
  applying $\log_\ell\mfs_p$.
\end{cor}

\begin{proof}
  Apply Proposition~\ref{prop:dhrv-trace} with
  $\psi_1:=\chi$ and $\psi_2:=\chi^{-1}$. Then
  $\psi_{12}=\one$ and
  $\psi_{12}'=\chi(\chi^{-1})^{c_K}
  =\chi^{1-c_K}=\xi$.
  Complex conjugation of Fourier coefficients gives
  $f^* = \theta_{\chi^{-1}}$. The change of variables
  $I\mapsto \overline{I}$ in
  Equation~\ref{eq:dhrv-class-theta} also gives
  $\theta_\chi = \theta_{\chi^{-1}}$. Hence the
  left-hand side of Equation~\ref{eq:dhrv-recovered}
  is the left-hand side of
  Equation~\ref{eq:dhrv-optimal-comparison}.

  By Lemma~\ref{lem:dhrv-order-optimal}, the order
  $\calm$ is $\xi$-optimal. Applying
  Proposition~\ref{prop:theta-opt} and repeating the
  passage to modular avatars preceding
  Equation~\ref{eq:theta-opt-standard} gives
  \[
    f^\opt(z,pz)
    =
    c_{\chi,\calm}
    \Theta_{\calm}^{\std}
    \paren{
      \sbrac{\one}_{\calm}
      \otimes
      \sbrac{\xi}_{\calm}}.
  \]
  Substitution proves Equation~
  \ref{eq:dhrv-optimal-comparison}.
\end{proof}

\begin{remark}
\label{rem:dhrv-generalization}
  Proposition~\ref{prop:dhrv-trace} is the full
  two-character identity, whereas
  Corollary~\ref{cor:dhrv-optimal} is its dual-pair
  specialization. The factor $4$ compares the traced
  product with the standard quaternionic theta lift,
  while $c_{\chi,\calm}$ compares that lift with the
  intrinsically normalized optimal form.
  For ramified $\chi$, Theorem~\ref{thm:opt-unique} supplies
  the dual-pair theta identity for the optimal form,
  while Theorem~\ref{thm:hv-reform} transfers its
  Harris--Venkatesh period to the newform. In this sense,
  the optimal-form construction extends the role of the
  dual-pair case of \cite[Theorem~2.2]{dhrv} to arbitrary
  conductor and removes the odd-discriminant hypothesis.
\end{remark}

We now return to the general notation and hypotheses of
Section~\ref{sec:proof-opt}.

%%%%%%%%%%%%%%%%%%%%%%%%%%%%%%%%%%%%%%%%%%%%%%%%%%%%%%%%%%%%%%%%

\subsection{Elliptic units}
\label{sec:elliptic-units}

In the imaginary dihedral setting,
$\Ad(\rho_f)\iso\eta\oplus\Ind_{G_K}^{G_\Q}(\xi)$,
where $\eta$ is the quadratic character of $K/\Q$.
Choose a basis $v_1,v_2$ of the induced representation
over $\Q(\chi)$ such that
$v_2=\rho_f(\widetilde c_K)v_1$. Then, for
$h\in\gal(E/K)$,
\[
  \rho_f(h)
  =
  \begin{pmatrix}
    \chi(h)&0\\
    0&\chi^{c_K}(h)
  \end{pmatrix},
  \qquad
  \rho_f(\widetilde c_K)
  =
  \begin{pmatrix}
    0&1\\
    1&0
  \end{pmatrix}.
\]
We use the $\xi^{-1}$-isotypical unit line
\[
  \calu_\xi
  :=
  \set{u\in\calo_{H_c}^\times\otimes\Z[\xi]
  \Mid
  u^\sigma=\xi(\sigma)^{-1}u
    \text{ for every }
  \sigma\in\gal(H_c/K)}.
\]

We next make the reduction maps precise. Since $p$ is
inert in $K$, the prime $p\calo_K$ is principal and has
trivial Artin symbol in $\gal(H_c/K)$; hence
$\calo_{H_c}/\mfp_c\iso\F_{p^2}$. The discrete logarithm
$\log_\ell:\F_p^\times\to\Z/\ell^t\Z$ extends uniquely
to a Frobenius-invariant map
$\log_{\ell,\mfp_c}:\F_{p^2}^\times\to\Z/\ell^t\Z$,
because $\F_{p^2}^\times/\F_p^\times$ has order $p+1$,
which is prime to $\ell$. We use this extension after
reducing ring-class units modulo $\mfp_c$.

For
$u\in\calu(\Ad(\rho_f))\otimes_{R_f}\Q(\chi)$, set
$x_\mfp:=2\rho_f(\Frob_\mfp)
-\Tr(\rho_f(\Frob_\mfp))\Id_{2\times2}$. Equivariance
implies that $u(x_\mfp)$ is fixed by $\Frob_\mfp$, so
its reduction lies in $\F_p^\times\otimes\Q(\chi)$.
Define
$\log_\ell\Reg_\mfp(u):=
\log_\ell(u(x_\mfp)\bmod\mfp)$. Thus no norm between
residue fields occurs in the adjoint regulator. The map
$\log_{\ell,\mfp_c}$ is used only for the intermediate
ring-class units; cf. \cite[Lemma~5.6]{dhrv} and
\cite[Theorem~2.5]{lecouturier-hv}.

Elliptic units are analogues of cyclotomic units
for imaginary quadratic fields
and were introduced by Robert \cite{robert}.
In general, they are constructed from singular values of modular functions
and defined in arbitrary class fields
(see \cite[Chapter 11]{kubert-lang}).
We will use a particular class of elliptic units
in $\calo_{H_c}^\times$ whose modulo-$p$ reduction was
studied by Darmon--Harris--Rotger--Venkatesh \cite[Section 5.1]{dhrv}.

\subsubsection*{The elliptic unit of Darmon--Harris--Rotger--Venkatesh}

We now recall the construction given in \cite[Section~5.1]{dhrv}
for the elliptic unit associated to $\xi$ and $\mfl$,
which lies in the ring class field $H_c$.
Consider the modular unit
$u_\lambda(z):=\Delta(z)/\Delta(\lambda z)$ on
$Y_0(\lambda)$; it is denoted $\Delta_N$ on $Y_0(N)$ in
\cite[Section~4.4]{dhrv}. Here
$\Delta(z)=q\prod_{n\geq1}(1-q^n)^{24}$ is the
Ramanujan delta function, with $q=e^{2\pi iz}$.
The curve $Y_0(\lambda)$ parametrizes degree-$\lambda$
isogenies $E_1\to E_2$, and the two projections are
$\pi_i(E_1\to E_2):=[E_i]$ for $i\in\set{1,2}$.

Let $Z(1) \subset X(1)$ be the set of isomorphism classes of
elliptic curves $E$ with $\End(E) = \calo_c$ and let
$Z_0(\lambda) \subset X_0(\lambda)$
be the subset consisting of points $\phi: E_1 \lra E_2$ with both
$\End(E_i) = \calo_c$. Then all points of $Z(1)$ and $Z_0(\lambda)$ are
defined over $H_c$. In particular, $u_\lambda(x) \in H_c^\times$ for all
$x \in Z_0(\lambda)$. We fix one point $x_c = \C / \calo_c \in Z(1)$.

The maps $\pi_i$ induce projections
$Z_0(\lambda)\to Z(1)$. Fixing the splitting
$\lambda=\mfl\overline{\mfl}$ gives a section
$\eta_\lambda:Z(1)\to Z_0(\lambda)$ defined by
$\eta_\lambda(E):=(E\to E/E[\mfl])$.
Darmon--Harris--Rotger--Venkatesh \cite[Definition 5.1]{dhrv}
define the elliptic unit
\begin{equation}
	\label{eq:uxip}
	u_{\xi, \lambda} := \sum_{\sigma \in \gal\paren{H_c / K}}
		u_\lambda\paren{\eta_\lambda\paren{x_c}}^\sigma \otimes
    \xi(\sigma) \in H_c^\times \otimes \Z[\xi].
\end{equation}
If $\xi$ is nontrivial, then $u_{\xi,\lambda}$ belongs
to $\calu_\xi\subset
\calo_{H_c}^\times\otimes\Z[\xi]$; see
\cite[Equation~103]{dhrv}.

\subsubsection*{An auxiliary-free elliptic unit}

By the Chebotarev density theorem, we may choose a split
prime $\lambda=\mfl\overline{\mfl}$, coprime to $c$,
such that $\xi(\mfl)$ generates $\Im(\xi)$. Retain the
integer $m(\xi)$ defined in the introduction. The
cyclotomic norm formula gives
$m(\xi)=\rmn_{\Q(\xi)/\Q}
(1-\xi(\overline{\mfl}))$.
Define
\begin{equation}
\label{eq:uxi}
  u_\xi := \frac{m(\xi)} {1 - \xi(\overline{\mfl})} u_{\xi,\lambda}.
\end{equation}

The scalar
$m(\xi)/(1-\xi(\overline{\mfl}))$ in
Equation~\ref{eq:uxi} belongs to $\Z[\xi]$. If
$\verts{\Im(\xi)}$ is not a prime power, then
$1-\xi(\overline{\mfl})$ is a unit in $\Z[\xi]$. If it
is a power of a prime $v$, then its norm is
$v=m(\xi)$, and the quotient is the product of the
remaining Galois conjugates of
$1-\xi(\overline{\mfl})$. Consequently,
$u_\xi\in\calu_\xi$.

We now prove that the unit $u_\xi$
is independent of the auxiliary prime,
making it slightly more natural
than the elliptic unit in Equation \ref{eq:uxip}.

\begin{prop}
\label{prop:uxip}
  The unit $u_\xi$ is independent of the choice of split
  prime $\lambda = \mfl\overline{\mfl}$ satisfying the
  preceding conditions.
\end{prop}

\begin{proof}
  It remains only to prove independence of $\lambda$.
	By the definition of $u_\xi$, we need to show that for any two
	primes $\lambda_1 \neq \lambda_2$ split in $K$,
	\[
		\paren{1 - \xi\paren{\overline{\mfl}_2}} u_{\xi, \lambda_1} =
			\paren{1 - \xi\paren{\overline{\mfl}_1}} u_{\xi, \lambda_2}
	\]
	where $\mfl_i$ is an invertible ideal in $\calo_c$
	and $\lambda_i \calo_c = \mfl_i \cdot \overline{\mfl}_i$.

	Consider the commutative diagram of isogenous elliptic curves,
	\[
		\begin{tikzcd}
			\C / \calo_c \arrow[r, "x_1"] \arrow[d, "x_2"]
				& \C / \mfl_1^{-1} \arrow[d, "x_4"] \\
			\C / \mfl_2^{-1} \arrow[r, "x_3"]
				& \C / \paren{\mfl_1\mfl_2}^{-1}.
		\end{tikzcd}
	\]
	By construction, the diagonal isogeny has square-free degree
	$\lambda_1 \lambda_2$ and thus  
	has a cyclic kernel isomorphic to
	$(\Z/\lambda_1\Z) \times (\Z/\lambda_2\Z)$.
	Then this diagonal isogeny defines a point $x$ on
	$X_0(\lambda_1\lambda_2)$ and
	is represented by a point $z \in \calh$ in the sense that 
	$x$ is represented by the $\lambda_1\lambda_2$-multiplication map,
	\[
		x:
			\C/(\Z + \Z z) \lra \C/(\Z + \Z \lambda_1\lambda_2 z).
	\]
	The two isogenies $x_1$ and $x_2$ are given by modular subgroups
	of $\ker(x)$ of order $\lambda_1$ and $\lambda_2$ respectively.
	This implies the representations,
	\begin{align*}
		x_1: \C/(\Z + \Z z) &\lra \C/(\Z + \Z \lambda_1 z), \\
		x_2: \C/(\Z + \Z z) &\lra \C/(\Z + \Z \lambda_2 z), \\
		x_3: \C/(\Z + \Z \lambda_2 z) &\lra \C/(\Z + \Z \lambda_1\lambda_2 z), \\
		x_4: \C/(\Z + \Z \lambda_1 z) &\lra \C/(\Z + \Z \lambda_1\lambda_2 z). \\
	\end{align*}

	Then we have points $x_1, x_3 \in X_0(\lambda_1)$ and
	$x_2, x_4 \in X_0(\lambda_2)$ with representatives
	in $\calh$: $z_1 = z_2 = z$, $z_3 = \lambda_2 z_1$,
	and $z_4 = \lambda_1 z_2$.
	By the definition of $u_{\lambda_1}$ and $u_{\lambda_2}$,
	we have the relation,
	\[
		u_{\lambda_1}\paren{x_1} u_{\lambda_2}\paren{x_4}
			= \frac{\Delta(z_1)}{\Delta\paren{\lambda_1 z_1}}
				\frac{\Delta(z_4)}{\Delta\paren{\lambda_2 z_4}}
			= \frac{\Delta(z)}{\Delta\paren{\lambda_1 \lambda_2 z}}
			= \frac{\Delta(z_2)}{\Delta\paren{\lambda_2 z_2}}
				\frac{\Delta(z_3)}{\Delta\paren{\lambda_1 z_3}}
			= u_{\lambda_2}\paren{x_2} u_{\lambda_1}\paren{x_3}.
	\]
	We now use Shimura reciprocity. With the Artin
  reciprocity convention of
  \cite[Equation~101]{dhrv}, the action of an
  invertible $\calo_c$-ideal $\mfb$ on a CM elliptic
  curve is
  \[
    \paren{\C/\mfa^{-1}}^{\Frob(\mfb)}
    \iso
    \C/\paren{\mfa\mfb}^{-1}.
  \]
  This action is compatible with the cyclic isogenies
  induced by inclusions of fractional ideals. Applying
  $\Frob(\mfl_2)$ to
  \[
    x_1:
    \C/\calo_c
    \longrightarrow
    \C/\mfl_1^{-1}
  \]
  therefore gives
  \[
    x_3:
    \C/\mfl_2^{-1}
    \longrightarrow
    \C/\paren{\mfl_1\mfl_2}^{-1}.
  \]
  Similarly, applying $\Frob(\mfl_1)$ to $x_2$
  gives $x_4$. Hence
  \begin{align*}
    x_3
      &= x_1^{\Frob\paren{\mfl_2}}, \\
    x_4
      &= x_2^{\Frob\paren{\mfl_1}}.
  \end{align*}
  Therefore,
	\[
		u_{\lambda_1}\paren{x_1}^{1 - \Frob\paren{\mfl_2}}
			= u_{\lambda_2}\paren{x_2}^{1 - \Frob\paren{\mfl_1}}.
	\]
	
	Now take the $\xi$-sum, as in
  Equation~\ref{eq:uxip}, to obtain
	\[
		\sum_{\sigma \in \gal\paren{H_c/K}}
			u_{\lambda_1}\paren{x_1}^{\paren{1 - \Frob\paren{\mfl_2}} \sigma}
				\otimes \xi(\sigma)
		= \sum_{\sigma \in \gal\paren{H_c/K}}
			u_{\lambda_2}\paren{x_2}^{\paren{1 - \Frob\paren{\mfl_1}}\sigma}
				\otimes \xi(\sigma).
	\]
	Unfolding these sums, we obtain,
	\[
		\paren{1 - \xi\paren{\overline{\mfl}_2}} u_{\xi, \lambda_1}
			= \paren{1 - \xi\paren{\overline{\mfl}_1}} u_{\xi, \lambda_2}.
	\]
\end{proof}

\subsubsection*{The corresponding adjoint unit}

We now separate the fixed and the
orientation-dependent parts of the comparison scalar.
The only factor that varies with the oriented prime
$\mfp$ is the local comparison factor at $p$.

\begin{lemma}
\label{lem:orientation-scalar}
  There is a scalar
  $c_\chi \in \Q(\chi)^\times$, independent of
  $p$ and $\mfp$, such that
  $c_{\chi,\mfp} = c_\chi \, \chi\paren{h_\mfp}$.
\end{lemma}

\begin{proof}
  Let $S$ be the set consisting of the archimedean
  place and the finite primes dividing $N$. For every
  place in $S$, fix once and for all the local
  quadratic-space models and comparison maps used in
  Proposition~\ref{prop:theta-opt}. At every finite
  prime $q \notin S \cup \set{p}$, use the unramified
  split model normalized so that its local character
  factor is $1$. The product of the factors at the
  fixed places is a scalar
  $c_\chi \in \Q(\chi)^\times$
  independent of $p$ and $\mfp$.

  It remains to compute the local factor at $p$. In
  the induced basis fixed above, set
  $J := \rho_f(\widetilde c_K)$,
  so that $Jv_1 = v_2$. The $\mfp_c$-orientation
  identifies the semilinear generator of the local
  quaternionic order with the Frobenius endomorphism
  of the reduction of the CM elliptic curve. Since
  \[
    \Frob_\mfp = \widetilde c_K h_\mfp,
  \]
  one has
  \[
    \rho_f\paren{\Frob_\mfp}v_1
    =
    \chi\paren{h_\mfp}v_2.
  \]
  Thus, relative to the fixed semilinear operator
  $J$, the oriented local generator acts on the
  $\chi$-eigenline with coefficient
  $\chi(h_\mfp)$.

  Under local Artin reciprocity, this coefficient is
  the local character factor
  $\chi_p(a_{\calo_\mfp,p})$ in the comparison scalar
  of Proposition~\ref{prop:theta-opt}. Therefore
  \[
    c_{\chi,\mfp}
    =
    c_\chi\chi\paren{h_\mfp}.
  \]
  This is the orientation calculation appearing in
  Lecouturier~\cite[Theorem~2.5]{lecouturier-hv} (cf. Darmon--Harris--Rotger--Venkatesh
  \cite[Lemma~5.6]{dhrv}).
\end{proof}

Let
\[
  E_{12}
  :=
  \begin{pmatrix}
    0 & 1 \\
    0 & 0
  \end{pmatrix},
  \qquad
  E_{21}
  :=
  \begin{pmatrix}
    0 & 0 \\
    1 & 0
  \end{pmatrix}.
\]
Using the trace pairing
$(A, B) \mapsto \Tr(AB)$ to identify
$\Ad(\rho_f)$ with its dual, define
\begin{equation}
\label{eq:fixed-adjoint-unit}
  \widetilde{u}_\xi
  := u_\xi \otimes E_{12}
    + \widetilde c_K(u_\xi) \otimes E_{21}.
\end{equation}

\begin{lemma}
\label{lem:fixed-adjoint-unit}
  The element $\widetilde{u}_\xi$ belongs to
  \[
    \calu\paren{\Ad(\rho_f)}
    \otimes_{R_f}\Q(\chi)
  \]
  and is independent of $p$ and $\mfp$. Moreover,
  \[
    \log_\ell\Reg_\mfp
    \paren{\widetilde{u}_\xi}
    =
    4\chi\paren{h_\mfp}
    \log_{\ell,\mfp_c}\paren{u_\xi}.
  \]
\end{lemma}

\begin{proof}
  For $h \in \gal(E/K)$, one has
  \begin{align*}
    h\paren{u_\xi \otimes E_{12}}
    &=
    \xi(h)^{-1}u_\xi
    \otimes
    \xi(h)E_{12} \\
    &=
    u_\xi \otimes E_{12}, \\
    h\paren{\widetilde c_K(u_\xi) \otimes E_{21}}
    &=
    \xi(h)\widetilde c_K(u_\xi)
    \otimes
    \xi(h)^{-1}E_{21} \\
    &=
    \widetilde c_K(u_\xi) \otimes E_{21}.
  \end{align*}
  Moreover, $\widetilde c_K$ exchanges the two summands in
  Equation~\ref{eq:fixed-adjoint-unit}. Hence
  $\widetilde{u}_\xi$ is invariant under
  $\gal(E/\Q)$ and belongs to the stated adjoint-unit
  line.

  Since $p$ is inert in $K$,
  $\Tr\paren{\rho_f\paren{\Frob_\mfp}} = 0$.
  Equation~\ref{eq:frob-cK-h} gives
  \[
    \rho_f\paren{\Frob_\mfp}
    =
    \begin{pmatrix}
      0 & \chi^{c_K}\paren{h_\mfp} \\
      \chi\paren{h_\mfp} & 0
    \end{pmatrix}.
  \]
  Furthermore,
  $\widetilde c_K=\Frob_\mfp h_\mfp^{-1}$.
  The Frobenius invariance of $\log_{\ell,\mfp_c}$ and the
  relation
  $u_\xi^{h^{-1}} = \xi(h)u_\xi$ give
  \[
    \log_{\ell,\mfp_c}
    \paren{\widetilde c_K(u_\xi)}
    =
    \xi\paren{h_\mfp}
    \log_{\ell,\mfp_c}\paren{u_\xi}.
  \]
  Evaluating Equation~\ref{eq:fixed-adjoint-unit} at
  $x_\mfp = 2\rho_f(\Frob_\mfp)$ now gives
  \begin{align*}
    \log_\ell\Reg_\mfp
    \paren{\widetilde{u}_\xi}
    &= 2\chi\paren{h_\mfp}
      \log_{\ell,\mfp_c}\paren{u_\xi}
      + 2\chi^{c_K}\paren{h_\mfp}
      \log_{\ell,\mfp_c}
      \paren{\widetilde c_K(u_\xi)} \\
    &= 4\chi\paren{h_\mfp}
      \log_{\ell,\mfp_c}\paren{u_\xi}.
  \end{align*}
\end{proof}

Define the fixed rational adjoint unit
\begin{equation}
\label{eq:u-f-opt}
  u_f^\opt
  :=
  \frac{c_\chi}{4}\widetilde{u}_\xi
  \in
  \calu\paren{\Ad(\rho_f)}
  \otimes_{R_f}\Q(\chi).
\end{equation}
Lemmas~\ref{lem:orientation-scalar} and
\ref{lem:fixed-adjoint-unit} give
\begin{equation}
\label{eq:orientation-regulator}
  \log_\ell\Reg_\mfp\paren{u_f^\opt}
  =
  c_{\chi,\mfp}
  \log_{\ell,\mfp_c}\paren{u_\xi}.
\end{equation}

\subsubsection{Relation to a higher Eisenstein element}

Let $E_{p + 1}$ be the level-one Eisenstein series of
weight $p + 1$. Darmon--Harris--Rotger--Venkatesh define
the higher Eisenstein element $\Sigma_1$ using the
reduction of $E_{p + 1}^{12}/\Delta^{p + 1}$
at the supersingular points; see
\cite[Definition~4.6 and Equation~93]{dhrv}.
Darmon--Harris--Rotger--Venkatesh \cite[Proposition~5.2]{dhrv},
with the $\mfp_c$-orientation
made explicit, gives
\[
  \paren{1 - \xi\paren{\overline{\mfl}}}
  \pair{\Sigma_1, \sbrac{\xi}_\mfp}
  =
  -\frac{1}{6}
  \log_{\ell,\mfp_c}\paren{u_{\xi,\lambda}}.
\]
Substituting Equation~\ref{eq:uxi} directly yields
\begin{equation}
\label{eq:sigma-xi}
  \pair{\Sigma_1, \sbrac{\xi}_\mfp}
  =
  -\frac{1}{6m(\xi)}
  \log_{\ell,\mfp_c}\paren{u_\xi}.
\end{equation}

\subsection{A chain of equalities}

Let $\Sigma_0$ be the Eisenstein element in
$\Z\sbrac{\Pic(B_p)}$.
Let $\paren{\Theta_\mfp^{\std}}^*$
be the adjoint of the standard theta lift
in Equation~\ref{eq:Theta-standard}.
Under the bijection of
Proposition~\ref{prop:oriented-maximal},
Darmon--Harris--Rotger--Venkatesh
\cite[Theorem~5.4]{dhrv} give
\[
  \paren{\Theta_\mfp^{\std}}^*
  \paren{(1 \otimes \log_\ell)\mfs_p}
  =
  \frac{1}{2}
  \paren{
    \Sigma_1 \otimes \Sigma_0
    +
    \Sigma_0 \otimes \Sigma_1}
  \pmod{\Sigma_0 \otimes \Sigma_0}.
\]
Moreover,
$\pair{\Sigma_0, \sbrac{\xi}_\mfp} = 0$ and
$\pair{\Sigma_0, \sbrac{\one}_\mfp} = h\paren{\calo_c}$.
Equation~\ref{eq:theta-opt-standard} therefore gives
\begin{align*}
  \log_\ell\mfs_p\paren{f^\opt(z, pz)}
  &= c_{\chi, \mfp}
    \pair{
      \sbrac{\one}_\mfp
      \otimes
      \sbrac{\xi}_\mfp,
      \paren{\Theta_\mfp^{\std}}^*
      \paren{(1 \otimes \log_\ell)\mfs_p}} \\
  &= \frac{1}{2}h\paren{\calo_c}
    c_{\chi, \mfp}
    \pair{\Sigma_1, \sbrac{\xi}_\mfp} \\
  &= -\frac{h\paren{\calo_c}}{12m(\xi)}
    c_{\chi, \mfp}
    \log_{\ell, \mfp_c}\paren{u_\xi} \\
  &= -\frac{\sbrac{H_c : K}}{12m(\xi)}
    \log_\ell\Reg_\mfp\paren{u_f^\opt}.
\end{align*}
The last equality follows from
Equation~\ref{eq:orientation-regulator}. This proves
Theorem~\ref{thm:opt}.

%%%%%%%%%%%%%%%%%%%%%%%%%%%%%%%%%%%%%%%%%%%%%%%%%%%%%%%%%%%%%%%%%

\section{Proof of Theorem~\ref{thm:hv}}
\label{sec:proof-hv}

Set
$r_\xi:=\sbrac{H_c:K}/(12m(\xi))$. By
Theorem~\ref{thm:opt}, for every auxiliary prime $p$
inert in $K$ and every prime $\mfp$ of $E$ above $p$,
\begin{equation}
\label{eq:hv-opt}
  \sbrac{\Gamma(1):\Gamma_0(N)}
  \calp_{\hv,\ell}\paren{f^\opt}
  =
  -r_\xi\log_\ell\Reg_\mfp\paren{u_f^\opt}.
\end{equation}

The optimal vector $\varphi^\opt$ and adjoint unit
$u_f^\opt$ are naturally defined over $\Q(\chi)$,
rather than a priori over the newform coefficient field
$L_f$. We therefore fix a finite extension
$L/L_f$ containing $\Q(\chi)$, with ring of integers
$\calo_L$, and put $\calo_{L,N}:=\calo_L\sbrac{1/N}$ and
\[
  \Pi_{f, \Sigma, \calo_{L,N}}
  :=
  \Pi_{f, \Sigma, \calo_{f,N}}
  \otimes_{\calo_{f,N}}\calo_{L,N},
  \qquad
  \calu_L\paren{\Ad(\rho_f)}
  :=
  \calu\paren{\Ad(\rho_f)}
  \otimes_{R_f}\calo_L.
\]
The coefficient-field trace in the following lemma
then produces a unit in the original $R_f$-lattice.
The lemma records that the comparison with the
newvector is uniform in the auxiliary data.

\begin{lemma}
\label{lem:uniform-comparison-transfer}
  Let
  $\varphi \in \Pi_{f, \Sigma, \calo_{L,N}}$ satisfy
  $\calp(\varphi) \neq 0$. Let $\sca$ be a collection
  of auxiliary data $(p, \mfp, \ell)$. Suppose that there
  are $m_\varphi \in \N$ and
  $u_\varphi \in \calu_L\paren{\Ad(\rho_f)}$,
  independent of $(p, \mfp, \ell)$, such that
  \[
    m_\varphi
    \calp_{\hv, \ell}\paren{f_\varphi}
    =
    \log_\ell\Reg_\mfp\paren{u_\varphi}
  \]
  for every $(p, \mfp, \ell) \in \sca$. Then there
  are $m_f \in \N$ and
  $u_f \in \calu\paren{\Ad(\rho_f)}$, independent
  of $(p, \mfp, \ell)$, such that
  \begin{equation}
  \label{eq:hv-new-uniform}
    m_f\calp_{\hv, \ell}\paren{f^\new}
    =
    \log_\ell\Reg_\mfp\paren{u_f}
  \end{equation}
  for every $(p, \mfp, \ell) \in \sca$.
\end{lemma}

\begin{proof}
  Set $D:=\prod_{q\mid N}(q-1)$ and write
  \[
    \frac{\calp(\varphi)}{\calp(f^\new)}
    =
    \frac{a_\varphi}{b_\varphi}
  \]
  with nonzero $a_\varphi,b_\varphi\in\calo_L$. Choose
  $a_\varphi^\vee\in\calo_L$ and $n_a\in\N$ such that
  $a_\varphi^\vee a_\varphi=n_a$. Equation~
  \ref{eq:comparison-cross-multiplied} gives
  \[
    b_\varphi D^2
    \calp_{\hv,\ell}\paren{f_\varphi}
    =
    a_\varphi D^2
    \calp_{\hv,\ell}\paren{f^\new}.
  \]
  Combining this with the assumed identity for
  $f_\varphi$ gives
  \[
    n_aD^2m_\varphi
    \calp_{\hv,\ell}\paren{f^\new}
    =
    \log_\ell\Reg_\mfp
    \paren{a_\varphi^\vee b_\varphi D^2u_\varphi}.
  \]

  Let $F_f := \Q\paren{\chi_{\Ad(\rho_f)}}$, and
  choose $c_f \in \N$ such that
  $c_f\calo_{F_f} \subset R_f$. Set
  \begin{align*}
    u_f
    &:=
    c_f\Tr_{L/F_f}
    \paren{
      a_\varphi^\vee b_\varphi D^2u_\varphi}, \\
    m_f
    &:=
    c_f\sbrac{L : F_f}n_aD^2m_\varphi.
  \end{align*}
  Compatibility of the regulator with scalar
  multiplication and coefficient-field trace gives the
  required identity. All quantities in these definitions
  depend only on $f$, $L$, and $\varphi$, not on the
  auxiliary datum. Hence $m_f$ and $u_f$ are uniform over
  $\sca$.
\end{proof}

By Proposition~\ref{prop:opt-pairing-nonzero},
$\calp(\varphi^\opt) \neq 0$. Choose
$c_L^\opt \in \N$ such that
$\varphi_L^\opt := c_L^\opt\varphi^\opt$ belongs to
$\Pi_{f, \Sigma, \calo_{L,N}}$, and choose
$d_L^\opt \in \N$ such that
\[
  u_{f, L}^\opt
  :=
  -d_L^\opt r_\xi u_f^\opt
  \in
  \calu_L\paren{\Ad(\rho_f)}.
\]
Multiplying Equation~\ref{eq:hv-opt} by
$c_L^\opt d_L^\opt$ gives
\[
  d_L^\opt
  \sbrac{\Gamma(1) : \Gamma_0(N)} \,
  \calp_{\hv, \ell}\paren{f_{\varphi_L^\opt}}
  =
  \log_\ell\Reg_\mfp
  \paren{c_L^\opt u_{f, L}^\opt}
\]
for every inert auxiliary prime $p$.

Apply Lemma~\ref{lem:uniform-comparison-transfer} to
$\varphi_L^\opt$, the fixed unit
$c_L^\opt u_{f, L}^\opt$, and the collection of inert
auxiliary data. We obtain $m_f > 0$ and a fixed
$u_f \in \calu\paren{\Ad(\rho_f)}$ such that
\[
  m_f\calp_{\hv, \ell}\paren{f^\new}
  =
  \log_\ell\Reg_\mfp\paren{u_f}
\]
for every inert auxiliary prime.

The formula for $u_f$ in
Lemma~\ref{lem:uniform-comparison-transfer} uses only
scalar multiplication and coefficient-field trace.
Since $u_f^\opt$ lies in the $G_\Q$-stable summand
$\Ind_{G_K}^{G_\Q}(\xi) \subset \Ad(\rho_f)$, so
does $u_f$. 
Set $u := \sbrac{\Gamma(1) : \Gamma_0(N)} u_f$.
For inert $p$,
Equation~\ref{eq:hv-new-uniform} and the finite-level
comparison in Section~\ref{sec:hv-period} give
\[
  m_f\log_\ell\mfs_p
  \paren{\Tr_p^{Np}\paren{f(z)f^*(pz)}}
  =
  \log_\ell\Reg_\mfp\paren{u}.
\]

It remains to consider an auxiliary prime $p$ that
splits in $K$. Darmon--Harris--Rotger--Venkatesh
\cite[Section~1.3]{dhrv} show that
\[
  \Tr_p^{Np}\paren{f(z)f^*(pz)} = 0.
\]
On the regulator side, $\Frob_\mfp \in \gal(E/K)$,
so $x_\mfp$ is diagonal in the induced basis. The
summand
\[
  \Ind_{G_K}^{G_\Q}(\xi)
  \subset
  \Ad(\rho_f)
\]
is the off-diagonal subspace spanned by $E_{12}$ and
$E_{21}$. Since $u$ belongs to this summand, evaluation
on the diagonal vector $x_\mfp$ gives $u(x_\mfp) = 0$
in the additive module
$\calo_E^\times \otimes_\Z R_f$. Equivalently, its
multiplicative representative is the trivial unit, and
hence $\log_\ell\Reg_\mfp\paren{u} = 0$.

The same pair $(m_f, u)$ therefore satisfies the
Harris--Venkatesh identity at both inert and split
auxiliary primes. Every prime $p \nmid N$ is
unramified in $K$, because $\disc(K)$ divides $N$.
Thus every auxiliary prime is either inert or split,
which proves Theorem~\ref{thm:hv}.

%%%%%%%%%%%%%%%%%%%%%%%%%%%% References %%%%%%%%%%%%%%%%%%%%%%%%%%%%%%

%\cleardoublepage
%\phantomsection
%\addcontentsline{toc}{part}{References}

\begingroup
	\renewcommand{\section}[2]{\part#1{#2}}
	%\nocite{*}
	\phantomsection
	\addcontentsline{toc}{part}{References}
	\bibliography{bibliography}{}

\begin{thebibliography}{DHRV22}

\bibitem[AL70]{atkin-lehner}
A.~O.~L. Atkin and J.~Lehner.
\newblock Hecke operators on {$\Gamma \sb{0}(m)$}.
\newblock {\em Math. Ann.}, 185:134--160, 1970.

\bibitem[Ata22]{atanasov}
Stanislav Atanasov.
\newblock {\em Derived {H}ecke Operators on Coherent Cohomology of Unitary
  {S}himura Varieties}.
\newblock PhD thesis, Columbia University, 2022.

\bibitem[Bum97]{bump}
Daniel Bump.
\newblock {\em Automorphic forms and representations}, volume~55 of {\em
  Cambridge Studies in Advanced Mathematics}.
\newblock Cambridge University Press, Cambridge, 1997.

\bibitem[Cas73]{casselman}
William Casselman.
\newblock On some results of {A}tkin and {L}ehner.
\newblock {\em Math. Ann.}, 201:301--314, 1973.

\bibitem[CS80]{casselman-shalika}
William Casselman and Joseph Shalika.
\newblock The unramified principal series of {$p$}-adic groups. {II}. {T}he
  {W}hittaker function.
\newblock {\em Compositio Math.}, 41(2):207--231, 1980.

\bibitem[Del71]{deligne}
Pierre Deligne.
\newblock Formes modulaires et repr\'{e}sentations {$l$}-adiques.
\newblock In {\em S\'{e}minaire {B}ourbaki. {V}ol. 1968/69: {E}xpos\'{e}s
  347--363}, volume 175 of {\em Lecture Notes in Math.}, pages Exp. No. 355,
  139--172. Springer, Berlin, 1971.

\bibitem[DHRV22]{dhrv}
Henri Darmon, Michael Harris, Victor Rotger, and Akshay Venkatesh.
\newblock The derived {H}ecke algebra for dihedral weight one forms.
\newblock {\em Michigan Math. J.}, 72:145--207, 2022.

\bibitem[DI95]{diamond-im}
Fred Diamond and John Im.
\newblock Modular forms and modular curves.
\newblock In {\em Seminar on {F}ermat's {L}ast {T}heorem ({T}oronto, {ON},
  1993--1994)}, volume~17 of {\em CMS Conf. Proc.}, pages 39--133. Amer. Math.
  Soc., Providence, RI, 1995.

\bibitem[DR73]{deligne-rapoport}
Pierre Deligne and Michael Rapoport.
\newblock Les sch\'{e}mas de modules de courbes elliptiques.
\newblock In {\em Modular functions of one variable, {II} ({P}roc. {I}nternat.
  {S}ummer {S}chool, {U}niv. {A}ntwerp, {A}ntwerp, 1972)}, Lecture Notes in
  Math., Vol. 349, pages 143--316. Springer, Berlin, 1973.

\bibitem[DS74]{deligne-serre}
Pierre Deligne and Jean-Pierre Serre.
\newblock Formes modulaires de poids {$1$}.
\newblock {\em Ann. Sci. \'{E}cole Norm. Sup. (4)}, 7:507--530 (1975), 1974.

\bibitem[Eic55]{eichler-1955}
Martin Eichler.
\newblock Zur {Z}ahlentheorie der {Q}uaternionen-{A}lgebren.
\newblock {\em J. Reine Angew. Math.}, 195:127--151 (1956), 1955.

\bibitem[Eic57]{eichler-1957}
Martin Eichler.
\newblock Eine {V}erallgemeinerung der {A}belschen {I}ntegrale.
\newblock {\em Math. Z.}, 67:267--298, 1957.

\bibitem[Eme02]{emerton}
Matthew Emerton.
\newblock Supersingular elliptic curves, theta series and weight two modular
  forms.
\newblock {\em J. Amer. Math. Soc.}, 15(3):671--714, 2002.

\bibitem[Gel75]{gelbart}
Stephen~S. Gelbart.
\newblock {\em Automorphic forms on ad\`ele groups}, volume No. 83 of {\em
  Annals of Mathematics Studies}.
\newblock Princeton University Press, Princeton, NJ; University of Tokyo Press,
  Tokyo, 1975.

\bibitem[Gro87]{gross-1987}
Benedict~H. Gross.
\newblock Heights and the special values of {$L$}-series.
\newblock In {\em Number theory ({M}ontreal, {Q}ue., 1985)}, volume~7 of {\em
  CMS Conf. Proc.}, pages 115--187. Amer. Math. Soc., Providence, RI, 1987.

\bibitem[GV18]{galatius-venkatesh}
S{\o}ren Galatius and Akshay Venkatesh.
\newblock Derived {G}alois deformation rings.
\newblock {\em Adv. Math.}, 327:470--623, 2018.

\bibitem[HK92]{harris-kudla-1992}
Michael Harris and Stephen~S. Kudla.
\newblock Arithmetic automorphic forms for the nonholomorphic discrete series
  of {${\rm GSp}(2)$}.
\newblock {\em Duke Math. J.}, 66(1):59--121, 1992.

\bibitem[HK04]{harris-kudla-2004}
Michael Harris and Stephen~S. Kudla.
\newblock On a conjecture of {J}acquet.
\newblock In {\em Contributions to automorphic forms, geometry, and number
  theory}, pages 355--371. Johns Hopkins Univ. Press, Baltimore, MD, 2004.

\bibitem[Hor23]{horawa}
Aleksander Horawa.
\newblock Motivic action on coherent cohomology of {H}ilbert modular varieties.
\newblock {\em Int. Math. Res. Not. IMRN}, (12):10439--10531, 2023.

\bibitem[HV19]{hv}
Michael Harris and Akshay Venkatesh.
\newblock Derived {H}ecke algebra for weight one forms.
\newblock {\em Exp. Math.}, 28(3):342--361, 2019.

\bibitem[Ich08]{ichino}
Atsushi Ichino.
\newblock Trilinear forms and the central values of triple product
  {$L$}-functions.
\newblock {\em Duke Math. J.}, 145(2):281--307, 2008.

\bibitem[Jac72]{jacquet}
Herv\'{e} Jacquet.
\newblock {\em Automorphic forms on {${\rm GL}(2)$}. {P}art {II}}.
\newblock Lecture Notes in Mathematics, Vol. 278. Springer-Verlag, Berlin-New
  York, 1972.

\bibitem[JL70]{jacquet-langlands}
Herv\'{e} Jacquet and Robert~P. Langlands.
\newblock {\em Automorphic forms on {${\rm GL}(2)$}}.
\newblock Lecture Notes in Mathematics, Vol. 114. Springer-Verlag, Berlin-New
  York, 1970.

\bibitem[Kat73]{katz-1973}
Nicholas~M. Katz.
\newblock {$p$}-adic properties of modular schemes and modular forms.
\newblock In {\em Modular functions of one variable, {III} ({P}roc. {I}nternat.
  {S}ummer {S}chool, {U}niv. {A}ntwerp, {A}ntwerp, 1972)}, Lecture Notes in
  Math., Vol. 350, pages 69--190. Springer, Berlin, 1973.

\bibitem[Kat78]{kato}
Shin-ichi Kato.
\newblock On an explicit formula for class one {W}hittaker functions on
  classical groups over {$p$}-adic fields, 1978.
\newblock preprint, University of Tokyo.

\bibitem[KL81]{kubert-lang}
Daniel~S. Kubert and Serge Lang.
\newblock {\em Modular units}, volume 244 of {\em Grundlehren der
  Mathematischen Wissenschaften}.
\newblock Springer-Verlag, New York-Berlin, 1981.

\bibitem[KM85]{katz-mazur}
Nicholas~M. Katz and Barry Mazur.
\newblock {\em Arithmetic moduli of elliptic curves}, volume 108 of {\em Annals
  of Mathematics Studies}.
\newblock Princeton University Press, Princeton, NJ, 1985.

\bibitem[Lec23]{lecouturier-hv}
Emmanuel Lecouturier.
\newblock On triple product {L}-functions and a conjecture of
  {H}arris-{V}enkatesh.
\newblock {\em Int. Math. Res. Not. IMRN}, (22):19476--19506, 2023.

\bibitem[Maz77]{mazur}
Barry Mazur.
\newblock Modular curves and the {E}isenstein ideal.
\newblock {\em Inst. Hautes \'{E}tudes Sci. Publ. Math.}, (47):33--186 (1978),
  1977.
\newblock With an appendix by Mazur and M. Rapoport.

\bibitem[Mer96]{merel2}
Lo\"{\i}c Merel.
\newblock L'accouplement de {W}eil entre le sous-groupe de {S}himura et le
  sous-groupe cuspidal de {$J_0(p)$}.
\newblock {\em J. Reine Angew. Math.}, 477:71--115, 1996.

\bibitem[Oh22]{oh}
Gyujin Oh.
\newblock {\em Arithmetic of Higher Coherent Cohomology of {S}himura
  Varieties}.
\newblock PhD thesis, Princeton University, 2022.

\bibitem[PV21]{prasanna-venkatesh}
Kartik Prasanna and Akshay Venkatesh.
\newblock Automorphic cohomology, motivic cohomology, and the adjoint
  {$L$}-function.
\newblock {\em Ast\'{e}risque}, (428):viii+132, 2021.

\bibitem[Rob73]{robert}
Gilles Robert.
\newblock {\em Unit\'es elliptiques et formules pour le nombre de classes des
  extensions ab\'eliennes d'un corps quadratique imaginaire}, volume~36 of {\em
  Bulletin de la Soci\'et\'e{} Math\'ematique de France. Supplement M\'emoire}.
\newblock Soci\'et\'e{} Math\'ematique de France, Paris, 1973.
\newblock Suppl\'ement au Bull. Soc. Math. France, Tome 101.

\bibitem[Ser79]{serre-local-fields}
Jean-Pierre Serre.
\newblock {\em Local fields}, volume~67 of {\em Graduate Texts in Mathematics}.
\newblock Springer-Verlag, New York-Berlin, 1979.
\newblock Translated from the French by Marvin Jay Greenberg.

\bibitem[Shi72]{shimizu}
Hideo Shimizu.
\newblock Theta series and automorphic forms on {${\rm GL}_{2}$}.
\newblock {\em J. Math. Soc. Japan}, 24:638--683, 1972.

\bibitem[Shi76]{shintani}
Takuro Shintani.
\newblock On an explicit formula for class-{$1$} ``{W}hittaker functions'' on
  {$GL\sb{n}$} over {$p$}-adic fields.
\newblock {\em Proc. Japan Acad.}, 52(4):180--182, 1976.

\bibitem[Shi94]{shimura}
Goro Shimura.
\newblock {\em Introduction to the arithmetic theory of automorphic functions},
  volume~11 of {\em Publications of the Mathematical Society of Japan}.
\newblock Princeton University Press, Princeton, NJ, 1994.
\newblock Reprint of the 1971 original, Kan\^{o} Memorial Lectures, 1.

\bibitem[Tat84]{tate-stark-1984}
John Tate.
\newblock {\em Les conjectures de {S}tark sur les fonctions {$L$} d'{A}rtin en
  {$s=0$}}, volume~47 of {\em Progress in Mathematics}.
\newblock Birkh\"{a}user Boston, Inc., Boston, MA, 1984.
\newblock Lecture notes edited by Dominique Bernardi and Norbert Schappacher.

\bibitem[Ven19]{venkatesh-2}
Akshay Venkatesh.
\newblock Derived {H}ecke algebra and cohomology of arithmetic groups.
\newblock {\em Forum Math. Pi}, 7:e7, 119, 2019.

\bibitem[Vig89]{vigneras-1989}
Marie-France Vign\'{e}ras.
\newblock Repr\'{e}sentations modulaires de {${\rm GL}(2,F)$} en
  caract\'{e}ristique {$l,\;F$} corps {$p$}-adique, {$p\neq l$}.
\newblock {\em Compositio Math.}, 72(1):33--66, 1989.

\bibitem[Voi21]{voight}
John Voight.
\newblock {\em Quaternion algebras}, volume 288 of {\em Graduate Texts in
  Mathematics}.
\newblock Springer, Cham, 2021.

\bibitem[Wal85]{waldspurger}
Jean-Loup Waldspurger.
\newblock Sur les valeurs de certaines fonctions {$L$} automorphes en leur
  centre de sym\'{e}trie.
\newblock {\em Compositio Math.}, 54(2):173--242, 1985.

\bibitem[Whi03]{whittaker}
E.~T. Whittaker.
\newblock An expression of certain known functions as generalized
  hypergeometric functions.
\newblock {\em Bull. Amer. Math. Soc.}, 10(3):125--134, 1903.

\bibitem[YZZ13]{yzz}
Xinyi Yuan, Shou-Wu Zhang, and Wei Zhang.
\newblock {\em The {G}ross-{Z}agier formula on {S}himura curves}, volume 184 of
  {\em Annals of Mathematics Studies}.
\newblock Princeton University Press, Princeton, NJ, 2013.

\bibitem[Zha01]{zhang-asian}
Shou-Wu Zhang.
\newblock Gross-{Z}agier formula for {${\rm GL}_2$}.
\newblock {\em Asian J. Math.}, 5(2):183--290, 2001.

\bibitem[Zha23a]{zhang-hvs}
Robin Zhang.
\newblock The {H}arris--{V}enkatesh conjecture for derived {H}ecke operators
  {II}: a unified {S}tark conjecture, 2023.
\newblock arXiv:2305.08956.

\bibitem[Zha23b]{zhang-rs}
Robin Zhang.
\newblock The {H}arris--{V}enkatesh conjecture for derived {H}ecke operators
  {III}: local constants, 2023.
\newblock \textit{Ann. Math. Blaise Pascal}, to appear. arXiv:2301.00612.

\bibitem[Zha26]{zhang-real}
Robin Zhang.
\newblock The {H}arris--{V}enkatesh conjecture for derived {H}ecke operators
  {IV}: real dihedral forms, 2026.
\newblock In progress.

\end{thebibliography}
	\bibliographystyle{alpha}
\endgroup

\end{document}